\documentclass[12pt,letterpaper]{lsuetd}

\usepackage{amsthm,amsmath,amssymb,amsxtra,latexsym,mathrsfs,url}
\usepackage{ab}
\usepackage{thesis_coms}
\usepackage{verbatim} 
\usepackage{cite} 

\usepackage{setspace} 
\usepackage[all]{nowidow} 

\usepackage{graphicx} 
\usepackage{tikz-cd} 
\usepackage{wrapfig} 

\usepackage{scalerel,stackengine}
\stackMath
\newcommand\verywidehat[1]{%
\savestack{\tmpbox}{\stretchto{%
  \scaleto{%
    \scalerel*[\widthof{\ensuremath{#1}}]{\kern-.6pt\bigwedge\kern-.6pt}%
    {\rule[-\textheight/2]{1ex}{\textheight}}
  }{\textheight}%
}{0.5ex}}%
\stackon[1pt]{#1}{\tmpbox}%
}

\usepackage{scalerel,stackengine,amsxtra} 
\stackMath %
\newcommand\verywidetilde[1]{ 
\savestack{\tmpbox}{\stretchto{ %
  \scaleto{ %
    \scalerel*[\widthof{\ensuremath{#1}}]{\kern-0.1pt \sptilde \kern-0.3pt} %
    {\rule[-\textheight/2]{1ex}{\textheight}} 
  }{\textheight} %
}{1.5ex}} %
\stackon[-2pt]{#1}{\tmpbox} %
}

\allowdisplaybreaks 

\usepackage{setspace,graphics,dsfont,verbatim,paralist,indentfirst}
\makeatletter

\begin{document}
\renewcommand\@pnumwidth{1.55em}
\renewcommand\@tocrmarg{9.55em}
\renewcommand*\l@chapter{\@dottedtocline{0}{1.5em}{2.3em}}
\renewcommand*\l@figure{\@dottedtocline{1}{0em}{3.1em}}
\let\l@table\l@figure

\pagenumbering{roman}
\thispagestyle{empty}
\begin{center}
TOWARDS THEORY AND APPLICATIONS OF GENERALIZED~CATEGORIES TO AREAS OF TYPE~THEORY AND CATEGORICAL~LOGIC\\

\vfill
\doublespacing
A Dissertation \\
\singlespacing
Submitted to the Graduate Faculty of the \\
Louisiana State University and \\
Agricultural and Mechanical College \\
in partial fulfillment of the \\
requirements for the degree of \\
Doctor of Philosophy \\
\doublespacing
in \\
                                       
The Department of Mathematics \\
\singlespacing
\vfill

by \\
Lucius Traylor Schoenbaum \\
B.S., University of Georgia, 2005 \\
M.A., University of Georgia, 2008 \\
M.S., Louisiana State University, 2012 \\
\vspace{1pt}
May 2017
\end{center}
\pagebreak



\addcontentsline{toc}{chapter}{\hspace{-1.5em} {{\bf \large Acknowledgements}} \vspace{12pt}}


\chapter*{Acknowledgments}
\doublespacing
\vspace{0.55ex}
\vspace{5ex}
This work clearly 
shows the scars of an existence on earth. %
I can only hope that for this reason, the work has something of the flavor of life, even though it is really only mathematics. %

I'm grateful to many people for helping me along the way. 
I would like to thank my Mom and Dad for their love and support. I love you Mom. I love you Dad.
I would like to thank Cecilia, Kelsey, Adrian, Mandy, and Melanie for their support while I was working on this, and experiencing all the highs and lows of life as a graduate student. I love you all. 

I would like to thank the mentors I have had over the years, especially Robert Varley, Robert Rumely, David Edwards, Ed Azoff, Akos Magyar, and Ted Shifrin at UGA mathematics for supporting and encouraging me and teaching me so much. I also thank everyone at UGA philosophy, especially my advisor for my Master's degree in philosophy, Brad Bassler, who first exposed me 
to the closely 
coupled vision of logic, philosophy, mathematics, and physics 
that has been the foundation of my adult life. 
I owe a debt of gratitude to Gestur Olafsson, Jimmie Lawson, and Daniel Sage for their mentorship and support during my time in Baton Rouge, and for often giving me the chance  
to pursue questions of my own. 
Several other LSU professors helped me in ways big and small, and I am especially grateful to James Madden, Milen Yakimov, Charles Delzell, and Ricardo Estrada for their help and their knowledge. I'm also grateful to the LSU math department staff, LSU math's talented IT, and the custodians and campus employees who kept my office clean and safe. 

There are many others whose work is reflected in mine. 
I am grateful to all the researchers who have inspired me and lifted me through their work. 
I would like to take the opportunity to thank all those who have been reviewers of my papers, especially the one based on my master's thesis, for their critiques  that led to many improvements. %
I am grateful to 
the people of Louisiana 
for hosting me, and who by way of a three year Board of Regents Fellowship made the early part of my studies possible. 
They (most of the time
) provided an environment conducive to study and reflection. 
For that, I would also thank all the places and spaces of South Louisiana I have haunted, if only I could.

\singlespacing
\tableofcontents
\pagebreak

\renewcommand\@pnumwidth{1.55em}
\renewcommand\@tocrmarg{8.55em}


\renewenvironment{abstract}{{\hspace{-2.2em} \huge \textbf{\abstractname}} \par}{\pagebreak}
\addcontentsline{toc}{chapter}{\hspace{-1.5em} {\bf \large Abstract}}
\begin{abstract}
\vspace{0.55ex}
\vspace{5ex}
\doublespacing


%
Motivated by potential applications to theoretical computer science, in particular those areas where the Curry-Howard correspondence plays an important role, as well as by the ongoing search in pure mathematics for feasible approaches to higher category theory, we undertake a detailed study of a new mathematical abstraction, the generalized category. 
It is a partially defined monoid equipped with endomorphism maps defining sources and targets on arbitrary elements, possibly allowing a proximal behavior with respect to composition. 
We first present a formal introduction to the theory of generalized categories. We describe functors, equivalences, natural transformations, adjoints, and limits in the generalized setting. 
Next we indicate how the theory of monads extends to generalized categories, and discuss applications to computer science. In particular we discuss implications for the functional programming paradigm, and discuss how to extend categorical semantics to the generalized setting. 
Next, we present a variant of the calculus of deductive systems developed in \cite{LaK1c,LaK2}, 
and give a generalization of the Curry-Howard-Lambek theorem giving an equivalence between the category of typed lambda-calculi and the category of cartesian closed categories and exponential-preserving morphisms that leverages the theory of generalized categories. 
Next, we develop elementary topos theory in the generalized setting of ideal toposes, utilizing the formalism developed for the Curry-Howard-Lambek theorem. In particular, we prove that ideal toposes possess the same Heyting algebra structure and squares of adjoints that ordinary toposes do. 
Finally, we develop generalized sheaves, and show that such categories form ideal toposes. We extend Lawvere and Tierney's theorem relating $j$-sheaves and sheaves in the sense of Grothendieck to the generalized setting. 

\end{abstract}

\pagenumbering{arabic}
\addtocontents{toc}{\vspace{12pt} \hspace{-1.7em} {\bf \large Chapter} \vspace{-0.2em}}

\singlespacing

\setlength{\textfloatsep}{12pt plus 2pt minus 2pt}
\setlength{\intextsep}{6pt plus 2pt minus 2pt}
\chapter{Prelude: The 2-Category of Categories}\label{c.prelude}
\doublespacing
\vspace{10ex}


The following is a brief review of perhaps the most important elementary construction in category theory: the strict 2-category of categories.

Let $\sC, \sD$ be categories. Two natural transformations $\beta: G \natto H, \alpha: F \natto G$ between functors $F,G: \sC \to \sD$ may be composed via the rule
$$\beta \vertof \alpha (X) := \beta(X) \x \alpha(X)$$
where $(\x)$ denotes composition in $\sD$. This gives a category $\Nat(\sC, \sD)$. Identities in $\Nat(\sC, \sD)$ are given by $\id_F (X) := \id_X$. 

Given natural transformations $\alpha: F \natto G$ between functors $\sC \to \sD$, and $\beta: F' \natto G'$ between functors $\sD \to \sE$, we obtain a well-defined function $\Ob(\sC) \to \Mor(\sE)$ via
$$\beta \star \alpha (X) := \alpha ( \hat \beta (X) ) \x \bar \alpha ( \beta (X)),$$
where hats and bars are used as defined in section \ref{s.gencat} below. This can also be written
$$\beta \star \alpha = (\alpha \of \hat\beta) \vertof (\bar\alpha \of \beta)$$
Note that 
$$\bar \alpha (X) = \overline{\alpha(X)},$$
$$\hat \alpha (X) = \widehat{\alpha(X)}.$$

\prop[The Five Facts]
In the notation above, whenever expressions on both sides of the formula are defined, we have:
\enu
	\item $\beta \star \alpha = (\hat\alpha \of \beta) \vertof (\alpha \of \bar\beta).$ 
	\item $\beta \star \alpha$ is a natural transformation $G \of F \natto G' \of F'.$ 
	\item $(\gamma \star \beta) \star \alpha = \gamma \star (\beta \star \alpha).$
	\item If
		$$\begin{rightbrace}
			\alpha: F \natto G \\
			\beta: G \natto H 
		\end{rightbrace} : \sC \to \sD, $$
		$$\begin{rightbrace}
			\alpha' : F' \natto G' \\
			\beta' : G' \natto H' 
		\end{rightbrace} : \sD \to \sE, $$
		then
		$$(\beta' \vertof \alpha') \star (\beta \vertof \alpha) = (\beta' \star \beta) \vertof (\alpha' \star \alpha).$$
	\item If $\id_F^{\vertof}$ is the identity of $F$ with respect to the product $\vertof$ in $\Nat(\sC, \sD)$, then 
		$$\alpha \star \id^{\vertof}_F = \alpha,$$
		$$\id^{\vertof}_F \star \beta = \beta,$$
		whenever both sides are defined.
\Enu
\Prop
\prf
(1) 
\begin{align*}
	(\beta \star \alpha) (X)		&= (\alpha \of \hat \beta) \vertof (\bar \alpha \of \beta) (X) \\
						&= \alpha(\hat \beta(X)) \x \bar \alpha (\beta(X)) \\
						&= \alpha(\widehat{\beta(X)}) \x F(\beta(X)) \\
						&= G(\beta(X)) \x \alpha(\overline{\beta(X)}) \\
						&= \hat \alpha (\beta(X)) \x \alpha(\bar \beta (X)) \\
						&= (\hat \alpha \of \beta) \vertof (\alpha \of \bar \beta) (X).
\end{align*}

(2) by Fact 1.

(3) Apply the definition.

(4) by Fact 1 and since $\hat \alpha = \bar \beta, \widehat{\alpha'} = \overline{\beta'}$.

(5) direct calculation.
\Prf

{\em Remarks.} 
\enu
	\item We may write simply $\id_F$ or $1_F$ in light of Fact (5).
	\item Fact 4 is often referred to as the {\em interchange law}.
\Enu

An immediate consequence of the Five Facts is the following: {\em The category of categories is a strict two-category.} %
By ``the category of categories'' is meant the set of small categories, functors, and natural transformations in a fixed universe $\sU_{univ}$. 

$$\scriptscriptstyle{\blacksquare \quad \blacksquare \quad \blacksquare}$$
\vspace{-5pt}

\pagebreak

\singlespacing
\chapter{Introduction}
\doublespacing
\vspace{10ex}



\section{Overview and Motivation}\label{s.overviewandmotivation}

Category theory \cite{MacCW,KaSc1,BaWe1} has its origins in mathematics, and has since become a well-established area of foundations, with a rich interaction with computer science. It begins with the insight that diagrams and morphisms have a mathematical life unto themselves, independent of function theory, and independent of any use of points as arguments. %
The more one works with categories, in fact, %
the more one becomes cognizant of the view that a category is in fact a set of morphisms that are in some way algebraically structured. %
%
In  \cite{EiMa1}, the paper on natural transformations in which the elementary notions of category theory are introduced for the first time, Eilenberg and MacLane write 
\begin{quote}
\onehalfspacing 
It is thus clear that the objects play a secondary role, and could be entirely omitted from the definition of a category. However, the manipulation of the applications would be slightly less convenient were this done. 
\end{quote}
Thus two views have been known to category theorists since the beginning of the subject. 
The two approaches, the one-sorted definition describing a universe of pure maps, and the two-sorted definition including the objects that are in applications prior to the maps that they inspire, pull against one another in a way that seems, in practice, like a natural, irresolvable tension. %
The latter approach has proven to be the dominant one, while the former approach has made occasional appearances, for example in work by Ehresmann \cite{Ehresmann1}, Street \cite{StT1}, and 
more recently, in work by R. 
Cockett \cite{CockettConstellations}. %

The potential for generalization begins with the less-often-used one-sorted formulation, which necessitates an axiom requiring the source and target maps $s$ and $t$ to be trivial upon iteration: $ss = st = s, tt = ts = t$. This condition, however, is extraneous. Dropping it gives rise to a rather general notion, 
which may be %
weakened further via replacing some equalities with inequalities, as suggested by some kinds of applications \cite{SmPl1}. 
This is the jumping-off-point of our work in Chapter \ref{c.gencat}. 

There exists no literature on this generalization of category theory. 
This lack of attention to the abstraction is likely due in some part to a lack of knowledge about the robustness of the categorical theory existing in the general case. %
Having noticed signs that such a theory might be sufficiently strong to have potential applications, the author undertook the investigations that appear in this dissertation. %
With much literature in computer science devoted to subjects such as metaprogramming, dependent types, and other generalizations of type theory based on Church's typed lambda calculus (see for example \cite{Barendregt1}), there exists no lack of potential applications of this work to categorical semantics and other areas of theoretical computer science. %
However, the author, not having been trained in computer science himself, does not venture far outside of the mathematical parts of category theory in this work. %
We find in this theoretical investigation that in the generalized setting a considerable part of the fundamental theory of categories and categorical logic persists in a remarkably robust form. During our investigation a wealth of unforeseen new abstractions have arisen, and there is much at the time of this writing that is still unexplored, both in the realm of theory and in the realm of applications. 

\section{Summary of Contents}\label{s.outline}

\indent
In Chapter \ref{c.gencat} we introduce the main abstraction of our work, the notion of generalized category. We show that the standard tools of category theory carry over to the generalized setting, and we also present some negative results that may perhaps indicate to the interested reader some of the obstacles to further 
generalization along the same lines. 

In Chapter \ref{c.genm} we give a treatment of the theory of monads. We show that the theory of the Kleisli triple carries over to the generalized setting, however only by a more intricate construction than in the ordinary categorical setting. We also investigate the theory of algebras and the Tripleability theorems. Although the theory of $T$-algebras for monads seems to fit comfortably in the setting of generalized categories, the results of an investigation into a theory of generalized Eilenberg-Moore category were disappointing. However, the Tripleability theorems and other aspects of the theory of $T$-algebras can nevertheless be carried out via the usual one-categorical construction. 

In Chapter \ref{c.th} we treat cartesian closed categories and the work of Lambek \cite{LaK1c,LaK2} on the Curry-Howard Correspondence. The main result shown is that the correspondence carries over to the generalized setting. This involves developing a new syntactic abstraction, a generalized typed lambda calculus, to coincide with the generalized semantics of Chapter \ref{c.gencat}. This result lays the groundwork for applications of generalized categories in type theory, particularly for models of subtyping and higher-kinded type systems. As we also discuss, the generalized result also proves instructive when interpreting the usual Curry-Howard Correspondence in the one-categorical setting. 

In Chapter \ref{c.it} we use the results of Chapter \ref{c.th} develop the theory of elementary toposes in the setting of generalized categories, following \cite{sga4,MaMo1,JoN1,FrD1,LaSc1}. The abstraction studied is referred to as an ideal topos, as it extends the notion of ideal cartesian closed category introduced in Chapter \ref{c.th}. 
The theory of toposes has an important (perhaps decisively so) relationship with the theory of sheaves. We have therefore undertaken an investigation into a generalization of sheaves and show that they, too, possess a robust theory, one which is in certain respects different than both the theory of ideal elementary toposes and the theory of one-categorical sheaves that it extends. 
The flavor of the subject is a departure from previous chapters due to the obstacles that arise to setting up the full toolkit of sheaf theory in the generalized setting. This impasse is overcome using bipartite (generalized) categories, or what we call pointed profunctors. 
As a byproduct of our work, we obtain a working definition of the category of generalized sets. %
While we discuss many fundamental topics in topos theory (basic topos theory, sheaf theory, and a few remarks on 2-sheaves) we do not include a discussion of the internal language of generalized categories, as this would burden the reader with still more lengthy preliminaries along the lines of the treatment in Chapter \ref{c.th}. 

The work done in Chapter \ref{c.it} shows that the topos theoretic foundations for logic, topology, and geometry laid down by some of the greatest minds of the past one hundred years has a wide new realm of applicability. 
Our work in this Chapter is lacking, however, for it states theory only, and gives no applications. 
We conclude this introduction by voicing the hope that such applications might someday be forthcoming. %

\pagebreak

\singlespacing
\chapter{Generalized Categories}\label{c.gencat}
\doublespacing
\vspace{10ex}

\section{Beginnings}\label{s.gencat}


{\em Preliminaries.} 
We use notation $\source(f),\target(f)$, $\dom(f),\cod(f)$, and $\bar f, \hat f$, more or less interchangeably, to denote the source and target of an element of a generalized category. The lattermost notation may be used when it improves readability of formulas. %
We write composition $G \of F := (f \mapsto G(F(f)))$ and in general, for mappings $F$ and $G$ with common domain and codomain (in which concatenation is meaningful) we define the operation
$$G \vertof F := (f \mapsto G(f)F(f)),$$
the standard vertical composition operation \cite{MacCW}. %
In any context where it is meaningful, we use the standard arrow notation $f:a \to b$ to mean that an element $f$ is given, the source of $f$ is $a$, and the target of $f$ is $b$. %
The notation $\downarrow$ 
indicates that all composed pairs of elements in the expression or relation are in fact composable pairs. 

\subsection{Definition}\label{ss.gencat}

\dfn\label{d.gencat}
A {\em generalized category} is a structure $(\sC, \dleq, \source,\target, \cdot)$ where $\sC$ is a set, $\dleq$ is a relation on $\sC$, $\source$ and $\target$ are mappings $\sC \to \sC$, and $(\cdot)$ is a partially defined mapping $\sC \times \sC \to \sC$, denoted $a \cdot b$ or $ab$. These are required to satisfy
\enu
	\item $(\sC, \dleq)$ is a partially ordered set, 
			\label{ax.gencat-po}
	\item $ab$ $\downarrow$ if and only if $\source(a) \dleq \target(b)$. \label{ax.gencat-po-comp}
	\item If $(ab)c$ $\downarrow$ or $a(bc)$ $\downarrow$ then $(ab)c = a(bc)$. \label{ax.gencat-assoc}
	\item If $ab$ $\downarrow$ then $\source(ab) = \source(b)$ and $\target(ab) = \target(a)$. \label{ax.gencat-comp-st}
	\item (Element-Identity) For all $a \in \sC$, there exists $b \in \sC$ such that \label{ax.gencat-element-id}
		\enu
			\item $\source(b) = \target(b) = a$,
			\item if $bc$ $\downarrow$ then $bc = c$,
			\item if $cb$ $\downarrow$ then $cb = c$,
		\Enu
	\item (Object-Identity) Let $a \in \sC$ and $\source(a) = \target(a) = a$. Then \label{ax.gencat-object-id}
		\enu
			\item if $ba$ $\downarrow$ then $ba = b$. 
			\item If $ab$ $\downarrow$ then $ab = b$.
		\Enu
	\item (Order Congruences\footnote{These axioms are needed for the Kleisli construction in Chapter \ref{c.genm}.})
		\enu
			\item If $a \dleq b$ then $\source(a) \dleq \source(b)$ and $\target(a) \dleq \target(b)$. \label{ax.gencat-order1}
			\item $a \dleq b$ and $c \dleq d$ and $ac,bd$ $\downarrow$ implies $ac \dleq bd.$ \label{ax.gencat-order2}
			\item $a \dleq b$ implies $1_a \dleq 1_b$. \label{ax.order3}
		\Enu
\Enu
The element $c$ of axiom (\ref{ax.gencat-element-id}) %
is unique, and is denoted $1_a$ or $\id_a$, and called the {\em identity} on $a$. %
\Dfn

As a partially ordered set a generalized category resembles, but is weaker than, a domain \cite{GrZ+}, indeed motivation for the ordering comes from domain theory \cite{SmPl1, WaD1}. If $a \dleq b$, we say that $a$ {\em approximates} $b$, and $b$ {\em sharpens} $a$. 
When the ordering $\dleq$ is nontrivial, one may call $\sC$ a {\em proximal} generalized category. %
We often think of proximal categories as having at least a bottom element $\bottom$, but we do not assume this in the definition, since we would like, as a special case, for an ordinary one-category to be a generalized category. %
If the order given by $\dleq$ is discrete, we might say that the generalized category is {\em discrete}, and similarly for other order-theoretic attributes, but as this may lead to confusion with the notion of a discrete category (one with essentially no morphisms), we shall say instead that such a generalized category is a {\em sharp} generalized category. We allow ourselves to refer to a {\em proximal} generalized category whenever we wish to emphasize that we refer to a generalized category that is not assumed to be sharp. 

%
An {\em element} $f \in \sC$ is an element $f$ of the underlying set $\sC$. 
An {\em object} $a$ in $\sC$ is an element $a$ of $\sC$ such that $\source(a) = \target(a) = a$. 
We write $\Ob(\sC)$ for the set of objects. 
For $a \in \sC$, we define the {\em height} of $a$, denoted $\height(a),$
to be the maximum of the set of nonnegative integers $n$ such that there exists a sequence $\vec s$ of source and target operations of length $n$ such that $\vec s (i)$ is an object, %
unless there is an infinite sequence $\vec s$ of source and target operations such that no subsequence yields an object. In that case, we say that $\height(a) = \infty$. 

With this terminology, Definition \ref{d.gencat} says that in a generalized category with identities, every element $a$ has an identity $1_a$, and that if the element is an object, this identity is $a$ itself. %
If $a \in \sC$ has identity $1_a$ and is not an object, then $a \neq 1_a$. %

The maps $\source$ and $\target$ of the definition are called the  {\em source} or {\em domain} and {\em target} or {\em codomain} maps, respectively. We may sometimes denote the map $\source(a)$ by $\bar a$, and the map $\target(a)$ by $\hat a$. 

Given a generalized category $\sC$, any element of $\sC$ may be composed with other compatible elements, and it is equipped with a ``tail'' of fellow elements, defined by the $\source$ and $\target$ maps. We think of the product 
as developing from right to left, and we may write $c:a \to b$
when $\source(a) = b$, $\target(a) = c$. %
Note as an aside that if one pictures instead a representation $a = \sous{c}a_b$ of $a$, one has a picture of composition $\sous{c}a_b \,\sous{b}d_e = \sous{c}(ad)_e$. 
This notation can be iterated to 
$$a = \sous{\sous{g}c_f}a_{\sous{e}b_d}$$ 
In this manner one can visualize a binary tree. 

\subsection{An Alternative Approach}\label{ss.secondapproach}

We pause to make note of an alternative approach, and discuss why we choose the approach of Definition \ref{d.gencat}. %

\dfn\label{d.gencatalt}
A {\em generalized category} is a structure $(\sC, \dleq, \source,\target, \cdot)$ where $\sC$ is a set, $\dleq$ is a relation on $\sC$, $\source$ and $\target$ are operators (mappings) $\sC \to \sC$, and $(\cdot)$ is a partially defined binary operation $\sC \times \sC \to \sC$, denoted $a \cdot b$ or $ab$. These are required to satisfy
\enu
	\item $(\sC, \dleq)$ is a partially ordered set, 
	\item If $(ab)c$ $\downarrow$ or $a(bc)$ $\downarrow$ then $(ab)c = a(bc)$. 
	\item If $ab$ $\downarrow$ then $\source(ab) = \source(b)$ and $\target(ab) = \target(a)$.
	\item $ab$ $\downarrow$ if and only if $\source(a) \dleq \target(b)$.
	\item (Object-Identity) Let $a \in \sC$ and $\source(a) = \target(a) = a$. Then \label{ax.gencatalt-objid}
		\enu
			\item if $ba$ $\downarrow$ then $ba = b$. 
			\item If $ab$ $\downarrow$ then $ab = b$.
		\Enu
	\item (Order Congruences)
		\enu
			\item If $a \dleq b$ then $\source(a) \dleq \source(b)$ and $\target(a) \dleq \target(b)$. \label{ax.gencatalt-order1}
			\item $a \dleq b$ and $c \dleq d$ and $ac,bd$ $\downarrow$ implies $ac \dleq bd.$ \label{ax.gencatalt-order2}
			\item $a \dleq b$ implies $1_a \dleq 1_b$. \label{ax.gencatalt-order3}
		\Enu
\Enu
A generalized category is said to be equipped {\em with identities} if for every $a \in \sC$, if there exists $b \in \sC$ such that $\source(b) = a$ or $\target(b) = a$, then there exists $c \in \sC$ such that $c b$ $\downarrow$ implies $c b = b$, and $b c$ $\downarrow$ implies $b c = b$.
The element $c$ 
is unique, and is denoted $1_a$ or $\id_a$, and called the {\em identity} on $a$. %
An {\em element} $f \in \sC$ is an element $f$ of the underlying set $\sC$. %
An {\em object} $a$ in $\sC$ is an element $a$ of $\sC$ such that $\source(a) = \target(a) = a$. %
A {\em subject} $U$ in $\sC$ is an element $U$ of $\sC$ such that there exists $f \in \sC$ such that $\source(f) = U$, or there exists $f \in \sC$ such that $\target(f) = U$. %
\Dfn

The approach of Definition \ref{d.gencat} has the advantage of having fewer basic concepts than Definition \ref{d.gencatalt}. %
All elements are subjects and all elements have identities. %
This makes many steps of the development go smoothly. %
On the other hand, Definition \ref{d.gencat} creates so many identities that one sometimes wonders if they are better avoided after all. %
Thus one might seem to be at an impasse concerning whether Definition \ref{d.gencat} or Definition \ref{d.gencatalt} is more preferable. %
This ambivalence is resolved by the notion of ideal category introduced in Chapter \ref{c.th}. %
Ideal categories arise naturally in categorical logic. %
In such categories, and in particular in the generalized category of contexts $\boC\Lambda$, there are identities present %
just as Definition \ref{d.gencat} requires. %
This tipping of the scales is the reason why we favor Definition \ref{d.gencat} over Definition \ref{d.gencatalt}. %
In order to facilitate discussions about generalized categories in the sense of Definition \ref{d.gencat}, we say that {\em closing over $1_{()}$} is the obvious operation of ensuring (via free generation where needed) that axiom (\ref{ax.gencat-element-id}) is satisfied. 



\subsection{Resuming, from Definition \ref{d.gencat}}

\prop\label{p.duality}
Up to reversal of $\,\dleq$, Definition \ref{d.gencat} is symmetric in the source and target maps $\source$ and $\target$. Therefore every proof $\Phi$ about a generalized category $\sC$ continues to hold when, in all assumptions, definitions, and deduction steps, composition, the order $\dleq$, and the role of source and target are reversed. 
\Prop

Such a proof $\Phi'$ is said to be obtained from $\Phi$ ``by duality'' \cite{MacCW}. %
This simple fact has a profound effect on the entire subject. %
The generalized category formed by the operation of Proposition \ref{p.duality} is called the {\em opposite generalized category} $\sC^{op}$ of $\sC$. %

\begin{example}
Let $\sC$ be a category \cite{MacCW}. %
Then {\em the generalized category generated by $\sC$} is obtained from $\sC$ by identifying the identity $1_X$ of each object $X \in \sC$ with $X$, and closing over $1_{()}$. %
Considering a concrete example, such as the generalized category generated by the category of all groups, we may write $\id_X$ for $X$, with the identification $\id_X = X$ being understood. %
More formally, we define: %

\dfn\label{d.cat}
A generalized category $\sC$ is a {\em category} or {\em one-category} if the source and target of every nonidentity $f$ in $\sC$ is an object in $\sC$. 
\Dfn

We now have a rough ontology:

$$
\begin{tabular}{c|c}
\begin{tabular}{c}
sharp category \\
= category 
\end{tabular} 
	& 
\begin{tabular}{c}
proximal category
\end{tabular} 
			\\ \hline
\begin{tabular}{c}
sharp generalized category 
\end{tabular} 
	& 
\begin{tabular}{c}
proximal generalized category \\
= generalized category \\
\end{tabular}
\end{tabular}
$$
\end{example}

\begin{example}
In some instances it is possible to write down a generalized category explicitly. There is an empty generalized category, and $\sC = \set{a:a \to a}$, the trivial generalized category. %
More generally, any set $S$ is a generalized category after setting $\source(a) = \target(a) = a$ for $a \in S$, we say that the generalized category is {\em discrete} or a {\em zero-category}, or simply that it is a set. %
(Thus, sets and categories are examples of generalized categories.) 
Because of the identity axiom, other than finite sets there are no finite generalized categories. 
To amend language, we therefore define:

\dfn\label{d.figen}
A generalized category $\sC$ is {\em finitely generated} if there is a finite set $\sC$ such that the remainder of $\sC$ consists only of identities. 
\Dfn

There are many examples of generalized categories that are not ordinary categories, the simplest perhaps being
$\sC = \set{a: a \to a, b: a \to b}$.
Another simple example is 
$\sC = \set{a: b \to b, b: a \to a}$.
This generalized category is finite, but does not possess objects, moreover every element is a subject. 
A generalized category may also lack objects due to infinite descent, for example
$\sC = \set{a_n : a_{n-1} \to a_{n-1} \mid n \in \Z}.$ 
\end{example}

\begin{example}
Besides the aforementioned sources in domain theory, 
motivation for the proximal relation $\dleq$ in a proximal generalized category comes, via categorical semantics, from %
the subtyping relation in some type theoretical systems \cite{Reynolds1980,Pierce1}, a feature characteristically found in object-oriented languages. %
Subtyped type systems are often preorders, thus, we can access the semantics given by a generalized category by, for one thing, equating mutual subtypes. %
As in domain theoretic orders, subtyped type systems often include a bottom type; %
they may also include a global maximum type. %
Such a structuring of types creates a comfortable intuitive environment for type theory, and makes the type checker behave less rigidly. %
However, representation of data in such systems can demand trade-offs that make such systems less suitable for some kinds of applications. %
Moreover, in an industry-level type system, problems and subtleties may arise due to the need for subtyping rules to %
interact coherently with rules that govern other advantageous type features, such as records, recursive types, and polymorphism. %
In practice, therefore, subtyping produces both benefits as well as costs, 
and has been the focus of much research and discussion in computer science. %

One approach to implementing subtyping involves data type {\em coercion,} or the automated physical modification of stored data at run-time. %
Type-theoretically, condition (\ref{ax.gencat-comp-st}) of Definition \ref{d.gencat} corresponds to a type system in which there exists a coercive evaluation mechanism. %
\end{example}

\begin{example}\label{ex.binarytree}
Let $\sC$ be a generalized category, and consider 
the condition on $\sC$ that hom sets should contain a unique element or else be empty. 
To obtain a (possibly infinite) planar binary tree one adds the condition that source and target may not loop except trivially, that is, for every element $a \in \sC$, and for every finite sequence $(x_1, \dots, x_n)$ where $x_i$ is either $\source$ or $\target$ (source or target) if $x_n x_{n-1} \dots x_1 a = a$ then it is required that $\source a = \target a = a$, that is, or (using the terminology of trees) that $a$ is a leaf.  
Presheaves on such trees arise for example in database theory, see for example \cite{SpK3}. %
\end{example}

\exa
A {\em generalized (directed) graph} (cf. Chapter \ref{c.th}) is simply a triple $(\sA, \source, \target)$, where $\sA$ is a carrier set, %
and $\source, \target$ are maps $\sA \to \sA$. %
An element of $\sA$ is (synonymously) an {\em edge}. %
An {\em object} in a generalized graph is an element $a \in \sA$ such that $\source a = \target a = a$, %
that is, a common fixed point of the endomorphisms $\source$ and $\target$. 
Ordinary graphs correspond bijectively with 1-dimensional generalized graphs, 
where we say that generalized graph is {\em 1-dimensional} if $\source \source = \source \text{ and } \target \target = \target.$
With the obvious composition via compound paths, a generalized graph becomes a (sharp) generalized category. 

There are plentiful settings where generalized graphs may arise. For example, suppose that there is a system of goods $\sA_0$. %
The edges of $\sA$ are certificates (issued, say perhaps, by different governing bodies) that say that a good $a \in \sA_0$ may be exchanged for another good $b \in \sA_0$. %
Suppose it is accepted that a good is always exchangeable for itself. %
Now let's suppose that such certificates themselves may be exchanged, 
but that this requires that one has a higher-level certificate for this higher-level trade. %
If we imagine a certain impetus exists among those we imagine making the exchanges, we can expect that there will next arise trading for these certificates as well, 
giving rise to a generalized graph (in fact, a generalized deductive system, via a simple extension of Kolmogorov's reasoning about intuitionistic logic in \cite{Kolmogorov1}). 
\Exa

\exa\label{x.setf}
For a planar binary tree $\ft$, let
\begin{align*}
\troot(\ft) &\text{ is the root of $\ft$.} \\
\tleft(\ft) &\text{ is the tree given by the left descendant of the root, and its descendants.} \\
\tright(\ft) &\text{ is the tree given by the right descendant of the root, and its descendants.}
\end{align*}
From any category $\sC$ we can form a sharp generalized category $\sC f$ as follows: take the set $\sC f$ to be the set of all planar binary trees of morphisms in $\sC$, subject to a source-and-target condition 
$$\dom \troot(\dom \ff) = \dom \troot(\ff) ,$$
and
$$\cod \troot(\cod \ff) = \cod \troot(\ff),$$
where if $\ff$ be such a tree,
$$\cod \ff = \tleft(\ff),$$
the left descendent tree of $\ff$, and 
$$\dom \ff = \tright(\ff),$$
the right descendent tree of $\ff$. 
These conditions set up a recursive condition on elements of $\sC f$.
For $\fg, \ff \in \sC f$, we define $\fg \cdot \ff$ to be the tree $\fh$ with left descendent $\troot(\ff)$, right descendent $\troot(\fg)$ and root $\troot(g)\cdot \troot(f)$. 
This is a well-defined product, by the source-and-target condition above. It is checked that this is a (sharp) generalized category. An element of $\sC f$ may be visualized as
$$
\begin{tikzcd}
															& \bullet \ar[d, dashed, "\hat \ff"]					 		\\
	X \ar[d, dashed, "\bar \ff"] \ar[r, dashed, "\ff"] 							& Y 		\\
	\bullet															&
\end{tikzcd}
$$
%
Constructions on the original $\sC$ can be carried over to $\sC f$, for example, if $\sC$ has products (equalizers, coproducts, coequalizers), then so (respectively) does $\sC f$. If $\sC$ is (co)complete, however, it does not imply that $\sC f$ is (co)complete. 

\Exa

\section{Elementary Theory, Category of Invertibles}\label{s.elem}

We now define functors and hom sets: 

\dfn\label{d.functor}
A mapping $\sC \to \sC'$ between generalized categories is {\em functorial} or a {\em functor} if 
\enu
	\item $a \dleq b \eimplies F(a) \dleq F(b),$
	\item $F(\bar a) = \overline{F(a)},$
	\item $F(\hat a) = \widehat{F(a)},$
	\item $F(ab) = F(a)F(b), \text{ if } ab \,\downarrow,$
	\item $F(1_a) = 1_{F(a)}.$
\Enu
We thus have a category ${\bf GenCat}$ of generalized categories and functors. %
\Dfn

Functors are also called {\em covariant functors}. A {\em contravariant functor} from $\sC$ to $\sC'$ is a unital map satisfying
\enu
	\item if $a \dleq b$ then $F(b) \dleq F(a)$,
	\item $F(\bar a) = \widehat{F(a)},$
	\item $F(\hat a) = \overline{F(a)},$
	\item $F(ab) = F(b)F(a) \text{ if } ab \,\downarrow,$
\Enu
instead of the corresponding covariant relations. %

\dfn\label{d.homsets}
The sets 
$$\hom(a,b) = \set{c \in \sC \mid \bar c = a, \hat c = b},$$ 
for $a,b \in \sC$, are called the {\em hom sets} of $\sC$. 
\Dfn

\dfn
A {\em subcategory} of a generalized category $\sC$ is a subset $\sC'$ of $\sC$ whose order is inherited from $\sC$ closed under source, target, composition, and identities: if $a \in \sC'$, then $1_a \in \sC'$. 
A subcategory $\sC'$ is {\em full} if $a,b \in \sC'$ implies $\hom(a,b)$ is contained in $\sC'$. 
\Dfn

The composition of two functors is a functor, and functors send objects to objects. 

\dfn
Two generalized categories $\sC$ and $\sC'$ are {\em isomorphic} if there is an invertible functor (i.e., invertible as a mapping) $F$ from $\sC$ to $\sC'$. 
\Dfn

\prop\label{p.flattening}
There is a functor, flattening, from the category of generalized categories to the category of categories. 
\Prop
\prf
Let $\sC$ be a generalized category with identities. 
Let $\Ob(\sC_{flat})$ be $\set{[f] \mid f \in \sC}$, the objects of $\sC$ indexed by the elements of $\sC$. 
Let $\Mor(\sC_F)$ again be a set $\set{(f) \mid f \in \sC}$ indexed by the elements of $\sC$, and define source and target 
$s((f)) = [s(f)],$ 
$t((f)) = [t(f)].$ 
Then $\sC_{flat}$ is a category whose composition and identities are
$$(g) \cdot (f) := (gf),$$
$$
1_{[f]} = (1_f).
$$
Given a functor $F: \sC \to \sD$ in ${\bf GenCat}$, we immediately obtain a functor $\sC_{flat} \to \sD_{flat}$. 
\Prf

Note that $\sC_{flat}$ contains a flattening of the identity structure, even in cases where $\hom(a,a) = \set{1_a}$.  

There is also a category $flat \sC$, the further flattening of $\sC$ to a zero-category. %
It is defined by:
$$flat(f) := \begin{cases} 	(f), 	& \text{if $f = 1_g$ for some $g \in \sC$, } \\ 
					[f]	& \text{otherwise,}
	\end{cases}
$$
where $[f]$ is defined by $\source([f]) = \target([f]) = [f]$, and $(f): flat(\source(f)) \to flat(\target(f)).$





\dfn\label{d.invertible}
If $\sC$ is a generalized category, an element $a \in \sC$ is {\em invertible} if there exists $b \in \sC$ such that $ab = 1_{\hat a}$ and $ba = 1_{\bar a}$. 
\Dfn

\prop\label{p.invertible}
$\phantom{v}$
\enu
	\item The inverse $a^{-1}$ of an element $a$ of $\sC$ is unique if it exists. 
	\item $\widehat{a^{-1}} = \bar a$ and $\overline{a^{-1}} = \hat a$. (Even if $\sC$ is proximal.) 
	\item All objects $a$ are invertible: $a^{-1} = a$. 
	\item Functors send invertibles to invertibles: $F(\theta^{-1}) = F(\theta)^{-1}$.
\Enu
\Prop

There are a few ways a generalized category may be partitioned into equivalence classes:

\dfn\label{d.epimonicisoclass}
For $a, b \in \sC,$ we have the following equivalence relations:
\enu
	\item $a$ and $b$ are in the same {\em monic class}, or {\em subobject}, $a \sim_m b$, if there exists invertible element $\theta \in \sC$ such that $a\theta = b$. 
	\item $a$ and $b$ are in the same {\em epic class}, or {\em quotient}, $a \sim_e b$, if there exists invertible element $\theta \in \sC$ such that $\theta a = b$; 
	\item $a$ and $b$ are in the same {\em iso class}, $a \sim b$, if there exist invertible elements $\theta_1, \theta_2 \in \sC$ such that $\theta_1 a = b \theta_2$. 
\Enu
\Dfn

Let $\Theta$ denote the set of all invertible elements in $\sC$. Define the symbol 
$$a\Theta := \set{a\cdot\theta \mid \theta \in \Theta \eand a\cdot\theta \downarrow},$$
and define the symbols $\Theta a, \Theta a\Theta$, etc. similarly. Then for $a,b \in \sC$,
$b$ belongs to the monic class of $a$ if and only if $b \in a\Theta$, %
$b$ belongs to the epic class of $a$ if and only if $b \in \Theta a$, %
and 
$b$ belongs to the iso class of $a$ if and only if $b \in \Theta a \Theta$. %
This notation is useful for back-of-the-envelope calculations, but it can be misleading: it need not be true that $\Theta f \Theta = \Theta g \Theta$, even if $f$ and $g$ are invertible. %


\dfn\label{d.epimoniciso}
An element $m$ of a generalized category $\sC$ is {\em monic} if $mf,mg \downarrow$ and $mf = mg$ implies $f = g$. %
An element $e$ in $\sC$ is {\em epi} if $fe,ge \downarrow$ and $fe = ge$ implies $f = g$. %
We say $a$ is {\em isomorphic} to $b$, denoted 
$$a \iso b,$$
if there exists an invertible element $\theta$ with $\bar\theta = a, \hat\theta = b$. %
\Dfn

If $a$ is monic and $a \sim_m b$, then $b$ is monic, and the $\theta$ given by the definition is unique. Similarly, if $a$ is epic and $a \sim_e b$. 


For every $a,b \in \sC$, $a$ is isomorphic to $b$ iff $1_a$ is in the same iso class as $1_b$, that is,
$$a \iso b \quad \ifaoif \quad 1_a \sim 1_b.$$
For $a,b$ objects, this becomes: 
$$a \iso b \quad \ifaoif \quad a \sim b.$$




\prop\label{p.isoskel}
Let $\sC$ be a generalized category. Then the set of iso classes forms a sharp category. %
The objects of this category are the iso classes of invertible elements of $\sC$. 
\Prop
\prf
Let $\tilde \sC$ be the set of iso classes of $\sC$, %
let $\tilde a, \tilde b$, ... denote elements in $\tilde \sC$. 
Define
$$\tilde a \cdot \tilde b := 
\set{\theta_1 a \theta_2 b \theta_3 \mid \theta_1,\theta_2,\theta_3 \text{ invertible, and $\theta_1 a \theta_2 b \theta_3 \, \downarrow$}}.$$
This is a partially defined map $\tilde \sC \cross \tilde \sC \to \tilde \sC$. For $a \in \sC$, let 
$$\tilde \source (\tilde a) := \widetilde{1_{\source a}},$$
$$\tilde \target (\tilde a) := \widetilde{1_{\target a}}.$$
These operations are well-defined: if $a = \theta_1 b \theta_2$, then $\bar a$ is isomorphic to $\bar b$, so, say, $\theta 1_{\bar b} \theta^{-1} = 1_{\bar a}$, so $\widetilde{1_{\source a}} = \widetilde{1_{\source b}}$, and similarly for $\tilde \target$. %

We take the order $\dleq$ on $\tilde \sC$ to be trivial, and we check Definition \ref{d.gencat}. %
The first four conditions are immediate: for (\ref{ax.gencat-comp-st}),  %
if $\tilde a, \tilde b \in \tilde \sC$, then $ \tilde a \tilde b \,\downarrow$. This occurs 
if and only if $ \set{\theta \in \sC \mid \theta: \bar a \to \hat b \text{ is invertible}}$ is nonempty, 
if and only if $ 1_{\bar a} \sim 1_{\hat b}$, 
if and only if $\tilde \source(\tilde a) = \tilde \target(\tilde b)$. %
Next, we observe that if $\tilde a$ is an element of the form $\tilde\source \tilde b$ or $\tilde \target \tilde b$ 
in $\tilde \sC$, then it must be of the form $\widetilde{1_b}$ for some $b \in \sC$, and 
$$\tilde\target(\widetilde{1_b}) = \tilde \source(\widetilde{1_b}) = \widetilde{1_b},$$
so $\tilde{1_b}$ is an object. Next, we have
$$\widetilde{1_a} \cdot \tilde b = \set{\theta_1 1_a \theta_2 b \theta_3 } = \set{\theta_4 b \theta_3} = \tilde b,$$
and similarly, $\tilde b \widetilde{1_a} = \tilde b$ whenever the product is defined. %
So $\tilde \sC$ is a sharp generalized category, in fact a one-category, after closing over $1_{()}$. %
The second statement is merely the observation that $a$ is invertible if and only if $\tilde a = \widetilde{1_{\source a}} = \widetilde{1_{\target a}}.$
\Prf

\dfn\label{d.catinvertibles}
We refer to the category $\tilde \sC$ of Proposition \ref{p.isoskel} as the {\em category of invertibles} of $\sC$. %
\Dfn

The {\em skeleton} of a generalized category $\sC$ is any full subcategory such that each element of $\sC$ is isomorphic in $\sC$ to exactly one element of the subcategory. %
Skeletons are unique up to isomorphism \cite{MacCW}. 
In the case of a category $\sC$, 
the category of invertibles expresses exactly the same data as a skeleton, but in a different way: any iso class that is an object in the category of invertibles contains not a set of invertibles in $\sC$ that are pairwise isomorphic, but instead, the set of all the isomorphisms that relate them pairwise to one another. %
On the other hand, an iso class that is an arrow in the category of invertibles is a noninvertible arrow $f \in \sC$ well-defined up to a commutative square with invertible columns. 

Since every element has an identity, thus taking the category of invertibles is the same as the operation of flattening (Proposition \ref{p.flattening})
followed by taking the skeleton, yielding the description just made in the previous paragraph. Thus it is perhaps natural to think of it as the ``category of identities'' of the generalized category.

It is also the case that a functor $F$ lifts to a functorial map $\tilde F$ on the category of invertibles.
Indeed, define 
$$\tilde F : \tilde \sC \to \tilde \sC',$$
via 
$$\tilde F(\tilde a) := \widetilde{F(a)}.$$
This is well-defined, as a consequence of (2) (which depends on the unital property of $F$):
$$\tilde F (\theta_1 a \theta_2) = \verywidetilde{F(\theta_1) F(a) F(\theta_2)} = \widetilde{F(a)}.$$
So we check functoriality: we have
$$\widetilde{1_{\source(F(a))}} = \widetilde{1_{F(\source(a))}} = \widetilde{F(1_{\source(a)})} = \tilde{F}(\widetilde{1_{\source(a)}}) = \tilde{F}(\source(\tilde a)),$$
and
$$\widetilde{1_{\source(F(a))}} = \source(\widetilde{F(a)}) = \source(\tilde{F}(\tilde{a})).$$
Similarly, 
$$\target(\tilde{F}(\tilde{a})) = \tilde{F}(\target(\tilde{a})).$$
And
$$\tilde{F}(\tilde{a}\tilde{b}) = \tilde{F}(\widetilde{a \theta b}) = \widetilde{F(a\theta b)} = \verywidetilde{F(a) \theta' F(b)} = \tilde{F}(\tilde{a}) \tilde{F}(\tilde{b}).$$
Finally, $\tilde{F}$ is unital since $\tilde \sC$ and $\tilde{\sC'}$ are categories.

A notion weaker than isomorphism arises from considering the categories of invertibles. 

\dfn\label{d.equivalent}
Generalized categories $\sC$ and $\sC'$ are {\em equivalent} if their categories of invertibles are isomorphic. 
\Dfn

This definition appeals directly to a comparison of the categories of invertibles. 
%
Now consider two functors $F,G: \sC \to \sC'$ that both define the same functor $\tilde{\sC} \to \widetilde{\sC'}$ on the categories of invertibles of $\sC$ and $\sC'$. This can only mean that there exist a pair of functions $\theta_1,\theta_2:\sC \to \sC'$ such that $\call{a \in \sC} \theta_i(a)$ is invertible for $i=1,2$, and for all $a \in \sC$, 
$$ \theta_1 (a) F(a) = G(a) \theta_2 (a) \,\downarrow.$$
If this holds we may write 
$$F \iso G.$$

\prop\label{p.equivalent}
Two generalized categories $\sC$ and $\sC'$ are equivalent if either of the following two equivalent conditions are satisfied. 
\enu
	\item Their categories of invertibles are isomorphic via a pair $\tilde F, \tilde G$, where $\tilde G = \tilde F^{-1}$, that come from functors
	 $F: \sC \to \sC'$ and $G: \sC' \to \sC$. 
	\item There exist two functors $F, G$ from $\sC \to \sC'$ ($\sC' \to \sC$, respectively) satisfying 
		$$F \of G \iso \id_{\sC'},$$
		$$G \of F \iso \id_{\sC}.$$
\Enu
\Prop

We can consider properties that a functor $\tilde F$ on the category of invertibles has as an ordinary functor, and view them as properties of the underlying functor $F$:

\dfn\label{d.essinjsurj}
A functor $F:\sC \to \sC'$ is {\em essentially injective} if it satisfies one of the following equivalent conditions, 
\enu
	\item $\tilde F$ is injective. 
	\item For $a,b \in \sC$, $F(a) = F(b)$ implies $a \sim b$.
\Enu
and $F$ is {\em essentially surjective} if it satisfies one of the following equivalent conditions: 
\enu
	\item $\tilde F$ is surjective.
	\item For $\alpha \in \sC'$, there exists $a \in \sC$ with $F(a) \sim \alpha$.
\Enu
\Dfn




From our initial investigation of equivalences between generalized categories, we arrived at the notion of equivalence via a pair of functors $F$ and $G$. We could, however, view this machinery (the pair ($\theta_1, \theta_2$)) as instead relating the two functors, and extend it:

\dfn\label{d.mof}
Let $\sC, \sC'$ be generalized categories, let $F,G:\sC \to \sC'$ be two functors. We say that a {\em morphism of functors} \cite{KaSc1} from $F$ to $G$ is a pair $(\theta_1,\theta_2)$ of maps $\sC \to \sC'$ satisfying, for all $a \in \sC$,
\begin{equation}\label{eq.mof}
\theta_1 (a) F(a) = G(a) \theta_2 (a) \, \downarrow
\end{equation}
Note that here, $\theta_1$ and $\theta_2$ are no longer presumed to be invertible. 
We may write the morphism of functors with the notation 
$(\theta_1, \theta_2): F \Rightarrow G$. 
\Dfn

Note that the maps $\theta_1$ and $\theta_2$ are maps from $\sC$ to $\sC'$, not from $\Ob(\sC)$ to $\sC'$ (cf. \cite{MacCW}). 

\begin{example}
Let $A = (a_{ij})$ be a matrix with coefficients in a ring $R$, and let $f:R \to S$ be a ring homomorphism. One naturally sets $f(A) = (f(a_{ij}))$, and doing this, one sees that
\eqn\label{eq.fdet}
\det(f(A)) = f(\det(A)).
\Eqn
This relation can be interpreted by observing that $GL_n$ is a functor from the category of rings to the category of groups, and likewise for the mapping that sends a ring to its group of units, and a ring homomorphism to the pointwise-identical homomorphism on the respective groups of units. So if $f:R \to S$, and writing $F(f)$ for the map defined above extending $f$ to a map on $GL_n (R)$, and $G(f)$ for the map changing $f$ to a map on the group of units, we have
$$\det() \of F(f) = G(f) \of \det()$$
by rewriting equation (\ref{eq.fdet}). From this expression we can read off the morphism of functors: 
$$\theta_1(f) = \det : GL_n(S) \to S^{\cross},$$
$$\theta_2(f) = \det : GL_n(R) \to R^{\cross}.$$
We see that in this example, $\theta_1$ and $\theta_2$ come from a single map $\theta$ on the objects (rings). This is not only typical of categories, it is guaranteed to happen. Indeed, if we return to the general situation of Definition \ref{d.mof}, inserting $a = 1_b$ into equation (\ref{eq.mof}) gives 
$$ \theta_1 (1_b) = \theta_2 (1_b)$$
for $b \in \sC$, so in particular, for all objects $b$,
$$ \theta_1 (b) = \theta_2 (b).$$
Thus $\theta_1$ and $\theta_2$ are identical on objects, and since one-categories have no higher morphisms, this single map on objects completely characterizes $(\theta_1, \theta_2)$. 
\end{example}

In the terminology of section \ref{s.nat} that follows, this means that a morphism of functors between functors relating categories is always {\em natural}. %
In the setting of generalized categories, we might suppose that this naturality property is a condition special to one-categories, %
since it does not appear to have any a priori motivation. %
However, 
the theory that results from dropping the naturality condition appears to be significantly weaker: %
\enu
	\item There is no strict 2-category of non-natural transformations, functors, and generalized categories. 
		Here, the wheel turns on the tiniest of pedestals: in the notation of Chapter \ref{c.prelude}, the relations
		$$\bar \alpha (X) = \overline{\alpha(X)},$$
		$$\hat \alpha (X) = \widehat{\alpha(X)}$$
		hold only in the natural setting. So we do not prove Fact 1. 
	\item While there is a notion of non-natural adjunction, there is no hom set bijection. %
		A key step in the proof uses the naturality of the unit and counit maps. %
		This in turn is used to prove that left adjoints are right exact. %
	\item Because there is no adjoint hom set bijection, some theorems relating equivalences of categories with 
		properties of functors 
		no longer hold. In particular a full, faithful, essentially surjective functor might not define an equivalence. 
\Enu
For these reasons, we do not take the development any further until we introduce naturality in the next section.

\section{Naturality}\label{s.nat}



In this section we establish the second of the two notions of equivalence we consider, namely
natural equivalence. %
As already noted, {\em the distinction between natural and non-natural vanishes in the case of categories.} %
Under natural equivalence, we obtain a 2-category of generalized categories, and in particular, an interchange law (Theorem \ref{t.functorcategory}). We can also establish, using the final lynchpin that naturality provides, the hom set bijection associated with adjoint pairs (Theorem \ref{t.nathombij}). Consequently the familiar rule that an equivalence between categories is given by a fully faithful essentially surjective functor carries over to generalized categories (Theorem \ref{t.ffes}). The full and faithful properties are tied to the naturality condition, which gives rise to maps not only on individual elements, but on entire hom sets. 

\dfn\label{d.nattrans}
Let $\sC, \sD$ be generalized categories, let $F,G: \sC \to \sD$. %
Let $(\theta_1, \theta_2) :F \Rightarrow G$ be a morphism of functors. %
We say that $(\theta_1, \theta_2)$ is {\em natural} or that $(\theta_1, \theta_2)$ is a {\em natural transformation} if, 
for every $a,b \in \sC$, %
$$\theta_1(a) = \theta_1 (b)$$
whenever $\hat a = \hat b$, and
$$\theta_2 (a) = \theta_2 (b)$$ 
whenever $\bar a = \bar b$.
\Dfn

Thus, naturality means that the function $\theta_1(a)$ can be replaced with the function $\hat a \mapsto \theta_1(1_{\hat a})$ of the element $\hat a$, and $\theta_2$ can be replaced with the function $\bar a \mapsto \theta_2(1_{\bar a})$ of the element $\bar a$. But, as noted in section \ref{s.elem}, 
$\theta_1(1_b) = \theta_2(1_b)$ 
for all elements $b$. 
Hence a natural transformation reduces to a single map $\theta: \sC \to \sC'$, from which $\theta_1$ and $\theta_2$ are immediately derived: 
$$\theta_1 (a) := \theta(1_{\hat a}),$$
$$\theta_2 (a) := \theta(1_{\bar a}).$$
We refer to a natural transformation $(\theta_1, \theta_2)$ by referring to this map $\theta$. In terms of $\theta$ the defining relation of a morphism of functors becomes
$$\theta(\hat f \hspace{0.8pt}) \cdot F(f) = G(f) \cdot \theta(\bar f \hspace{0.8pt}) \, \downarrow.$$

\dfn\label{d.naturaleq}
Two generalized categories $\sC$ and $\sC'$ are {\em naturally equivalent} 
 if they are equivalent via natural transformations
$$\theta:F\of G \iso \id_{\sC'},$$
$$\theta':G\of F \iso \id_\sC.$$
\Dfn

Naturally equivalent generalized categories are, in particular, equivalent (Definition \ref{d.equivalent}). %
With the extra condition of naturality, 
the way is clear to extend many 
justly well-known results of one-category theory \cite{MacCW} to the generalized setting: 

\thm\label{t.functorcategory}
The system given by all of the 
generalized categories, 
functors, 
and natural tranformations
forms a strict 2-category. 
\Thm
\prf
We define the products 
$$\theta_1 \vertof \theta_2,$$
$$\theta_1 \star \theta_2$$ 
just as in Chapter \ref{c.prelude}, and proceed as in the one-categorical case. %
\Prf

We include the naturality condition when defining adjoints: %

\dfn\label{d.adjunction}
Let $\sC$ and $\sD$ be generalized categories. An 
{\em adjunction} 
$(F,G,\eta,\varepsilon)$ is a pair of functors 
\vspace{-1ex}
$$
\begin{tikzcd}
	\sC \arrow[r, "F", shift left] \arrow[r, leftarrow, shift right, "G" below] & \sD \\
\end{tikzcd}\vspace{-3.5ex} 
$$
together with 
natural transformations
$$
\eta: \id_\sC \to G\of F, \quad \varepsilon: F \of G \to \id_\sD,
$$
satisfying the identities
\begin{align}\label{eq.triid}
(G \of \varepsilon) \vertof (\eta \of G) &= 1_G, \\
(\varepsilon \of F) \vertof (F \of \eta) &= 1_F,
\end{align}
where $1_F$ is the mapping $f \mapsto 1_{F(f)}$. %
Given an adjunction $(F,G,\eta,\varepsilon)$, $\eta$ is called the {\em unit} and $\varepsilon$ is called the {\em counit} of the adjunction. 
A natural equivalence $(\theta, \theta')$ is an {\em adjoint equivalence} if $\theta$ and $\theta'$ are the unit and counit of an adjunction. 
\Dfn

\thm\label{t.nathombij}
Let $\sC, \sD$ be generalized categories, and let $F,G: \sC \to \sD$ be functors. The following are equivalent:
\enu
	\item $(F,G,\eta,\varepsilon)$ forms an adjunction 
$
\begin{tikzcd}
	\sC \arrow[r, "F", shift left] \arrow[r, leftarrow, shift right, "G" below] & \sD \\
\end{tikzcd}\vspace{-1.5ex} 
$.
	\item 
For every $f$ in $\sC$ and $g$ in $\sD$, there is a 
bijection of sets
\begin{equation}\label{eq.adjbij}
\hom(F(f), g) \iso \hom(f, G(g)),
\end{equation}
that is natural in $f$ and $g$. This means that if $\phi_{f,g}$ is the bijection (\ref{eq.adjbij}), then 
for every $k: g \to g'$, and $h:f' \to f$, the following diagrams commute: 
$$
\begin{tikzcd}
\hom(F(f),g) \arrow[r, "\phi_{f,g}"] \arrow[d, "k_{*}"] & 
\hom(f,G(g)) \arrow[d, "G(k)_{*}"] \\
\hom(F(f),g') \arrow[r, "\phi_{f,g'}"] &
\hom(f,G(g')) 
\end{tikzcd}
\quad\quad
\begin{tikzcd}
\hom(F(f),g) \arrow[r, "\phi_{f,g}"] \arrow[d, "F(h)^{*}"] &
\hom(f, G(g)) \arrow[d, "h^{*}"] \\
\hom(F(f'), g) \arrow[r, "\phi_{f',g}"] &
\hom(f', G(g))
\end{tikzcd}
$$
Equivalently $\phi$ satisfies
\[
u \cdot F(v): F(f) \to g \eimplies \phi(u \cdot F(v)) = \phi(u) \cdot v,
\]
\[
v' \cdot v: F(f) \to g \eimplies \phi(v' \cdot v) = G(v') \cdot \phi(v).
\]
\Enu
\Thm
\prf
The proof is formally the same as in the one-categorical case (see \cite{MacCW}). 
\Prf

\dfn\label{d.faithfulfull}
Let $\sC,\sD$ be generalized categories, $F: \sC \to \sD$ a functor. 
For $a,b \in \sC$, let $F_{a,b}$ be the mapping on the domain $\hom(a,b)$ given by $f \mapsto F(f)$. %
We say that $F$ is {\em faithful} if %
for all $a,b$, $F_{a,b}$ is injective, 
and 
we say that $F$ is {\em full} if %
for all $a,b$, $F_{a,b}$ is surjective. 
\Dfn

Thus for example full means: if $\alpha, \beta$ in $\sD$ are of the form $F(a), F(b)$, for $a,b \in \sC$, and if $\gamma: \alpha \to \beta$, then $\gamma$ is of the form $F(c)$ for $c \in \sC$. 

\thm\label{t.ffes}
Let $\sC, \sD$ be generalized categories, and let $F: \sC \to \sD$ be a functor. The following are equivalent:
\enu
	\item $F$ is a natural equivalence,
	\item $F$ is a natural adjoint equivalence,
	\item $F$ is full, faithful, and essentially surjective.
\Enu
\Thm
\prf
The proof, much the same as in the one-categorical case, is left to the reader. 
\Prf


\section{Limits}\label{s.lim}

In this section we establish the elements of the theory of limits and colimits in sharp generalized categories. 
We consider limits with respect to mappings $I \to \sC$ as in Definition \ref{d.limit} that are weaker than functors. This, for example, allows us to form the shape of a product or coproduct of any set of elements in a generalized category. 

\dfn\label{d.limitfunctor}
Let $\sC, \sC'$ be generalized categories. A {\em functor up to objects} from $\sC$ to $\sC'$ is a map $F:\sC \to \sC'$ satisfying, for every $a,b \in \sC$,
\enu
	\item $F(ab) = F(a)(b)$,
	\item $F(a)$ is an identity in $\sC'$ if and only if $a$ is an identity in $\sC$,
	\item $F(\source(a)) = \source(F(a))$ unless $a$ is an object of $\sC$,
	\item $F(\target(a)) = \target(F(a))$ unless $a$ is an object of $\sC$. 
\Enu
\Dfn

\dfn\label{d.cone}
Let $\sC$ be a generalized category, $I$ a generalized category (the index of a cone needs only be a set, but in practice it is always a (generalized) category). A {\em cone} in $\sC$ with index $I$ is a map $\sigma:I \to \sC$ such that 
$$\text{for all $i,j \in I$, } \overline{\sigma(i)} = \overline{\sigma(j)}.$$ %
Dually, {\em cocone} in $\sC$ with index $I$ is a map $\sigma:I \to \sC$ such that 
for all $i,j \in I$, $\widehat{\sigma(i)} = \widehat{\sigma(j)}.$ %
A cone or cocone is {\em finitely generated} if the index set $I$ is finitely generated (Definition \ref{d.figen}). %
This common source is the {\em vertex} of the cone, and the {\em vertex} of a cocone is the common target. %
Given a cone or cocone $\pi$, we may refer to $\pi(i)$ for some $i \in I$ as a {\em member} of the cone. %
\Dfn

\dfn\label{d.limit}
Let $\sC, I$ be generalized categories. Let $\alpha: I \to \sC$ be a functor, possibly only a functor up to objects. %
A cone is said to be {\em over (or below) the base $\alpha$} if %
\enu
	\item $\widehat{\pi(i)} = \alpha(i)$, for all $i \in I$, %
	\item for all $i \in I$, $\pi(\hat i) = \alpha(i) \pi(\bar i)$. %
\Enu
A {\em limit} of $\alpha$ is a cone $\pi:I \to \sC$ below the base $\alpha$ such that for any cone $\tilde \pi:I \to \sC$ over the same base $\alpha$, there is a unique $\lambda \in \sC$ such that $\tilde \pi = \pi \vertof \lambda$. (Here, $\pi \vertof \lambda$ is the map defined by $(\pi \vertof \lambda)(i) = \pi(i) \cdot \lambda.$) 

Dually, a cocone is said to be {\em over (or below) the base $\alpha$} if 
\enu
	\item $\overline{\pi(i)} = \alpha(i)$, for all $i \in I$,
	\item for all $i \in I$, $\pi(\bar i) = \pi(\hat i) \alpha(i)$.
\Enu
A {\em colimit} of $\alpha$ is a cocone $\pi:I \to \sC$ such that for any cone $\tilde \pi:I \to \sC$ over the base $\alpha$, there is a unique $\lambda \in \sC$ such that $\tilde \pi = \lambda \vertof \pi$. Here, $\lambda \vertof \pi$ is the map defined by $(\lambda \vertof \pi)(i) = \lambda \cdot \pi(i)$, as before. 
\Dfn

Thus a cone fits a pattern as in the following Figure:
$$
\begin{tikzcd}[column sep=1.5em]
\alpha(\bar i) \arrow{rr}{\alpha(i)} && \alpha(\hat i) \\
 & \text{(vertex)} \arrow{ul}{\pi(\bar i)} \arrow{ur}[swap]{\pi(\hat i)} 
\end{tikzcd}
$$


The word limit is often used to refer to the domain of the cone, and similarly colimit is used to refer to the codomain of the cocone. 
The terms {\em product}, {\em equalizer}, {\em coproduct}, {\em coequalizer}, etc. retain their meaning from ordinary categories, referring to limits based on diagrams $\alpha: I \to \sC$ of the same shape as in the one-categorical case, and where $\alpha$ may be a functor only up to objects. %
We follow standard terminology and say that a generalized category {\em has finite limits} if there is a limit cone for every finitely generated diagram $\alpha: I \to \sC$, and dually for colimits. 


We denote the set of limits of the functor $\alpha:I \to \sC$ by $\lim(\alpha,I)$ or just $\lim \alpha$. We denote the colimit $\colim(\alpha,I)$ or simply $\colim(\alpha)$. %

If $\sC$ is a generalized category, there exist (finitely generated) diagrams $J \to \sC$ that cannot be defined and do not exist in an ordinary category. 
However, we still have:

\thm\label{t.prodeqsuffice}
Let $\sC$ be a generalized category. For $\sC$ to have all finite limits, it suffices that $\sC$ has all finite products and equalizers. 
\Thm
\prf
We proceed by induction on the height of finitely generated diagrams $\alpha: I \to \sC$. 
A finitely generated diagram of height $0$ is a finite product, hence it has a limit cone in $\sC$ by hypothesis. %
Suppose that all finitely generated diagrams of height $k \geq 0$ have a limit cone, and 
let $\alpha: I \to \sC$ be a diagram of height $k + 1$. %
Define 
$$\alpha^{\leq k}$$
to be $\alpha$ restricted to the generalized category $I^{\leq k}$ formed by taking the collection of all elements of $I$ of height $\leq k$, along with all identities of $I$. 
It is easy to see that $I^{\leq k}$ is closed under composition, thus it is a generalized category. Therefore $\alpha^{\leq k}$ is a diagram on $\sC$, and by hypothesis, has a limit cone $\sigma^{\leq k}$ with vertex, say, $L^{\leq k}$. %
%
%
%
Consider $flat(I^{\leq k})$, the flattening of $I^{\leq k}$ to a zero-category (section \ref{s.elem}). %
The diagram $flat(\alpha^{\leq k}): flat(I^{\leq k}) \to \sC$ induced by $\alpha^{\leq k}$ %
is a diagram of height zero, so it has a limit cone $\sigma^{\leq k, flat}$, with vertex, say, $L^{\leq k, flat}$. 
The cone $\sigma^{\leq k}$ on $I^{\leq k}$ induces a cone on $flat(I^{\leq k})$, 
so there exists a universal arrow 
$$u_1: L^{\leq k} \to L^{\leq k, flat}.$$
Now let $I^{k+1, flat}$ be the flattened (to a zero category) elements of $I$ of height $k+1$. %
The diagram $\alpha$ induces a diagram $\alpha^{k+1, flat}$ on $I^{k+1, flat}$, %
defined by 
$$\alpha^{k+1, flat} (i) := \target(\alpha(i)).$$
This diagram (of height zero) has a limit cone $\sigma^{k+1, flat}$ 
with vertex, say, $L^{k+1, flat}$. %
For $i \in I$ of height $k+1$, let $\pi_i$ be the element in $\sC$ which is the projection 
$$\pi_i : L^{\leq k, flat} \to \target(\alpha(i)),$$
coming from the diagram $\sigma^{\leq k, flat}$ on $I^{\leq k, flat}$ (where our notation hides this fact about $\pi_i$). 

The previous cone $\sigma^{\leq k, flat}$ with vertex $L^{\leq k , flat}$ itself has projection arrows to the elements $\target(\alpha(i))$ 
as $i$ ranges over $\alpha^{k+1, flat}$. 
Therefore, there is a universal arrow 
$$u_2: L^{\leq k, flat} \to L^{k+1, flat}.$$
Moreover, for each $i$ of height $k+1$, there is also a projection arrow to the element $\source(\alpha(i))$, and 
composing each of these projection arrows with $\alpha(i)$ gives a second cone with the same vertex $L^{\leq k, flat}$ %
on the diagram $\alpha^{k+1, flat}$. 
So we may again find a universal arrow
$$u_3: L^{\leq k, flat} \to L^{k+1, flat},$$
by applying the universal property of the limit with vertex $L^{k+1, flat}$ a second time. %
We compose $u_2$ and $u_3$ with $u_1$ to form parallel arrows, and take the equalizer:
\[
\begin{tikzcd}
L \arrow[r, "e"] &
L^{\leq k} \arrow[r, "u_{1}"] &
L^{\leq k {,} flat}
\arrow[r, shift left, "u_{2}"]
\arrow[r, shift right, swap, "u_{3}"] &
L^{k+1 {,} flat}
\end{tikzcd}
\]
Now we define, for $i$ in $I$ of height $\leq k+1$,
$$\sigma^{\leq k+1} (i) := \pi_i \cdot u_1 \cdot e.$$
We claim that this is a limit cone for the diagram $\alpha^{\leq k+1}: I^{\leq k+1} \to \sC$. 
Since we pass through $e$ to reach $L^{\leq k+1}$, $\sigma^{\leq k+1}$ satisfies 
$\sigma^{\leq k+1} (\hat i) = \alpha^{\leq k+1} (i) \cdot \sigma^{\leq k+1} (\hat i)$, 
hence is a limit cone. 
Suppose that 
$\tilde \sigma^{\leq k+1}: I^{\leq k+1} \to \sC$ is a diagram with vertex, say, $\tilde L$ satisfying 
$\tilde \sigma^{\leq k+1} (\hat i) = \alpha^{\leq k+1} (i) \tilde \sigma^{\leq k+1} (\bar i).$ 
Then $\tilde \sigma^{\leq k+1}$ restricts to a cone on $\alpha^{\leq k}$, hence there is a universal arrow
$$\tilde e: \tilde L \to L^{\leq k}.$$
Because $\tilde \sigma^{\leq k+1}$ has the limit property even at the height $k+1$, $\tilde \sigma^{\leq k+1}$ satisfies $u_2 \cdot u_1 \cdot \tilde e = u_3 \cdot u_1 \cdot \tilde e$, and thus $\tilde e$ factors through $e$ uniquely, as desired. 
\Prf

\dfn
Let $F: C \to C'$ be a functor. Then $F$ {\em preserves limits} or is {\em left exact} if for every functor $\alpha:I \to C$, 
$$F(\lim(\alpha)) \lies \lim(F\of \alpha).$$
Dually, $F$ {\em preserves colimits} or is {\em right exact} if for every functor $\alpha:I \to C$,
$$F(\colim(\alpha)) \lies \colim(F\of \alpha).$$
$F$ is said to {\em create limits} if for every element $\pi \in \lim(F \of \alpha)$, there exists a unique $\pi' \in \lim(\alpha)$ such that $F(\pi') = \pi$. Dually, $F$ is said to {\em create colimits} if for every element $\pi \in \colim(F \of \alpha)$, there exists a unique $\pi' \in \colim(\alpha)$ such that $F(\pi') = \pi$.
\Dfn

For example, the hom functor 
$$b \mapsto \hom(-,b)$$
preserves limits. 
Dually, the contravariant hom functor 
$$a \mapsto \hom(a,-)$$
preserves colimits. These functors may be extended to generalized categories.  

\thm\label{t.hlalex}
Let $F: \sC \to \sD$ be a functor between generalized categories $\sC$ and $\sD$. 
Then if $F$ has a left adjoint $G: \sD \to \sC$, then it is left exact. 
\Thm
\prf
Like the proof for categories, the proof for generalized categories relies on naturality of the adjoints via the bijection (\ref{eq.adjbij}). 
\Prf

The dual statement to \ref{t.hlalex} is immediate: a functor with a right adjoint is right exact. %

\section{Globular Sets}\label{s.globularsets}

It is worthwhile to remark on the relationship between generalized categories and globular sets. Recall that
a globular set is a presheaf of shape $\mathbb{G}$ (that is, a functor $\mathbb{G}^{\text{op}} \to \Set$), where $\mathbb{G}$ is the category of natural numbers $n \geq 0$ together with maps
$$
\begin{tikzcd}
0 \ar[r, shift left=1ex, "\sigma_{0}"] \ar[r, shift right=1ex, "\tau_{0}"']
	& 1 \ar[r, shift left=1ex, "\sigma_{1}"] \ar[r, shift right=1ex, "\tau_{1}"']	
	& 2 \ar[r, shift left=1ex, "\sigma_{2}"] \ar[r, shift right=1ex, "\tau_{2}"']
	& ...
\end{tikzcd}
$$
subject to the relations 
$
\sigma_{i+1} \of \sigma_{i} = \tau_{i+1} \sigma_{i}, \quad \tau_{i+1} \of \tau_{i} = \sigma_{i+1} \of \tau_{i},
$
for $i \geq 0$.


\dfn
Let $\sC$ be a generalized category. A {\em $k$-cell} in $\sC$ is an element $f$ of $\sC$ such that for every $k$-element sequence $\vec s$ of operations $\source$ and $\target$ that satisfy when applied to $f$, 
\enu
	\item $\source^k f$ and $\target^k$ are objects, and $\source^{k-n} f$ and $\target^{k-n}$ are not objects, for all $0 \leq n \leq k$, 
	\item $\source \target f = \source \source f$ and $\target \source f = \target \target f$,
	\item $\source f$ and $\target f$ is are $(k-1)$-cells.
\Enu
An element $f$ of a generalized category $\sC$ is {\em cellular} if $f$ is a $k$-cell for some $k \geq 1$, %
and a generalized category $\sC$ is {\em cellular} if every element of $\sC$ is cellular. 
\Dfn


\prop
There is an equivalence (given by a forgetful-free adjunction) between sharp, cellular generalized categories and the category of globular sets. 
\Prop
\prf
To prove this, we must be sure clarify the statement: when referring to sharp, cellular generalized categories, we refer not to the full subcategory but to the category whose morphisms $F: \sC \to \sD$ are subject to the extra condition 
\enu
	\item for all $a \in \sC$, $s(F(a)) = F(a)$ implies $s a = a$.
\Enu
This says we cannot map $k$-cells for $k > 0$ to $0$-cells. %
Then let $\dim(a) := \min\set{n \mid \source^n a = \source^{n+1} a }$. %
Define a mapping 
$$\sC \mapsto (n \mapsto \set{a \in \sC \mid \dim a = n}).$$ %
to the category of globular sets, for a sharp cellular generalized category $\sC$. This is the desired equivalence. 
\Prf

Examples of noncellular generalized categories are abundant, for example arising from the theory of trees and related notions, see for example \cite{ElBlTi1}. 


\section{Conclusion}\label{s.gencat-conc}

There are numerous advanced notions of category theory that have not yet made an appearance in our development, 
for example, ends, coends, monads, Kan extensions, to name only a few. 
Our investigation has yielded the following observations: %
There exists a generalization of category theory. 
More precisely, there exists a theory of functors, natural transformations, adjoint pairs, limits, and colimits for generalized categories. %
Still more precisely, there are two generalizations that are combined into one larger one:  
First, by allowing an approximate operation of composition (i.e., proximal categories),
and second, by allowing generalized higher cells. 
We have seen that the structure of limits and natural transformations is similar to the structure as it arises in ordinary one-categories, 
so that, surprisingly, perhaps, the proximal structure has little effect on aspects of the theory. %
We have investigated a notion of non-natural transformation suggested by the one-categorical case where naturality is not a necessary assumption, and we have found that the device of non-natural equivalence is not sufficiently strong. Therefore, we have argued that naturality must be an explicit assumption in the generalized setting. 
Thus, we have both extended the boundaries of category theory, and made note of some limits to further extensions of the new boundaries we have drawn, which strengthens the case for our particular approach. %
With the foundations developed here, it is possible to go further. 

\pagebreak
\singlespacing
\chapter{Monads and Generalized Categories}\label{c.genm}
\doublespacing
\vspace{10ex}

\section{Introduction}\label{s.genm-i}

The monad abstraction, which arose out of pure mathematics, has proven to be an enduring and versatile notion through a wealth of 
 applications to algebra, logic, and computer science in recent 
 decades. %
Generalized categories have been introduced due to their inherent interest as mathematical objects---they are a class of abstractions far more general than ordinary categories, and yet, ordinary category theory continues to apply to them---and due to their own potential applications to computer science, where they are suggested as a model for higher order programming, in particular higher-kinded types. %
Therefore it is natural to place these two notions, the monad abstraction and the generalized category, side by side to see what can be said about their relationship. 

Order-enriched categories (suitable ones, with the property of $\omega$-completeness) are the foundation of the categorical domain theory developed by Smyth, Plotkin, Wall and others \cite{SmPl1,WaD1}. A special case of a generalized category is an order-enriched category, thus their work carries over to the present setting; one of our goals is to further clarify this relationship. The motivation for going to a further level of generality is twofold: first, the theoretical material is in fact well-behaved in the wider setting, indicating that it is appropriate to do so---in other words, generalized categories have an associated mathematical theory---and second, the observation that combining type coercion with the order of approximation (the ordering motivating the study of order-enriched categories in the context of domain theory) yields a transtive relation, opening a door to applications. Another door to applications is via the Curry-Howard-Lambek correspondence and the categorical semantics of Moggi. 



\section{Monads in Generalized Categories}\label{s.mon}

\dfn\label{d.monad}
Let $\sC$ be a generalized category. A {\em monad} on $\sC$ is a structure $(T,\eta,\mu)$, where $T: \sC \to \sC$ is a functor, and $\eta$ and $\mu$ are (order-preserving) natural transformations $\id_\sC \to T$ and $T^2 \to T$, respectively, such that the following hold:
\enu
	\item $\mu \vertof (T \of \mu) = \mu \vertof (\mu \of T)$ \label{ax.monmagic1}
	\item $\mu \vertof (T \of \eta) = \mu \vertof (\eta \of T) = 1_T$, \label{ax.monmagic2}
\Enu
where $1_T$ denotes the mapping $f \mapsto 1_{T(f)}$. %
A monad is said to satisfy the {\em monic condition} \cite{MoI1} if
for all $x,y$ in $\sC$, 
$$\eta(\hat x) x \dleq \eta(\hat y) y \eimplies x \dleq y.$$
It follows that $\eta(x)$ is a monic for all $x \in \sC$. %
\Dfn

A monad is often referred to 
simply by the symbol for the functor $T$, with $\eta = \eta_T$, $\mu = \mu_T$ understood. 

\exa\label{x.mon}
Consider the generalized category $\Set f$, that is, the construction $\sC f$ of Example \ref{x.setf} applied to the category of sets. If we fix a group $G$, the usual monad on $\Set$ generated by $G$ extends directly to $\Set f$, giving a monad on $\Set f$: 
$$T(\ff) = \id_G \times \ff$$
$$
\eta(\ff) = \text{ the tree $\fh$ with $\tleft(\fh) = T(\ff)$, $\tright(\fh) = \ff$, $\troot(\fh) = x \mapsto (1_G, x)$ }
$$
$$
\mu(\ff) = \text{ the tree $\fh$ with $\tleft(\fh) = T^2 (\ff)$, $\tright(\fh) = T(\ff)$, $\troot(\fh) = (g,h,x) \mapsto (gh, x)$ }.
$$
\Exa


Now we establish some background needed for the Kleisli contruction: 
a generalized category equipped with an adjunction is also equipped with a monad. 

\prop\label{p.monadj}
Let $\sC$ be a generalized category, and let $(F,G,\eta,\epsilon)$ be an (order-preserving) adjunction on $\sC$. Then there is a monad on $\sC$ given by $T = G \of F$ with $\eta_T = \eta, \mu_T = G \of \epsilon \of F$.
\Prop
\prf
The proof is standard, but the definitions are not. %
The unit $\eta$ and counit $\epsilon$ are indeed order-preserving natural transformations. Hence so is $\eta_T$ and $\mu_T$, which we denote $\eta$ and $\mu$ for the rest of the proof. 
We still have horizontal composition $\hozof$ given by $\beta \star \alpha = (\beta \of \bar{\alpha}) \vertof (\hat{\beta} \of \alpha)$, and the interchange law. These imply (a sharp expression)
$$\epsilon \hozof \epsilon = \epsilon \vertof (F \of G \of \epsilon) = \epsilon \vertof (\epsilon \of F \of G).$$
So applying $G$ and $F$ on the left and right side respectively, 
\begin{align*}
G \of (\epsilon \vertof (F \of G \of \epsilon)) \of F 
	&= G \of (\epsilon \vertof (\epsilon \of F \of G) ) \of F
\end{align*}
and so
\begin{align*}
(G \of \epsilon \of F) \vertof (G \of F \of G \of \epsilon \of F)
	&= (G \of \epsilon \of F) \vertof (G \of \epsilon \of F \of G \of F),
\end{align*}
or $\mu \vertof (T \of \mu) = \mu \vertof (\mu \of T).$ %
We need to check if $\mu \vertof (T \of \eta) = \mu \vertof (\eta \of T) = 1_T.$
This follows similarly from the triangular identities. 
Therefore $(T, \eta, \mu)$ is an (order-preserving) monad. 
\Prf

We will now verify that, as in the ungeneralized case, every monad comes from an adjunction. First, the Kleisli construction may be extended to the generalized setting:

\dfn[Kleisli Construction]\label{d.kleisliconst}
If $\sC$ is a generalized category, and $T = (T,\eta, \mu)$ is a monad on $\sC$, set
$$\sC_T := \set{(y,f) \mid T(y) = \hat f, \,(\star)\,}$$
where $(\star)$ is a condition (detailed in the proof) that may be safely ignored on a first reading, and
\begin{equation*}
\begin{cases}
	\, \overline{(y,f)} = (\hat{\bar{f}}, \eta(\hat{\bar{f}})\bar f), &\\
	\, \widehat{(y,f)} = (\hat y , \eta(\hat y) \cdot y), &\\
	\, (z,g) \cdot (y,f) = (z, \mu(z) \cdot T(g) \cdot f), \hspace{20pt}\text{ if } \bar g \dleq y.  &\\
\end{cases}
\end{equation*}
\Dfn

\thm\label{t.kleisli}
Let $\sC$ be a generalized category with identities. Then if $T$ satisfies the monic condition, then $\sC_T$ is a generalized category (not necessarily with identities). There exists an order-preserving adjunction
$$
\begin{tikzcd}
	\sC \arrow[r, "F", shift left] \arrow[r, leftarrow, shift right, "G" below] & \sC_T^0 \\
\end{tikzcd}\vspace{-3.5ex}
$$
where $\sC_T^0$ is the image of $F$ in $\sC_T$, such that the monad generated by $(F,G,\eta, \epsilon)$ is $T$. 
\Thm
\prf
Define a poset structure by setting
$$(z,g) \dleq (y,f) := (z \dleq y, g \dleq f).$$
Now suppose that $(z,g), (y,f)$ are two elements of $\sC_T$ such that $\overline{(z,g)} = \widehat{(z,g)}.$ Then 
$$(\hat{\bar{g}}, \eta(\hat{\bar{g}}) \bar g) = (\hat y, \eta(\hat y) y).$$
So $\eta(\hat z) \bar g = \eta(\hat y) y.$ By the monic condition, 
$$\bar g = y.$$ 
So we have
\begin{align*}
((w,h)\cdot (z,g)) \cdot (y,f) 
	&= (w,\mu(w)T(h)\cdot g) \cdot (y,f) \\
	&= (w,\mu(w) \cdot T(\mu(w)) T^2 (h) \cdot T(g)\cdot f) \\
	&= (w, \mu(w) \cdot \mu(T(w))\cdot T^2 (h) \cdot T(g) \cdot f), \text{ and $T(w) = \hat h$} \\
	&= (w, \mu(w) \cdot T(h) \cdot \mu(\bar h) \cdot T(g) \cdot f) \text{ and $\bar h = z$} \\
	&= (w, \mu(w) \cdot T(h) \cdot (\mu(z) \cdot T(g) \cdot f)) \\
	&= (w,h) \cdot (z, \mu(z) \cdot T(g) \cdot f) \\
	&= (w,h) \cdot ((z,g) \cdot (y,f)), 
\end{align*}
after making the choice of $y$ to represent $(y,f)$. 

Next,
\begin{align*}
\overline{(z,g) \cdot (y,f)}
	&= \overline{(z, \mu(z) \cdot T(g) \cdot f)} \\
	&= (\verywidehat{\overline{\mu(z) \cdot T(g) \cdot f}}, \eta( \widehat{\overline{\text{same}}}) \cdot \overline{\text{same}}) \\
	&= (\widehat{T(z)}, \eta(\widehat{T(z)})\cdot T(z)) \\
	&= (\hat{\bar{g}}, \eta(\hat{\bar{g}}) \bar{g}) = \overline{(z,g)}, 
\end{align*}
and
\begin{align*}
\verywidehat{(z,g) \cdot (y,f)} 
	&= \verywidehat{(z,\mu(z) \cdot T(g) \cdot f) } \\
	&= (\hat z, \eta(\hat z) \cdot z) \\
	&= \widehat{(z,g)}.
\end{align*}

Now if $\overline{(z,g)} \not\dleq \widehat{(y,f)},$ they do not compose in $\sC_T$ by definition. On the other hand, if $\overline{(z,g)} \dleq \widehat{(y,f)},$ we must guarantee that they may be composed. If we have
$$
\overline{(z,g)} \dleq \widehat{(y,f)},
$$
then
$$
\eta(\hat{\bar{g}}) \cdot \bar{g} \dleq \eta(\hat{y}) \cdot y$$
so 
$$T(\eta(\hat{\bar{g}})) \cdot T(\bar{g}) \dleq T(\eta(\hat{y})) \cdot T(y)$$
so 
$$\mu(\hat{\bar{g}}) \cdot T(\eta(\hat{\bar{g}})) \cdot T(\bar{g}) \dleq \mu(\hat{y}) \cdot T(\eta(\hat{y}) \cdot T(y)$$
so $T(\bar g) \dleq \hat f$, as desired, using monotonicity of $\mu$. 

Now if $(y,f) \in \sC_T$ is an object, then $(y,f) = (y, \eta(y))$ for some $y \in \Ob(\sC)$. For then 
$$(y,f) = (\hat{\bar{f}}, \eta(\hat{\bar{f}}) \cdot \bar{f}) = (\hat{y}, \eta(\hat{y})\cdot y)$$
So 1. $y = \hat{\bar{f}}$ and 2. $f = \eta(\hat{\bar{f}}) \bar{f}$ and 3. $y = \hat{y}$ and 4. $f = \eta(\hat y) \hat y$. 3. and 4. with $\eta$-monotonicity gives $\bar f = y$. Taking source in 2. gives $y = \bar y$. So with 3., $y$ is an object in $\sC$. Now 4. shows that $f = \eta(y)$. %
Now we can check that, using $\eta$-monotonicity (twice!), $(y, \eta(y))$ satisfies axiom (\ref{ax.gencat-object-id}) of Definition \ref{d.gencat}.

Now, if $(z,g) \dleq (y,f),$ then $\hat{\bar{g}} \dleq \hat{\bar{f}}$. So $\eta(\hat{\bar{g}}) \dleq \eta( \hat{\bar{f}}).$ And $\bar{g} \dleq \bar{f},$ so $\eta(\hat{\bar{g}}) \bar g \dleq \eta(\hat{\bar{f}}) \bar{f}$. Similarly, $\eta(\hat z) z \dleq \eta(\hat y) y.$ So $\overline{(z,g)} \dleq \overline{(y,f)}$, and $\widehat{(z,g)} \dleq \widehat{(y,f)}.$ 

If $(z,g) \dleq (y,f),$ and $(u,k) \dleq (w,h),$ it is similarly verified that $(z,g) \cdot (u,k) \dleq (y,f) \cdot (w,h)$. 

The last axiom to check is $a \dleq b$ then $1_a \dleq 1_b$. The condition $(\star)$ in Definition \ref{d.kleisliconst} is the following: 
$$(\star): (y,f) = 1_{(x,u)} \eimplies \there v \text{ such that } (x,u) = (\hat{v}, \eta(\hat{v}) \cdot v).$$
This condition specifies that we throw away all identities in $\sC_T$ that are not in the image of $F_T$, defined below. We will use this now, and we will use it again below. 
Let $(y,f) \dleq (z,g)$, and suppose that both $(y,f)$ and $(z,g)$ have identities $1_{(y,f)}$ and $1_{(z,g)}$, respectively. By the condition $(\star)$ we may write $(y,f) = (\hat u, \eta(\hat u) u)$ and $(z,g) = (\hat v, \eta(\hat v) v)$ for some $u,v \in \sC$. Then it is not difficult to see that $1_{(y,f)} = (u,\eta(u))$ and $1_{(z,g)} = (v, \eta(v))$. Now 
we apply the full strength of the monic condition to conclude that $u \dleq v$, from which axiom (Identities) of Definition \ref{d.gencat} follows. 
So $\sC_T$ is a generalized category. 

Next, define mappings
$$F_T (f) := (\hat f, \eta(\hat f) \cdot f),$$
$$G_T(y,f) := \mu(y) \cdot T(f).$$
We have 
\begin{align*}
F_T (g \cdot f) 
	&= (\widehat{g \cdot f}, \eta(\widehat{g \cdot f}) \cdot g \cdot f) \\
	&= (\hat g, \eta(\hat g) \cdot g \cdot f) \\
	&= (\hat g, \id_{T(\hat g)} T(g) \cdot \eta(\bar g) \cdot f), \text{ and $\bar g = \hat f$} \\
	&= (\hat g, \mu(\hat g) \cdot T(\eta(\hat g)) \cdot T(g) \cdot \eta(\hat f) f) \\
	&= (\hat g, \mu(\hat g) \cdot T(\eta(\hat g) \cdot g) \cdot \eta(\hat f)\cdot f) \\
	&= (\hat g, \eta(\hat g) \cdot g) \cdot (\hat f, \eta(\hat f) \cdot f) \\
	&= F_T(g) \cdot F_T(f).
\end{align*}
Now if $f \dleq g$, then $\hat f \dleq \hat g$, and by substitutivity (Definition \ref{d.gencat}) and monotonicity of $\eta$, %
$(\hat f, \eta(\hat f)f) \dleq (\hat g, \eta(\hat g) g).$ So $F_T$ is monotonic. We can also check that 
$$F_T (\hat f) = (\hat{\hat{f}}, \eta(\hat{\hat{f}}) \cdot \hat{f}) = \verywidehat{(\hat{f}, \eta(\hat{f}) \cdot f)} = \widehat{F_T (f)}$$
and 
$$F_T (\bar f) = (\hat{\bar{f}}, \eta( \hat{\bar{f}}) \bar{f}) = \overline{F_T (f)}.$$


Now suppose that $f \in \sC$ is a subject. Then it has an identity $1_f$. The image $F(1_f)$ of this element is $(f,\eta(f))$, the identity of $F(f)$. 

So $F_T$ is a functor $\sC \to \sC_T$. The image $\sC_T^0$ of $\sC$ under $F$ is a generalized category with identities, 
for it is easy to see in general that the image under a functor of an element that is not a subject is not a subject.

We also have 
\begin{align*}
G_T ((z,g) \cdot (y,f)) 
	&= G_T (z, \mu(z)T(g) f) \\
	&= \mu(z) \cdot T(\mu(z) \cdot T(g) \cdot f) \\
	&= \mu(z) \cdot T(g) \cdot \mu(\bar g) \cdot T(f) \text{ and $\bar g = y$} \\
	&= G_T (z,g) \cdot G_T (y,f).
\end{align*}
And
\begin{align*}
G_T (\overline{(y,f)} 
	&= G_T (\hat{\bar{f}}, \eta(\hat{\bar{f}}) \bar{f}) \\
	&= \mu (\hat{\bar{f}}) \cdot T(\eta ( \hat{\bar{f}}) \cdot \bar{f}) \\
	&= T(\bar{f}) \\
	&= \overline{T(f)} \\
	&= \overline{\mu(y) T(f)} \\
	&= \overline{G_T (y,f)}.
\end{align*}
Similarly
\begin{align*}
G_T (\widehat{(y,f)}) 
	&= \mu(\hat y) \cdot T(\eta(\hat y))T(y) \\
	&= T(y) \\
	&= \verywidehat{\mu(y)T(f)} \\
	&= \verywidehat{G_T (y,f)}.
\end{align*}
For monotonicity, if $(y,f) \dleq (z,g)$,
$$G_T (y,f) = \mu(y) T(f) \dleq \mu(z) T(g) = G_T (z,g).$$
For identities, we apply $(\star)$ again:
$$G_T (1_{(y,f)}) = G_T (v, \eta(v)) = \mu(v)T(\eta(v)) = 1_{T(v)}$$
and indeed, $T(v) = G_T (\hat v, \eta(\hat v) v) = G_T (y,f).$ So $F_T$ and $G_T$ are functors. We shall make no distinction between these and their appropriate restrictions to $\sC_T^0$.

Next, define a map
$$\phi: \hom(F(f), (z,g)) \to \hom(f, G(z,g))$$
by $\phi(u,k) := k$. Clearly, $\phi$ is order-preserving. 
Next, $\phi$ is injective, for if $\phi(u,k) = \phi(u', k)$, then $\eta(z) \cdot u = g$ and $\eta(z) \cdot u' = g$. So $\eta(z) \cdot u = \eta(z) \cdot u'$. Hence $u = u'$ by the monic condition. 

It is surjective if it can be shown that for every $k \in \hom(f,G(z,g))$ there is a $u$ such that $\phi(u,k) = k$, or in other words, there is a $u$ such that %
$\eta(z) \cdot u = g$. 
If $(z,g)$ is in the image of $F$, then we obtain the desired $u$. So $\phi$ is surjective onto $\sC_T^0$. 

Now we check naturality of $\phi$. We have
\begin{align*}
\phi((w,h) \cdot (u,k)) 
	&= \phi(w, \mu(w T(h) \cdot k) \\
	&= \mu(w) T(h) \cdot k \\
	&= G(w,h) \cdot \phi(u,k),
\end{align*}
and
\begin{align*}
\phi(w,l_1) \cdot F(l_2)) 
	&= \phi((v,l_1) \cdot (\hat{l_2}, \eta( \hat{l_2}) l_2) \\
	&= \phi((v, \mu(v) \cdot T(l_1) \cdot \eta(\hat{l_2}) l_2)) \text{ and $\bar{l_1} = \hat{l_2}$} \\
	&= \mu(v) \cdot T(l_1) \cdot \eta(\bar{l_1}) \cdot l_2 \\
	&= \mu(v) \eta(\hat{l_1}) \cdot l_1 \cdot l_2 \text{ and $\hat{l_1} = T(v)$} \\
	&= l_1 l_2 \\
	&= \phi(v, l_1) \cdot l_2.
\end{align*}
So we have an order-preserving adjunction between $\sC$ and $\sC_T^0$. Applying Proposition \ref{p.monadj}, we obtain a monad $(T^\phi, \eta^\phi, \mu^\phi)$ on $\sC$, and it remains to show that this monad is precisely the original monad $T$ given on $\sC$. %
First, we have
\begin{align*}
T^\phi (f) 
	&= G_T \of F_T (f) \\
	&= G_T (\hat f, \eta(\hat f) f) \\
	&= \mu(\hat f) \cdot T (\eta(\hat f) \cdot f) \\
	&= T(f).
\end{align*}
Next,
$$\mu^\phi = G \of \epsilon \of F = \phi^{-1}(1_{\mu(y) T(f)}),$$
so
\begin{align*}
\mu^\phi(f) 
	&= G(\epsilon(F(f))) \\
	&= G(\epsilon(\hat f, \eta(\hat f) f)) \\
	&= G(\phi^{-1}(1_{G(\hat f, \eta(\hat f) f)})) \\
	&= G(\phi^{-1}(1_{T(f)})),
\end{align*}
and $\phi^{-1}(1_{T(f)}) = (f, 1_{T(f)})$ because $\phi$ is injective and $\phi(f, 1_{T(f)}) = 1_{T(f)}.$ Therefore this is
\begin{align*}
	&= G(\phi^{-1} (1_{T(f)})) \\
	&= G(f, 1_{T(f)}) \\ 
	&= \mu(f).
\end{align*}
Finally,
\begin{align*}
\eta^\phi (f)
	&= \phi(1_{F(f)}) \\
	&= \phi(1_{(\hat f, \eta(\hat f) f)}) \\
	&= \eta(f)
\end{align*}
again since $\phi$ is injective.
\Prf

\subsection{Monads and Triples}\label{ss.montriple}

For applications, e.g. \cite{MoI1, WaR1, AcAdMiVe1}, a monad is often thought of as a {\em triple} in the sense we now define. We will verify that the usual interchangability between the two notions holds, but only in a certain sense. 

\dfn\label{d.triple}
$\sC$ a generalized category. A {\em triple}, or {\em Kleisli triple}, on $\sC$ is a pair $(\eta, ()^*)$ %
where $\eta$ and $()^*$ are monotonic 
maps $\sC \to \sC$. 
For $f \in \sC$, set the abbreviation
$$T(f) := (\eta(\hat f) \cdot f)^*.$$
The mappings $\eta, ()^*$ are required to satisfy
\enu
	\item[$(0^1)$] $\widehat{f^*} = \hat f, \,\, \overline{f^*} = T(\bar f)$
	\item[$(0^2)$] $\widehat{\eta(f)} = T(f), \,\, \overline{\eta(f)} = f$
	\item[$(1)$] $\eta(f)^* = 1_{T(f)}$ \label{ax.etastar}
	\item[$(2)$] $f = f^* \cdot \eta(\bar f)$ \label{ax.stareta}
	\item[$(3)$] $(g^* \cdot f)^* = g^* \cdot f^*$
\Enu
\Dfn
Note that by axiom (1), $T$ restricts to a mapping $\Ob(\sC) \to \Ob(\sC)$. 

\exa
Triples appear in applications to programming language theory, where the mapping $T()$ is considered to send a type $A$ to its corresponding type $T(A)$ of {\em computations} of that type, in a system where programs are regarded not as pure functions, but as functions with complications such as failure to terminate, indeterminacy, continuations, or side effects. Correspondingly, $\eta: A \to T(A)$ is regarded as the inclusion of values into computations of type $A$, and $(f)^*:T(A) \to T(B)$ is thought of as the extension of $f:A \to T(B)$ to $T(A)$. 
\Exa

\prop\label{p.montriple}
In the generalized setting, a Kleisli triple gives rise to a monad, and conversely (though not on the Kleisli category $\sC_T$) assuming the following {\em hypothesis}: if $T(a) = T(b)$ then $\mu(a) = \mu(b)$. 
\Prop
\prf
Let $\sC$ be equipped with the triple $(\eta, ()^*)$. Then 
\begin{align*}
T(g \cdot f) 
	&= (\eta(\widehat{gf}) gf )^* \\
	&= (\eta(\hat g) \cdot g \cdot f)^* \\
	&= (\eta(\hat g) \cdot g \cdot f)^* \\
	&= ((\eta(\hat g) \cdot g)^* \cdot \eta(\overline{\eta(\hat g)\cdot g}) \cdot f)^* \\
	&= ((\eta(\hat g) \cdot g)^* \eta(\hat f) \cdot f)^* \\
	&= (\eta(\hat g) g)^* (\eta(\hat f) f)^* \\
	&= T(g) T(f).
\end{align*}
and
\begin{align*}
\overline{T(f)} 
	&= \overline{(\eta(\hat f) \cdot f)^*} \\
	&= T(\overline{\eta(\hat f) \cdot f}) \\
	&= T(\bar f).
\end{align*}
and
\begin{align*}
\widehat{T(f)} 
	&= \verywidehat{(\eta(\hat f) \cdot f)^*} \\
	&= \verywidehat{\eta(\hat f) \cdot f)} \\
	&= \widehat{\eta(\hat f)} = T(\hat f).
\end{align*}
and
$$ T(1_f) = (\eta(\widehat{1_f})1_f)^* = \eta(f)^* = 1_{T(f)}.$$
Now if $f \dleq g$, then 
$$T(f) = \eta(\hat f)\cdot f \dleq \eta(\hat g) \cdot g = T(g),$$
using monotonicity of $\eta$. 
So $T$ is a functor. %
Now for $f \in \sC$ define
$$\mu(f) = (1_{T(f)})^*.$$ %
To show that $(T,\eta,\mu)$ is a monad on $\sC$, it remains to check that $\eta$ and $\mu$ are natural transformations, and axioms (\ref{ax.monmagic1}) and (\ref{ax.monmagic2}) of Definition \ref{d.monad}. %
We leave this to the reader. %
We use the monotonicity of $()^*$ to prove the monotonicity of $\mu$. %


Now, %
let $(T,\eta, \mu)$ be a monad on $\sC$. %
We define a very abstract generalized category that contains some things that can be easily constructed, and some things that cannot be. Let
$$\sC'_T = \set{(f,\vec y) \mid \vec y = (y_0, y_1, \dots), \widehat{s^i(f)} = T(y_i) }.$$
Put
\begin{equation}
\begin{cases}
	\, \widehat{(f, \vec y)} = (T(y_0), \hat{y_0}, \hat{\bar{y_0}}, \hat{\bar{\bar{y_0}}}, \dots) & \\
	\, \overline{(f, \vec y)} = (\bar f, y_1, y_2, \dots) & \\
	\, (g, \vec z) \cdot (f, \vec y) = (g \cdot f, z_0, y_1, \dots) & 
\end{cases}
\end{equation}
and
\begin{equation}
\begin{cases}
	\, \eta(f, \vec y) = (\eta(f), f, \vec y) & \\
	\, (f, \vec y)^* = (\mu(y_0)T(f), y_0, T(y_1), T(y_2), \dots). & 
\end{cases}
\end{equation}
Let $(f, \vec y) \dleq (g, \vec z)$ if $f \dleq g$ and for all $i \geq 0$, $y_i \dleq z_i$, so $\sC'_T$ is a poset. 
Composition is associative since
\begin{align*}
(h, \vec v) \cdot ((g, \vec z) \cdot (f, \vec y))
	&= (h,\vec v) \cdot (g \cdot f, z_0, y_1, \dots) \\
	&= (h \cdot g \cdot f, v_0, y_1, \dots) \\
	&= (h \cdot g, v_0, z_1, \dots) \cdot (f, \vec y) \\
	&= ((h, \vec v) \cdot (g, \vec z)) \cdot (f, \vec y) 
\end{align*}
and indeed $\overline{(g,\vec z) \cdot (f, \vec y)} = (\bar f, y_1, \dots) = \overline{(f,\vec y)},$ and $\widebar{(g,\vec z) \cdot (f, \vec y)} = \widehat{(g,\vec z)}.$ %
If $\overline{(g,\vec z)} \dleq \widehat{(f,\vec y)},$ then 
$$(\bar g, y_1, \dots) \dleq (T(y_0), \hat{y_0}, \hat{\bar{y_0}}, \dots),$$
so $\bar g \dleq \hat f$, and hence $(g,\vec z) \cdot (f, \vec y)$ is defined. (The converse may be imposed in the definition of $\sC'_T$.) %
Axioms (\ref{ax.gencat-order1}) and (\ref{ax.gencat-order2}) of Definition \ref{d.gencat} are trivial. %

It can be checked that objects of $\sC'_T$ are in one-to-one correspondence with objects of $\sC$ in the image of $T$. These in turn have well-behaved identities, so axiom (\ref{ax.gencat-object-id}) of Definition \ref{d.gencat} is satisfied. %
So $\sC'_T$ is a generalized category. %
It can also be easily verified by the reader that an element $(f, \vec y)$ of $\sC'_T$ has an identity iff it is of the form $(T(u), \hat u, \hat{\bar{u}}, \dots)$ for some $u \in \sC$. In other words, there is a one-to-one correspondence between elements of $\sC$ possessing an identity and elements of $\sC$ in the image of $T$.\footnote{%
In fact, the ``ingoing'' hom set of an element {\em not} of this form is empty. Therefore these elements may be thought of as copresheaves over a base.} %

Now, we check that we have a triple on $\sC'_T$. First, we check that 
$$T(f,\vec y) = (T(f), T(y_0), T(y_1), \dots)$$
Next,
\begin{align*}
\widehat{(f,\vec y)^*} 
	&= \verywidehat{(\mu(y_0)T(f), y_0, T(\dots))} \\
	&= (T(y_0), \hat{y_0}, \hat{\bar{y_0}}, \dots) \\
	&= \widehat{(f, \vec y )}
\end{align*}
and
\begin{align*}
\overline{(f, \vec y)^*} 
	&= \overline{(\mu(y_0)T(f), y_0, T(\dots))} \\
	&= (\overline{\mu(y_0)T(f)}, T(y_1), T(y_2), \dots) \\
	&= (T(\bar f), T(y_1), \dots) \\
	&= T(\bar f, y_1, \dots) \\
	&= T(\overline{(f, \vec y)}) 
\end{align*} 
and 
\begin{align*}
\widehat{\eta(f, \vec y)} 
	&= \verywidehat{(\eta(f), f, \vec y)} \\
	&= (T(f), \hat f, \hat{\bar{f}}, \hat{\bar{\bar{f}}}, \dots) \\
	&= (T(f), T(y_0), T(y_1), \dots) \\
	&= T(f, \vec y)
\end{align*}
and similarly, $\overline{\eta(f,\vec y)} = (f, \vec y)$. Next, we have
\begin{align*}
\eta(f, \vec y)
	&= (\eta(f), f, \vec y)^* \\
	&= (\mu(f)T(\eta(f)), f, T(\vec y)) \\
	&= (1_{T(f)} , f, T(y_0), T(y_2), \dots) \\
	&= (1_{T(f)}, f, \hat f, \hat{\bar{f}}, \dots) \\
	&= 1_{(T(f), \hat f, \hat{\bar{f}}, \dots)} \\
	&= 1_{T(f, y_0, y_1, \dots)}.
\end{align*}
and
\begin{align*}
(f,\vec y)^* \cdot \eta(\overline{(f, \vec y)})
	&= (\mu(y_0)T(f), y_0, T(y_1), T(y_2), \dots) \cdot (\eta(\bar f), \bar f, y_1, y_2, \dots) \\
	&= (\mu(y_0) \cdot T(f) \cdot \eta(\bar f), y_0, y_1, \dots) \\
	&= (f, \vec y).
\end{align*}
and finally
\begin{align*}
((g, \vec z)^* \cdot (f, \vec y))^*
	&= ((\mu(z_0)T(g), z_0, T(z_1), T(z_2), \dots) \cdot (f, y_0, y_1, \dots))^* \\
	&= ((\mu(z_0)T(g)f, z_0, y_1, \dots))^* \\
	&= (\mu(z_0) \cdot T(\mu(z_0)\cdot T(g) \cdot f), z_0, T(y_1), T(y_2), \dots) \\
	&= (\mu(z_0) \cdot T(\mu(z_0) \cdot T^2 (g) \cdot T(f), \text{same}) \\
	&= (\mu(z_0) \cdot \mu(T(z_0)) \cdot T^2 (g) \cdot T(f), \text{same}), \text{ and $T(z_0) = \hat g$} \\
	&= (\mu(z_0) \cdot T(g) \cdot \mu(\bar g) \cdot T(f), z_0, T(y_1), T(y_2), \dots).
\end{align*}
Now	applying the hypothesis in the statement of the Proposition, we conclude that $\mu(\bar g) = \mu(y_0)$. Hence this expression is $(g,\vec z)^*(f,\vec y)^*$, and $(\eta, ()^*)$ is a triple on $\sC'_T$, as desired. 
\Prf

Can the hypothesis be removed from Proposition \ref{p.montriple}? %

Is there also a triple defined on $\sC_T$? %

What is $(M(\sC'_T))'_T$, where $M()$ is the construction of the first half of the proof?

It can be verified that, via the proof of \cite{MacCW} with few changes, there is a unique comparison functor from $\sC_T$ to $\sD$, given any adjunction 
\begin{tikzcd}
	\sC \arrow[r, "F", shift left] \arrow[r, leftarrow, shift right, "G"'] & \sD \\
\end{tikzcd}
that induces the monad $T$.

\subsection{Eilenberg-Moore Algebras}
Once we have investigated Kleisli categories, it is natural to ask what occurs when we look at Eilenberg-Moore algebras (as well as ordinary algebras and coalgebras) in the generalized setting. %
The results are so far inconclusive. %
Several key steps in results on Eilenberg-Moore algebras depend on the 
observation that objects in $\sC^T$ possess a shared object-morphism character. %
Thus, it might appear that the theory of algebras is in some sense already a native of the generalized setting, 
and that (interestingly) the category of algebras is a noncellular generalized category, manifesting as a sort of ``refraction,'' if you like, of the input category. 
However, the latter point of view runs into difficulties that the author could not redeem. At several points, uniqueness results apparently break down, resulting in a weaker theory than the theory of algebras based on the 1-dimensional category $\sC^T$. 

On the other hand, the 
1-dimensional definition succeeds, as usual, and takes its usual role as the terminal object in the category of those adjunctions that define $T$ in $\sC$, and supplies the usual Monadicity (Tripleability) theorems. 

\section{Conclusion and Future Work}\label{s.conc}

In this chapter we have 
presented the theory of monads as evidence of the existence of a robust mathematical framework generalizing category theory. %
We have not mentioned several topics, e.g., strengths, Linton's correspondence, and distributive laws, %
but we have given indication of some applications of this framework. 

It may be noticed that the mathematics of generalized categories is in some ways conceptually simpler, and in other ways conceptually more challenging compared to ordinary category theory. The inherent emphasis on recursion and corecursion in the framework suggests potential connections between generalized categories, iterative equations, and categorical fixpoint theory. 



\pagebreak
\singlespacing
\chapter{A Generalization of the Curry-Howard Correspondence}\label{c.th}
\doublespacing
\vspace{10ex}

\section{Introduction}\label{s.th-i}

In a series of papers \cite{LaK1c,LaK2}, Lambek developed an extension of the Curry-Howard correspondence \cite{Howard1969} to the domain of categorical logic. %
Lambek's extension has since become a cornerstone of programming language theory, particularly in the functional programming paradigm. It has also been influential in logic. %
This paper 
is devoted to a generalization of the Curry-Howard-Lambek correspondence which makes use of the tools provided by generalized categories. %
Those who agree with Philip Wadler \cite{WaR5} that, as a general rule, semantics should guide development in logic and programming language theory may take interest in this product of a generalization on the semantic side. %
Those with a pure interest in category theory might note some features of our approach, for example, we show (section \ref{s.cat}) that using the framework of generalized categories, a cartesian functor between cartesian closed categories may be ``promoted'' to a cartesian closed functor. To the best of our knowledge this construction is at least somewhat new. 

Lambek in his work makes extensive use of deductive systems \cite{LaK1c}. %
A short discussion of the intuition for this notion (which may be unfamiliar) affords the opportunity to provide some intuition for the notion of generalized category. However, the reader is free to ignore this discussion if he or she wishes; nothing in the main body of the paper depends on it. %
A deductive system is just enough machinery to allow the question: from a given point $a$ of the deductive system $\sA$, can I travel to another point $b \in \sA$ via a valid path? %
A conceptual picture of this is the following. Suppose that there is a system of goods $\sA_0$. %
The edges of $\sA$ are certificates (issued, say perhaps, by different governing bodies) that say that a good $a \in \sA_0$ may be exchanged for another good $b \in \sA_0$. %
It is accepted that a good is always exchangeable for itself. %
Now let's suppose that such certificates themselves may be exchanged, 
but that this requires that one has a higher-level certificate for this higher-level trade. %
If we imagine a certain impetus exists among those we imagine making the exchanges, we can expect that there will next arise trading for these certificates as well. %
Let us make two simple observations:
\enu
	\item The resulting deductive system is not necessarily {\em cellular}, in the sense that the economy is liberalized to the extent that certificates may be good for exchange of different {\em kinds} of goods and certificates. For example, a certificate may be for a good, in return for a certificate good for a certificate in return for a good.
	\item There need not be, in the abstract, any {\em goods} at all. The system could be one of certificates for certificates for certificates, and so on. This observation may be utilized to clean up the abstract formalism: a system with no atoms is conceptually simpler and the easiest one to work with while developing elementary principles. 
\Enu
These two observations suggest, via the intuition, a generalization of category theory that we outline in section \ref{s.gencat}. 

Some work during intermediate stages is necessary in order to accomplish our aim. %
Under the usual Curry-Howard correspondence, types are interpreted as propositions which are true only when they are inhabited by a term. 
It is based on the types-as-targets view of categorical semantics, which limits the applicability of generalized categories to type theory. %
If we consider the alternative {\em types-as-paths} view, in which a proposition depends on both a source and a target, we find a calculus that is not only amenable to the generalized setting, but also fits well with the Lambek equational theory of cartesian closed categories \cite{LaK2}. %
The types-as-path view is motivated by the notion that a type is like the blueprint of a bridge between two points, or (in the logical intuition), a {\em conjecture}. %
Using the intuition from programs, on the other hand, the type-as-path is an {\em approximation} or {\em abstraction} of any choice of concrete transformation between two different kinds of data. This supports our approach, since this is how types are often viewed in applications, see for example \cite{Pierce1}. %
The types, which we write $a \vdash b$, when viewed categorically, assume the role of exponential objects. %
We are able to give this description a precise formal treatment by combining (1) the contributions of Lambek and (2) the framework of generalized category theory. 

In section \ref{s.th}, we develop ideal cartesian closed  categories, the notion we take of cartesian closed category in the generalized setting. These come equipped with an ideal of types, in the sense discussed above. In section \ref{s.poly}, we introduce polynomial categories, by closely following Lambek \cite{LaK2}, and in section \ref{s.lam} we define a notion of generalized type theoretic system (lambda calculus) corresponding to the semantics we have introduced, and verify that the anticipated equivalence holds. %
In all that we have done we have closely followed 
the well-established work of Lambek and others. However, our work lays the foundation for many possible avenues for further development in areas such as proof theory, programming language semantics, topos theory, and homotopy type theory. We discuss some topics for future work in section \ref{s.conc}. 

\vspace{12pt}
\section[Generalized Deductive Systems and Ideal Cartesian Closed Categories]{\onehalfspacing Generalized~Deductive~Systems and
Ideal~Cartesian~Closed~Categories}\label{s.th}

\subsection{Generalized Deductive Systems and Generalized Graphs}\label{ss.gendedsys}

\dfn\label{d.gengraph}
A {\em generalized graph} is a triple $(\sA, s, t)$, where $\sA$ is a set, %
$s, t$ are maps $\sA \to \sA$. %

A {\em morphism} $\Phi: \sA \to \sB$ of generalized graphs is a mapping $\Phi$ from $\sA$ to $\sB$ such that 
\enu
	\item $\Phi(s(a)) = s(\Phi(a))$
	\item $\Phi(t(a)) = t(\Phi(a))$
\Enu
This gives a category $\bf{Graph}$ of generalized graphs.
\Dfn

An element of $\sA$ is (synonymously) an {\em edge}. %
An {\em object} in a generalized graph is an element $a \in \sA$ such that $sa = ta = a$. %
We say that a {\em subject} in a generalized graph is an element $a \in \sA$ such that there is an element $f \in \sA$ such that either $sf = a$ or $tf = a$. %
We write $\Ob(\sA)$, $\Sb(\sA)$ for the set of objects and subjects of $\sA$, respectively. %
We say that generalized graph is {\em 1-dimensional} if 
$$ss = s \text{ and } tt = t.$$
Ordinary graphs correspond bijectively with 1-dimensional generalized graphs. %

Recall that in an algebraic system $(A, f)$ in which $A$ is a carrier set where equality $(=)$ is defined and a unary operation $f$ is defined (a mapping $A \to A$), we say that $f$ is {\em substitutive} if for all $a,b \in A$
$$a = b \eimplies fa = fb.$$
(The word {\em congruence} also arises frequently in connection with this property.) %
The source and target operations in a generalized graph are not assumed to be substitutive. %
(In fact, there is no notion of equality defined in the language of generalized graphs until we come to Definition \ref{d.icat}.) %
This comes with the advantage that we can apply inductive pattern-matching in proofs about elements in a generalized graph (and we may even do so constructively, if they are finitely generated in some finite language), %
though yet another hypothesis is needed if these patterns matchings are to be exhaustive in $\sA$. (Such a hypothesis will apply to polynomials in section \ref{s.poly}.) %

\dfn\label{d.dedsys}
A {\em generalized typed deductive system} or simply a {\em generalized deductive system} is a structure 
$$(\sA, s, t, \cdot, \vdash, \sV),$$ 
where $(\sA, s, t)$ is a generalized graph, 
($\cdot$) is a partially-defined operation $\sA \times \sA \to \sA$ on $\sA$,
$(\vdash)$ is an operation $\sA \times \sA \to \sA$, 
and $\sV$ is a subset of $\sA$, satisfying 
\enu
	\item for all $a,b \in \sA$, $b \cdot a \text{ is defined iff } ta = sb$
	\item $s(ab) = s(b)$ and $t(ab) = t(a)$
	\item $s(a \vdash b) = a$ and $t(a \vdash b) = b$. \label{ax.vdashst}
	\item for every $a \in \sA$, $a \vdash a \in \sV$. \label{ax.identityvalid}
	\item if $a,b \in \sV$, and $a \cdot b$ $\downarrow$, then $a \cdot b \in \sV$. \label{ax.compositionvalid}
	\item for every $a,b \in \sA$, if there exists $u \in \sV$ with $\bar u = a$ and $\hat u = b$, then $a \vdash b \in \sV$. \label{ax.inhabitance}
\Enu
A morphism of generalized typed deductive systems $\phi: \sA \to \sB$ is a morphism of generalized graphs satisfying
\enu
	\item $\phi(a \cdot b) = \phi(a) \cdot \phi(b),$
	\item $\phi(a \vdash b) = \phi(a) \vdash \phi(b),$
	\item if $a \in \sV$, then $\phi(a) \in \sV$.
\Enu
This gives a category ${\bf DedSys}$ of generalized (typed) deductive systems.
\Dfn

Since we can now compose edges, we shall refer elements $a$ of a deductive system $\sA$ as edges or {\em paths} (there is no actual distinction between the two terms, except in case products in $(\cdot)$ are freely generated on a basis in $\sA$.) %
The elements of $\sV$ may be thought of as {\em valid paths} of $\sA$. 
In the set of edges going from $a$ to $b$, the unique edge $a \vdash b$ is called the {\em type} with source $a$ and target $b$. %
We may use the notation $a \dashv b$ interchangeably to denote $b \vdash a$, thus $a\dashv b \threeline b \vdash a$. 
Finally, when using axiom (\ref{ax.inhabitance}) we call $u$ a {\em witness} and say that the type $a \vdash b$ is {\em inhabited} if there is found such a $u$. %
We may write $1_a$ in place of $a \vdash a$. %


In our work it is possible to ignore the role of $\sV$, but its presence suggests generalizations of the calculus, for example $\sV$ 
might be useful in a model of concurrency, or be impacted by modal operators.

\dfn\label{d.pidedsys}
An {\em positive intuitionistic generalized deductive system} is a generalized deductive system 
$$(\sA, s, t, \cdot, \vdash, \sV)$$ 
equipped with the additional structure
$$(\trut, \wedge, \langle,\rangle, ()^*)$$
consisting of:
\enu
	\item A distinguished element $\trut \in \sA$,
	\item A mapping $\wedge: \sA \times \sA \to \sA$,
	\item A partially defined mapping $\langle,\rangle : \sA \times \sA \to \sA$,
	\item A partially defined mapping $()^* : \sA \to \sA$
\Enu
subject to the following axioms,
\enu
	\item $\langle a,b \rangle$ is defined if and only if the source of $a$ and $b$ are identical.
	\item $a^*$ is defined if and only if the source of $a$ is the wedge of two subjects in $\sA$.
\Enu
the following source and target conditions:
\enu
	\item $\hat\trut = \bar\trut = \trut$, 
	\item $s(a \wedge b) = s(a) \wedge s(b)$ and $t(a \wedge b) = t(a) \wedge t(b)$, 
	\item $s(\langle a , b \rangle) = s(a)$ and $t(\langle a, b \rangle) = t(a) \wedge t(b)$
	\item $s(a^*) = proj_1(s(a))$ and $t(a^*) = proj_2 (s(a)) \vdash t(a)$, where $proj_1$ and $proj_2$ are the projections on wedge ($\wedge$) products.
\Enu
and the following {\em rules}, or {\em validities:} For all $a,b \in \sV$,
\enu
	\item $\trut \in \sV$,
	\item $a \wedge b$ is valid if $a$ and $b$ are valid,
	\item $\langle a, b \rangle$ is valid if $a$ and $b$ are valid,
	\item $a^*$ is valid if $a$ is valid,
	\item for every pair of subjects $a,b \in \sA$, the following types are valid: \label{ru.validatoms}
		\enu
			\item $\trut \vdash a$
			\item $a \wedge b \vdash a$ 
			\item $a \wedge b \vdash b$ 
			\item $(a \vdash b) \wedge a \vdash b$ 
		\Enu
\Enu
A morphism $f: \sA \to \sB$ of positive intuitionistic generalized deductive systems is a morphism of generalized deductive systems satisfying 
\enu
	\item $F(\trut) = \trut,$
	\item $F(a \wedge b) = F(a) \wedge F(b),$
	\item $F(\langle f, g \rangle ) = \langle F(f), F(g) \rangle,$
	\item $F(f^*) = F(f)^*.$
\Enu
This gives a category ${\bf p.i.DedSys}$ of positive intuitionistic generalized deductive systems. 
\Dfn

In order to form complex expressions out of simple ones, it is convenient to have names for individual elements of $\sA$. %
For example, we choose (applying rule \ref{ru.validatoms}) valid elements of $\sA$
	$$ter_a : a \to \trut$$
	$$\pi_{a,b} : a \wedge b \to a$$
	$$\pi_{a,b} : a \wedge b \to b$$
	$$\epsilon_{a,b} : (a \vdash b) \wedge a \to b$$
Note that these elements may themselves be types, even though we usually think of the types as valid due to the existence of a witness and use of axiom \ref{ax.inhabitance} of Definition \ref{d.dedsys}. %
By {\em term} (or {\em global element}) of a deductive system we refer to any element of a deductive system whose source is $\trut$. 

One may use the deductive system to show the validity of Heyting's axioms for intuitionistic logic (those that do not contain the $\lor$ and $\bot$ connectives), showing that any type that may be interpreted as a valid proposition of intuitionistic logic has a witness. The following types, for example, are inhabited.
\enu
	\item $a \vdash a \wedge a$
	\item $a \vdash a \wedge \trut$
	\item $((a \wedge b) \vdash c) \vdash (a \vdash (b \vdash c))$
	\item $(a \vdash (b \wedge c)) \dashv \vdash ((a \vdash b) \wedge (a \vdash c))$
\Enu
where $\dashv \vdash$ denotes that the type is bi-inhabited (or there is a valid path going in either direction). %

\subsection{Categories Equationally Defined}\label{s.cat}
Lambek \cite{LaK1c,LaSc1} observed that categories are obtained from deductive systems via a set of equational axioms. %
In this section we will develop Lambek's formalization in the setting of generalized categories. %
It is clear that any ordinary category (or generalized category) can be made into a ``typed deductive'' category. %
Simply take all arrows to be valid and introduce $\vdash$ as a free operation 
Observe that if composition $(\cdot)$ is viewed as multiplication and $\wedge$ is viewed as an additive product on the subjects of $\sC$, the set of elements of the form $a \vdash b$ behaves like a (ring-theoretic) ideal in the category. Thus if we are thinking of a category, we may think next of introducing an ``ideal of types'' to the category. %
This demands we introduce a further technicality, a set of constants. %

\dfn\label{d.icat}
A {\em ideal category} or {\em ideal generalized category} is a structure $(\sC, \vdash, \sV)$ consisting of a generalized category $\sC$ (section \ref{s.gencat}), %
a distinguished subset of elements $\sV \lies \sC$, 
a distinguished subset of elements $\sK \lies \sC$,
and an operation $\vdash: \sC \times \sC \to \sC$ such that 
\enu
	\item $\verywidehat{f \vdash g} = f,$
	\item $\overline{f \vdash g} = g.$
	\item $f \cdot (\bar f \vdash g) = \hat f \vdash g,$ unless $\bar f = g$, in which case $f \cdot (\bar f \vdash g) = f$, or unless $f \in \sK$ or $\bar f \vdash g \in \sK$. \label{ax.vdash1}
	\item $(f \vdash \hat g) \cdot g = f \vdash \bar g,$ unless $f = \hat g$, in which case $(f \vdash \hat g) \cdot g = g$, \label{ax.vdash2}
	\item if $g \cdot f$ $\downarrow$, and $g,f \in \sV$, then $g \cdot f \in \sV$,
	\item $f \vdash f \in \sV$ for all $f \in \sC$,
	\item (witnesses) $u \in \sV$ implies $\hat u \vdash \bar u \in \sV$.
	\item $\vdash$ is substitutive (Section \ref{ss.gendedsys}) in both arguments. 
\Enu
A functor $F: \sC \to \sD$ between generalized ideal categories is an ordinary functor (section \ref{s.gencat}) which preserves validity and $\vdash$:
\enu
	\item $f \in \sV_\sC \eimplies F(f) \in \sV_\sD,$
	\item $F(f \vdash g) = F(f) \vdash F(g)$. \label{ax.deductivefunctor2}
\Enu
This defines a category ${\bf IdealCat}$. 
\Dfn

By axioms \ref{ax.vdash1} and \ref{ax.vdash2}, for $f \in \sC$, $f \vdash f$ is the identity of $f$, which may be denoted $1_f$. %
In particular, all elements (including identities) of an ideal category have identities. 
The identities, types, and constants figuring here will arise again in Section \ref{s.lam}, where we encounter the symbols $\settl{x}{x}$ and $\settl{x}{y}$. %



\dfn\label{d.iccc}
An {\em ideal cartesian closed category} is a ideal category with identities $\sC$ that is equipped with a structure
$$(\trut, \wedge, \langle,\rangle, ()^*)$$
where
\enu
	\item $\trut$ is a distinguished valid element in $\sC$,
	\item $\wedge$ is an operation $\sC \times \sC \to \sC$,
	\item $\langle, \rangle$ is a partially defined operation $\sC \times \sC \to \sC$
	\item $()^*$ is a partially defined operation $\sC \to \sC$
\Enu
which satisfies the conditions:
\enu
	\item $\trut \in \sK$, and $\sK$ is closed under $\wedge, \langle,\rangle, $ and $()^*$,
	\item the structure
		$$(s,t,\sV,\vdash, \cdot, \trut, \wedge, \langle,\rangle, ()^*)$$
		defines a positive intuitionistic deductive system on $\sC$.
	\item for all $a \in \sC$, if $f: a \to \trut$ then $f = (a \vdash \trut)$. \label{ax.terminalar}
	\item For every pair $(a,b)$ of subjects of $\sC$, there exists a good pair $(\pi, \pi')$ for $(a,b)$. 
	\item For every good pair $(\pi, \pi')$ for any pair of subjects $(a,b)$, there is a good evaluation $\epsilon = \epsilon_{\pi, \pi'}$ for $(\pi, \pi')$. 
\Enu
Here, if $(a,b)$ is a pair of subjects of $\sC$, then a pair $(\pi, \pi')$ of elements of $\sC$ are a {\em good pair for $(a,b)$} if
			\enu
				\item $\pi$ and $\pi'$ are valid,
				\item $\pi: a \wedge b \to a$, and $\pi': a \wedge b \to b$,
				\item if $\pi \langle f, g \rangle$ $\downarrow$ then $\pi \langle f, g \rangle = f,$
				\item if $\pi' \langle f, g \rangle$ $\downarrow$ then $\pi' \langle f, g \rangle = g,$
				\item if $\langle \pi f, \pi' f \rangle$ $\downarrow$ then $\langle \pi f, \pi' f \rangle = f$, 
				\item if $f \cdot \pi$ and $g \cdot \pi'$ $\downarrow$, then $\langle f \cdot \pi, g \cdot \pi' \rangle = f \wedge g$. 
			\Enu
and a {\em good evaluation} for a good pair $(\pi, \pi')$ for a pair of subjects $(a,b)$ is an element $\epsilon = \epsilon_{\pi, \pi'}$ of $\sC$ that satisfies, for every $c \in \sC$ and every good pair $(\pi_{c,a}, \pi'_{c,a})$ for $(c,a)$,
			\enu
				\item $\epsilon$ is valid,
				\item $\epsilon: (a \vdash b) \wedge a \to b$,
				\item if $\epsilon \cdot \langle f^* \cdot \pi_{c,a}, \pi'_{c,a} \rangle$ $\downarrow$ then $\epsilon \cdot \langle f^* \cdot \pi_{c,a}, \pi'_{c,a} \rangle = f,$ 
				\item if $(\epsilon \cdot \langle f \cdot \pi_{c,a}, \pi'_{c,a} \rangle )^*$ $\downarrow$ then $(\epsilon \cdot \langle f \cdot \pi_{c,a}, \pi'_{c,a} \rangle )^* = f$. 
			\Enu
A morphism $F: \sC \to \sD$ between ideal cartesian closed categories $\sC$ and $\sD$ is a functor of ideal categories satisfying
\enu
	\item $F(\trut) = \trut$,
	\item $F(a \wedge b) = F(a) \wedge F(b),$
	\item $F(\langle a , b \rangle) = \langle F(a), F(b) \rangle,$
	\item $F$ sends a good pair in $\sC$ to a good pair in $\sD$. 
\Enu
Thus we have a category ${\bf ICCC}$ of ideal cartesian closed categories.
\Dfn

Axiom \ref{ax.terminalar} is relevant when the possibility exists that the element $a \vdash \trut$ might be a constant. 
%
We continue to use the notation of deductive systems in a category $\sC$. Note that many authors write $\times$ for the binary product, which we continue to denote $\wedge$, and $1$ for the terminal object, which we continue to denote $\trut$. This seems appropriate as we will never stray far from the point of view provided by deductive systems and the lambda calculus. 

Note that morphisms of ideal cartesian closed categories are stronger maps than ordinary functors between categories that happen to be cartesian closed. For ordinary categories, these functors are sometimes called {\em cartesian functors}. %
It is easy to see that a good pair $(\pi, \pi')$ for a pair of subjects $(a,b)$ is unique if it exists. Hence a good evaluation $\epsilon \threelines \epsilon_{\pi, \pi'}$ depends only on $(a,b)$ and may be denoted $\epsilon_{a,b}$. 
Similarly, we often write $\pi \threelines \pi_{a,b}$ and $\pi' \threelines \pi'_{a,b}$. It follows that $F(\pi_{a,b}) = \pi_{F(a), F(b)}$, and similarly for $\pi'$.

\prop\label{p.basiciccc}
The following hold in ideal cartesian closed categories:
\enu
	\item $\langle f, g \rangle \cdot h = \langle f \cdot h , g \cdot h \rangle$
	\item $1_a \wedge 1_b = 1_{a \wedge b}$
	\item $(f \wedge g) \cdot (f' \wedge g') = (f \cdot f') \wedge (g \cdot g')$
	\item $\epsilon_{a,b}^* = 1_{a \vdash b}$
	\item $f^* \cdot g$ $\downarrow$ implies $f^* \cdot g = (f \cdot \langle g \cdot \pi, \pi' \rangle)^*$, where $(\pi, \pi')$ is the obvious good pair.
\Enu
\Prop

\prop
A morphism $F: \sC \to \sD$ between ideal cartesian closed categories preserves the evaluation $\epsilon$ and adjoint operation $()^*$. %
\Prop
\prf
By functoriality, we have $F((f: a \wedge b \to c)^*) = F(f^*) : F(a) \to F(b \vdash c) = F(f^*) : F(a) \to (F(b) \vdash F(c))$. But this latter expression is $F(f)^*$, so $F(f^*) = F(f)^*$. %
It follows that evaluations $\epsilon$ are also preserved. Indeed, if $(a,b)$ are chosen and $(\pi, \pi')$ is a good pair for $(a,b)$, then choose a good evaluation $\epsilon = \epsilon_{a,b}$ for $(\pi, \pi')$. %
Then 
$$F(\epsilon)^* = F(\epsilon^*) = F(1_{a \vdash b}) = 1_{F(a) \vdash F(b)}.$$
Hence 
\begin{align*}
\epsilon_{F(a), F(b)} 	&= \epsilon_{F(a), F(b)} \cdot \langle 1_{F(a) \vdash F(b)} \pi_{F(a) \vdash F(b) ,a}, \pi'_{F(a) \vdash F(b),a} \rangle  \\
				&= \epsilon_{F(a), F(b)} \cdot \langle F(\epsilon_{a,b})^* \pi_{c,a}, \pi'_{c,a}\rangle  \\
				&= F(\epsilon_{a,b}),
\end{align*}
by the good evaluation properties of $\epsilon_{a,b}$. 
\Prf


Next we present a few ways to produce ideal cartesian closed categories:

\prop\label{p.cattoicat}
There is an (in general, nonconstructive) functor from the category ${\bf CCC}$ of cartesian closed categories to the category {\bf ICCC}. 
\Prop
\prf
Let $F: \sC \to \sD$ be a functor in the category of cartesian closed categories (of the ordinary sort). We carry out the following construction on both $\sC$ and $\sD$; first take $\sC$. Take any new pair of identifiers $\vdash$ and $\wedge$. For each object $X$ of $\sC$, form, via recursion, the collections of triples
$$\sA_X = \set{(Y_1, Z_1, \vdash) \mid \text{ there exists } Y, Z \in \Ob(\sC) \text{ such that } X = Z^Y \eand Y_1 \in \sC_Y, Z_1 \in \sC_Z}$$
$$\sB_X = \set{(Y_1, Z_1, \wedge) \mid \text{ there exists } Y, Z \in \Ob(\sC) \text{ such that } X = Y \wedge Z \eand Y_1 \in \sC_Y, Z_1 \in \sC_Z}$$
$$\sC_X = \sA_X \cup \sB_X.$$
We take
$$\Ob(\tilde{\sC}) = \Ob(\sC) \cup \bigcup_{X \in \Ob(\sC)} \sC_X,$$
and for each $V \in \Ob(\tilde{\sC})$ we assume given from the construction of the $\sC_X$'s a function $\deflate(V)$ defined by
$$\deflate(V) = \begin{cases} V, & \eif V \in \sC, \\ Z^Y, &\eif V \in \sA_X \text{ for some $X$}, \\ Y \wedge Z, & \eif V \in \sB_X \text{ for some $X$}. \end{cases}$$
For every $U,V \in \tilde \sC$, define
$$\hom(U,V) := \hom(\deflate(U), \deflate(V)),$$
with composition and identities defined in the obvious way, in particular 
$$\deflate(f \cdot g) := \deflate(f) \cdot \deflate(g),$$
where $\deflate(f)$ for a morphism $f$ is defined in the obvious way analogous to $\deflate()$ on objects. 
The reader can now check that the symbols in Definition \ref{d.iccc} may be introduced and the axioms verified, and that we may extend $F$ to a functor $\tilde F: \tilde \sC \to \tilde \sD$ that satisfies the conditions of Definitions \ref{d.icat} and \ref{d.iccc}. 
\Prf

Another result that gives examples of ideal cartesian closed categories is:

\prop
Let $\sE$ be a generalized category of generalized presheaves over a generalized category $\sC$. Then $\sE$ is an ideal cartesian closed category. 
\Prop
\prf 
See Chapter \ref{c.it}. 
\Prf

The adjunction that holds in a cartesian closed category, because the mappings $- \times X$ and $-^X$ are no longer functors in the generalized setting. However, we do have:

\prop
\enu
	\item there is a bijection  
		$$\hom(c \wedge b, a) \bij \hom(c, a^b)$$
	\item there is a bijection
		$$\hom(a,b) \bij \hom(\trut, b^a)$$
\Enu
Let $a \iso b$ denote that there exists a pair of elements $f:a \to b$ and $g:b \to a$ such that $fg = 1_b$ and $gf = 1_a$. Then in an ideal cartesian closed category
\enu
	\item $(a \wedge b) \vdash c \iso (a \vdash b) \vdash c$
	\item $a \vdash (b \wedge c) \iso (a \vdash b) \wedge (a \vdash c)$
\Enu
\Prop
\prf
See \cite{LaSc1}. 
\Prf

Given $f:a \to b$ we write 
$$\name f $$ 
for the induced term $1 \to a \vdash b$, called (Lawvere's terminology) the {\em name} of $f$. 
%

Finally, we relate deductive systems to categories as follows:

\prop
Every deductive system $\sA$ on which there is defined an equivalence relation denoted $=$, and a distinguished subset $\sK$ of constants in $\sA$, with respect to which the following statements are satisfied:
\enu
	\item $f \cdot (\bar f \vdash g) = \hat f \vdash g,$ unless $\bar f = g$, in which case $f \cdot (\bar f \vdash g) = f$, unless $f$ is constant or $\bar f \vdash g$ is constant,
	\item $(f \vdash \hat g) \cdot g = f \vdash \bar g,$ unless $f = \hat g$, in which case $(f \vdash \hat g) \cdot g = g$, 
	\item $(hg)f = h(gf)$ for all composable $f,g,h \in \sA$,
	\item $a = b$ implies $s(a) = s(b)$,
	\item $a = b$ implies $t(a) = t(b)$,
	\item $a = b$ implies $ca = cb$ and $ac = bc$, for all composable $c$,
	\item $a = b$ implies $a \vdash c = b \vdash c$ and $c \vdash a = c \vdash b$, for all $c$ in $\sA$,
\Enu
is an ideal generalized category (in particular, a generalized category), taking $\sV$ to be the valid paths in $\sA$. 
\Prop
\prf
We check the axioms of Definition \ref{d.gencat} and see that they may be verified using axioms and rules of Definitions \ref{d.dedsys} and \ref{d.icat}. %
\Prf

The notion of a cartesian closed category cannot be extended to the generalized setting: the mapping $X \mapsto X \times Y$ is a functor only when $X$ is an object. %
Our approach is to allow 
the mapping on the other side, $Z \mapsto Z^Y$, to fail to be a functor as well. %
This is possible thanks to Lambek's formalization: %
We are able, by following Lambek, to derive a calculus of cartesian closed categories in the generalized setting, in spite of the weaker underlying structure. %




\section{Polynomials and Lambda-Calculi}\label{s.poly}

Adding variables to a deductive system with a positive intuitionistic structure reduces, by the Deduction Theorem (Theorem \ref{t.dedt}), validity of all paths to the validity of paths from a terminal object. Therefore the focus shifts from the space to the {\em polynomials over the space}, in the sense we now define. %

\subsection{Polynomials Systems and Polynomial Categories}\label{ss.poly}
The notion of indeterminate may be applied in this setting just as it may be applied in the setting of groups, rings, and fields. However, 
we must assign a source and target to each new indeterminate. It is convenient to let the source of every indeterminate be $1$, the (fixed choice of) terminal object. This does not mean we cannot substitute a variable with a different source for the indeterminate---substitution of, say, $a$ for $x$ in $\phi(x)$ is allowed whenever $x$ and $a$ have the same target; the source of $a$ is irrelevant. %
In this sense, it is more correct (but less convenient) to say that an indeterminate simply does not have a source. %
We denote an indeterminate over a deduction system $\sA$ by symbols $x, y, z, $ etc. 
For now, we require that the target of $x, y, \dots$ is in $\sA$. (In particular, it cannot itself be a polynomial). 
A more general system might allow indeterminates over polynomials and make use of the notion of {\em telescope} \cite{deBruijn1}, but we will have no need for this added generality.

\dfn\label{d.poly}
Let $\sA$ be a positive intuitionistic deductive system. %
Let $x$ be an indeterminate with target $\hat x$ in $\sA$.
We write $\sA[x]$ for the positive intuitionistic deduction system freely 
generated on the set $\sA \cup \set{x}$. %
This means that
\enu
	\item Operations on $\sA$ of Definition \ref{d.pidedsys} are extended from $\sA$ to $\sA[x]$ by free generation on expressions $\phi$ containing any instance of $x$: 
		$$\phi ::= f \,\mid\, x \,\mid\, 
		\phi \vdash \phi \,\mid\, \phi \wedge \phi \,\mid\, \langle \phi, \phi \rangle \,\mid\, \phi^* $$
	where $f$ can be any element of $\sA$, and $x$ is any indeterminate. Expressions so generated that do not contain any instance of $x$ are thrown out, and the set of all elements of $\sA$ is then added back in. 
	\item The valid elements of $\sA[x]$ are $x$, those of $\sA$, and those generated from $x$ and those of $\sA$ using the validities of Definition \ref{d.pidedsys}.
\Enu
There is an obvious embedding of $\sA$ in $\sA[x]$, via which we will usually view $\sA$ as a subset of  $\sA[x]$. 
\Dfn

We call elements of $\sA[x]$ synonymously {\em polynomials over $\sA$}. %
We write $\phi, \psi, \dots$ to denote polynomials in $\sA[x]$. %
We do not normally write the variable $x$ as in $\phi(x),$ etc. as many authors do, but this should not lead to any confusion as long as it is understood what may depend on $x$.
When we iterate to form $\sA[x][y]$, etc., we again require that the source and target of indeterminates be in $\sA$. Given indeterminates $x_1, x_2, \dots, x_n$, we denote by $\sA[x_1, \dots, x_n]$ or $\sA[\vec x]$ the iterated construction $(\dots((\sA[x_1])[x_2]) \dots [x_n])$.


We could define a ``proof'' to be a valid path from the terminal object $\trut$ in a positive intuitionistic deductive system (say). Then we could ask what structure might allow us to ``discharge'' assumptions, as is done in natural deduction systems (see for example \cite{TrSc1}). 
To refine the question, one may consider a proof $\phi$ of $f \in \sA[x]$, for $f \in \sA$. 
This would be a path through the deductive system that is allowed to ``use'' the ``assumption'' $x$. %
In logic, the following result is, by long tradition, known as the Deduction Theorem. It is interpreted as an introduction rule when the construction of polynomials is interpreted as establishing a context. Note that polynomials do not necessarily have an element of $\sA$ as source and target, so the quantifiers on $a$ and $b$ are a significant part of the statement. 

\thm\label{t.dedt}
Let $\sA$ be a positive intuitionistic deductive system. 
Then for all $a,b \in \sA$, $a \vdash b$ is valid in $\sA[x]$ if and only if $\hat x \wedge a \vdash b$ is valid in $\sA$. 
\Thm
\prf
The proof is just as in \cite{LaSc1}, except that we must add clauses for the operations $\wedge$ and $\vdash$. 
Note that several steps depend on the existence of identities on the subjects of $\sA$, as assumed in definition \ref{d.pidedsys}. 
First, 
let $f$ be a valid path from $\hat x \wedge a$ to $b$ in $\sA$. Then since $\phi := \langle (x \cdot (a \vdash \trut) , 1_a \rangle $ is a valid path from $a$ to $\hat x \wedge a$ in $\sA[x]$, we obtain a witness $f \cdot \phi$ of the type $a \vdash b$ in $\sA[x]$, as desired. 

Now suppose $\phi$ is a valid path from $a$ to $b$ in $\sA[x]$. Suppose that for all polynomials in $x$ $\phi_<$ of length strictly less than $\phi$, there is a witness of $\hat x \wedge \overline{\phi_<} \vdash \widehat{\phi_<}$, denoted 
$$\kappa_x (\phi_<).$$
Now we proceed by cases:
\enu
	\item if $\phi \in \sA$, then $\phi \cdot \pi'_{\hat x, a}$ validates $\hat x \wedge a \vdash b$. \label{case.constant}
	\item if $\phi = x$, then $\pi_{\hat x, a}$ validates $\hat x \wedge a \vdash b$. 
	\item if $\phi = \psi \vdash \chi$ for some $\psi, \chi \in \sA[x]$, then 
			$a$ is identical to $\psi$ and $b$ is identical to $\chi$, hence this case reduces to case (\ref{case.constant}).
	\item if $\phi = \psi \cdot \chi$ for some $\psi, \chi \in \sA[x]$, then 
			$$\kappa_x \psi \cdot \langle \pi_{\hat x, a}, \kappa_x \chi \rangle$$
			is the desired witness. ($\chi \cdot \kappa_x \psi$ doesn't work, because $x$ is still not eliminated.)	
	\item if $\phi = \psi \wedge \chi$ for some $\psi, \chi \in \sA[x]$, then 
			$$\langle \kappa_x (\psi) \cdot \pi_{\bar{\psi}, \bar{\chi}}, \kappa_x (\chi) \cdot \pi'_{\bar{\psi}, \bar{\chi}} \rangle $$
			is the desired witness. (The alternative witness $\kappa_x \psi \wedge \chi_x \cdot \lambda$, where $\lambda$ is a munging factor, gives a definition of $\kappa_x$ under which one does not prove Theorem \ref{t.fCompT}.)
	\item if $\phi = \langle \psi, \chi \rangle$ for some $\psi, \chi \in \sA[x]$, then
			$$\langle \kappa_x \psi, \kappa_x \chi \rangle$$
			is the desired witness.
	\item if $\phi = \psi^*$ for some $\psi \in \sA[x]$, then
			$$(\kappa_x (\psi )\cdot \alpha)^*$$
			is the desired witness, where $\alpha$ is the associator.
\Enu
Proceeding by induction on the length of polynomials $\phi$ in $\sA[x]$ if necessary, we obtain in each case the desired witness of $\hat x \wedge a \vdash b$. 
\Prf

We denote the witness of $\hat x \wedge a \vdash b$ derived by pattern matching on $\phi: \bar \phi \to \hat \phi$ in the second half of the preceding proof by
$$\kappa_x (\phi): \hat x \wedge \bar \phi \to \hat \phi.$$
Now we pass from deductive systems to (ideal) categories. When we do so, it is necessary to ensure that the polynomial system over an indeterminate remains in our category. Hence we fix the following definition: 

\dfn\label{d.polyicat}
Let $\sC$ be an ideal cartesian closed category. %
Let $x$ be an indeterminate in $\sC$. 
To define the symbol
$$\sC(x),$$ 
observe that $\sC$ is equipped with the structure 
$$(s, t, \cdot, \vdash, \sI, \sV)$$
of a positive intuitionistic deductive system, when regarded as a generalized graph. 
Take $\sK_{\sC(x)}$ to be the set of constant polynomials.\footnote{This definition restricts behavior of terminal arrows $\phi \vdash \psi$ for polynomials $\phi$ and $\psi$, but it will not make a difference for our purposes.} 
Now take the polynomial system $\sC[x]$ of Definition \ref{d.poly}, %
and then take the smallest equivalence relation $=_x$ of paths in $\sC[x]$ satisfying the conditions:
\enu
	\item If $f = g$ in $\sC$, then $f =_x g$ in $\sC(x)$,
	\item $(\phi \vdash \hat \psi ) \cdot \psi =_x (\phi \vdash \bar \phi)$ unless $\phi =_x \bar \psi$, in which case $(\phi \vdash \hat \psi) \cdot \psi = \psi$,
	\item $\psi \cdot (\bar \psi \vdash \phi) =_x (\hat \psi \vdash \phi) $ unless $\phi =_x \hat \psi$, in which case $\psi \cdot (\bar \psi \vdash \phi) =_x (\hat \psi \vdash \phi),$ unless $\psi \in \sK$ or $\bar \psi \vdash \phi \in \sK$, 
	\item For all $\phi, \psi \in \sC[x]$, if $(\chi \cdot \psi) \cdot \phi$ is defined, then 
		 $$(\chi \cdot \psi) \cdot \phi =_x \chi \cdot (\psi \cdot \phi),$$
	\item Composition $(\cdot)$, combination $\langle, \rangle$, and the turnstile $(\vdash)$ in $\sC(x)$ is substitutive in both arguments:
		\enu
			\item if $\phi =_x \psi$ then $\phi \vdash \chi =_x \psi \vdash \chi$ and $\chi \vdash \phi =_x \chi \vdash \psi$,
			\item if $\phi =_x \psi$ and $\phi \cdot \chi$ $\downarrow$ then $\phi \cdot \chi =_x \psi \cdot \chi$ and if $\chi' \cdot \phi$ $\downarrow$ then $\chi' \cdot \phi =_x \chi' \cdot \psi$,
			\item if $\phi =_x \psi$ and $\langle \phi, \chi \rangle$ $\downarrow$ then $\langle \phi, \chi \rangle =_x \langle \psi, \chi \rangle$ and $\langle \chi, \phi \rangle =_x \langle \chi, \psi \rangle$,
		\Enu 
	\item For all $\phi: a \to \trut$, $f =_x a \vdash \trut$,
	\item For all pairs $(a,b) \in \sC$ (viewed as a deductive system), if the unique good pair for $(a,b)$ is $(\pi_{a,b}, \pi'_{a,b})$ and any good evaluation $\epsilon_{a,b}$ is taken, then these are required to satisfy their usual equational properties in expressions involving $x$: 
		\enu
			\item if $\pi_{a,b} \langle \phi, \psi \rangle$ $\downarrow$ then $\pi_{a,b} \langle \phi, \psi \rangle =_x \phi$,
			\item if $\pi'_{a,b} \langle \phi, \psi \rangle$ $\downarrow$ then $\pi'_{a,b} \langle \phi, \psi \rangle =_x \psi$,
			\item if $\langle \pi_{a,b} \cdot \phi, \pi'_{a,b} \cdot \phi \rangle$ $\downarrow$ then $\langle \pi_{a,b} \cdot \phi, \pi'_{a,b} \cdot \phi \rangle =_x \phi$,
			\item if $\langle \phi \cdot \pi_{a,b} , \psi \cdot \pi'_{a,b} \rangle$ $\downarrow$ then $\langle \phi \cdot \pi_{a,b}, \psi \cdot \pi'_{a,b} \rangle = \phi \wedge \psi$,
			\item if $\epsilon \langle \phi^* \cdot \pi_{c,a}, \pi'_{c,a} \rangle$ $\downarrow$ then $\epsilon \langle \phi^* \cdot \pi_{c,a}, \pi'_{c,a} \rangle = \phi$, 
			\item if $(\epsilon \langle \phi \cdot \pi_{c,a}, \pi'_{c,a} \rangle )^*$ $\downarrow$ then $(\epsilon \langle \phi \cdot \pi_{c,a}, \pi'_{c,a} \rangle )^* = \phi$. 
		\Enu
\Enu
\Dfn

The construction of $\sC(x)$ is thus carried out closely following Lambek. 
By iterating the construction of Definition \ref{d.polyicat} we may define general polynomial systems $\sA[\vec{x}]$ and general polynomial categories $\sC(\vec{x})$. 
A {\em polynomial over $\sC$} is an element of $\sC(\vec x)$ for any sequence of indeterminates $\vec x$. %

The following properties are established in \cite{LaSc1} for ordinary cartesian closed categories. The proof in our setting is similar when source and target do not depend on $x$, but in general requires a recursive step:

\lem\label{l.fCompT-iccc}
Let $\sC$ be an ideal cartesian closed category. Then $\sC(x)$ is an ideal closed category, and moreover:
\enu
	\item For every ideal cartesian closed category $\sD$, for every $F: \sC \to \sD$, and for every $a: F(\bar x) \to F(\hat x)$ in $\sD$, %
		there exists a unique functor $\theta: \sC(x) \to \sD$ satisfying %
			$$ \theta(x) = a, \quad \theta(f) = F(f) \,\, \text{ for all } f \in \sC.$$
	\item As a consequence of (1), for every $a \in \sC$, there is a unique functor $S_x^a: \sC(x) \to \sC$ (called {\em substitution of $a$ for $x$}) satisfying
			$$S_x^a (x) = a, \qquad S_x^a (f) = f \,\, \text{ for all } f \in \sC.$$
\Enu
\Lem



\thm\label{t.fCompT}
Let $\sC$ be an ideal cartesian closed category, %
let $\phi \in \sC(x),$ where $\phi: \trut \to \hat \phi$. 
Then %
there exists a unique element $g:\hat x \to \hat{\phi}$ in $\sC$, such that
$$\phi = g \cdot x$$
in $\sC(x)$.
\Thm
\prf
The proof we give, following Lambek, proceeds by passing through $\sC[x]$, the polynomial generalized positive deductive system over $\sC$, and then verifying that one is able to mod out by $=_x$. 
%
First we show that $\kappa_x \phi$ has a new behavior because of $=_x$:

\lem\label{l.fCompT}
$\kappa_x \phi$ is a well-defined element of $\sC(x)$, satisfies
$$\kappa_x \phi \cdot \langle x , \trut \rangle = \phi,$$
and is the unique element of $\sC(x)$ that does so. 
\Lem
\prf 
One must check that
$$\eif \phi =_x \psi, \ethen \kappa_x \phi =_x \kappa_x \psi.$$
This requires checking each of the relations
We need only check the new case created by $\wedge$; the other cases can be checked as in \cite{LaSc1}. This follows from the definition of $\kappa_x$: for any $\phi, \psi$ in $\sC(x)$ we have $\phi \wedge \psi =_x \langle \phi \cdot \pi , \psi \cdot \pi' \rangle$. We verify that
\begin{align*}
\kappa_x (\phi \wedge \psi) 		&= \langle \kappa_x (\phi) \pi , \kappa_x (\psi) \pi' \rangle \\
						&= \langle \kappa_x (\phi \cdot \pi) , \kappa_x (\psi \cdot \pi') \rangle \\ 
						&= \kappa_x (\langle \phi \cdot \pi, \psi \cdot \pi' \rangle).     
\end{align*}
from the definition of $\kappa_x$ for this case. 
The uniqueness of the choice of $\ksi(\phi)$ is the result of the following calculation in $\sC(x)$ \cite{LaSc1, LaK2}:
\begin{align*}
\hspace{170pt}
\kappa_x \phi 	&=_x \kappa_x (\tilde f \cdot \langle x , \trut \rangle) \\
			&=_x \tilde f \cdot \kappa_x ( \langle x , \trut ) \\
			&=_x \tilde f \cdot \langle \kappa_x x , \kappa_x \trut \rangle \\
			&=_x \tilde f \cdot \langle \pi_{\hat x, \trut}, \trut \cdot \pi'_{\hat x, \trut} \rangle \\
			&=_x \tilde f. \hspace{230pt}\qedhere
\end{align*}
\Prf

Now we finish the proof of Theorem \ref{t.fCompT}. %
We define the element $g$ in $\sC(x)$ to be
$$g := \kappa_x \phi \cdot \beta,$$
where $\beta$ is just the obvious munging term, in fact $\beta \threeline \langle 1_{\hat x}, \hat{x} \vdash \trut \rangle$. %
Indeed, we have
\begin{align*}
g \cdot x 	
		&= \kappa_x \phi \cdot \langle 1_{\hat x}, \hat{x} \vdash \trut \rangle \\
		&= \kappa_x \phi \cdot \langle x, \trut \vdash \trut \rangle \\
		&= \kappa_x \phi \cdot \langle x, \trut \rangle \\
		&= \phi
\end{align*}
by Lemma \ref{l.fCompT}. %
For uniqueness of $g$, suppose that $\tilde g \in \sC$ satisfies $\tilde g \cdot x = \phi$ in $\sC(x)$. We calculate
\begin{align*}
\kappa_x (\phi) \cdot \beta 	
		&= \kappa_x (\tilde g \cdot x) \cdot \beta \\
		&= \kappa_x (\tilde g \cdot x) \cdot \langle 1_{\hat x}, \hat x \vdash \trut \rangle \\
		&= \kappa_x (\tilde g) \cdot \langle \pi_{\hat x, \trut}, \kappa_x x \rangle \cdot \langle 1_{\hat x} , \hat x \vdash \trut \rangle \\
		&= \tilde g \cdot \pi'_{\hat x, \hat x} \langle \pi_{\hat x, \trut} , \pi_{\hat x, \trut} \rangle \cdot \langle 1_{\hat x} , \hat x \vdash \trut \rangle \\
		&= \tilde g \cdot \pi_{\hat x, \trut} \cdot \langle 1_{\hat x}, \hat x \vdash \trut \rangle \\
		&= \tilde g \cdot 1_{\hat x} \\
		&= \tilde g.
\end{align*}
But $\kappa_x (\phi) \cdot \beta = g$ by definition of $g$. So $g = \tilde g$, and $g$ is unique. 
\Prf


From Theorem \ref{t.fCompT} we define notation (to resemble a counit) $\varepsilon_x \phi: \hat x \to \hat \phi$ by
$$\varepsilon_x \phi := g = \kappa_x (\phi) \cdot \beta.$$
Theorem \ref{t.fCompT} has the following corollary:

\cor\label{c.fCompT}
Let $\sC$ be an ideal cartesian closed category, %
and let $\phi \in \sC(x)$ have source $\trut$. %
Then there exists a unique element $h: \trut \to (\hat x \vdash \hat \phi)$ such that 
$$\phi =_x \epsilon \cdot \langle h, x \rangle$$
in $\sC(x)$. %
\Cor
\prf
This is obtained by taking the name of the element $g$ of Theorem \ref{t.fCompT}: that is, take
\[
h = \name{g}. \qedhere
\]
\Prf

From Corollary \ref{c.fCompT} we define notation $\lambda_x \phi : \trut \to (\hat x \vdash \hat \phi)$ by
\[
\lambda_x \phi := h = \name{\kappa_x (\phi) \cdot \langle 1_{\hat x} , \hat x \vdash \trut \rangle} 
\]


As an aside, we observe from the proofs of Theorem \ref{t.dedt} and \ref{t.fCompT} that $\wedge$'s identity in categories suggests whether the symbol may be sugared out of generalized deduction systems entirely. This would mean $\langle, \rangle$ would be defined as a basic operation subject to an equational axiom: 
$$t(\langle a, b \rangle) = \langle \hat a \pi, \hat b \pi' \rangle.$$
In this case $\pi$ and $\pi'$ must satisfy a self-referential axiom:
$$s(\pi) = s(\pi') = \langle a \cdot \pi, b \cdot \pi' \rangle.$$


%

%

%

\section[Typed Lambda Calculus and the Main Correspondence]{\onehalfspacing Typed Lambda Calculus and the Main Correspondence}\label{s.lam}

In this section we will finally observe what happens on the syntactic side of the correspondence after generalizing semantics. As it turns out,  types acquire a richer structure and simultaneously assume the role of function constants. %
By a {\em generalized lambda calculus} (Definition \ref{d.lam}) we refer to the simplest such type system possible: 
we do not make mention of natural numbers objects (see \cite{LaSc1}), Boolean types, or other features that may appear in applications of lambda calculus. 

%
%
%
%
%
%
%
%
%
%
The next definition is not used in the sequel. It is included in order to establish a basis for defining variables before making Definition \ref{d.lam}.

\dfn\label{d.prelam}
A {\em pre-generalized typed lambda calculus} is a structure
$$(\Lambda, \sT_\Lambda, \sS_\Lambda, s,t, \cdot, \vdash, 
\trut, \wedge, \ty, \name{}, *, ()^\cdot, \pi, \pi', \wr, \langle,\rangle, \lambda, \sV_\Lambda)$$
where
\enu
	\item $\Lambda$ is a set,
	\item $\sT_\Lambda$ and $\sS_\Lambda$ are disjoint subsets of $\Lambda$ and $\sT_\Lambda \cup \sS_\Lambda = \Lambda$,
	\item $\sV_\Lambda$ is a subset of $\Lambda$,	
	\item the system
		$$(\sT_\Lambda, s,t,\cdot,\vdash, \sV')$$
		is an ideal category, where $\sV' = \sV_\Lambda \cap \sT_\Lambda$,
	\item $\trut$ is a designated element of $\sT_\Lambda$,
	\item $\wedge$ is a mapping $\sT_\Lambda \times \sT_\Lambda \to \sT_\Lambda$,
	\item $\name{}$ is a mapping $\sT_\Lambda \to \sS_\Lambda$,
	\item $\ty$ is a mapping $\sS_\Lambda \to \sT_\Lambda$,
	\item and in $\sS_\Lambda$:
	\enu
		\item $*$ is a designated element of $\sS_\Lambda$,
		\item $()^\cdot$ is a mapping $\Lambda \to \sS_\Lambda$,
		\item $\pi$, $\pi'$ are partially defined mappings $\sS_\Lambda \to \sS_\Lambda$,
		\item $\wr$ and $\langle,\rangle$ are partially defined mappings $\sS_\Lambda \times \sS_\Lambda \to \sS_\Lambda$,
		\item $\lambda$ is a mapping $\sX \times \sS \to \sS$, where $\sX$ is defined below,
	\Enu
\Enu
subject to the conditions
\enu
	\item 
		$\hat \trut = \bar \trut = \trut,$
 	\item 
		for all $s \in \sS_\Lambda$, $\pi(s)$ $\downarrow$ iff $\pi'(s)$ $\downarrow$ iff there exist $A,B \in \sT_\Lambda$ such that $\ty(s) = A \wedge B$, 
	\item 
		$s \wr t$ $\downarrow$ iff there exist $A,B \in \sT_\Lambda$ such that $\ty(s) = A \vdash B$ and $\ty(t) = A$,
	\item 
		$\langle s,t \rangle$ $\downarrow$ iff $\ty(s) = \ty(t)$,
\Enu
typing conditions
\enu
	\item 
		$\ty(\name{A}) = \bar A \vdash \hat A$,
 	\item 
		$\ty(*) = \trut$,
 	\item 
		for all $\alpha \in \Lambda$, $\ty(\alpha^\cdot) = \ty(\alpha),$
 	\item 
		if $s \in \sS_\Lambda$ and $\ty(s) = A \wedge B$, then $\ty(\pi(s)) = A$ and $\ty(\pi'(s)) = B$,
	\item 
		if $s \wr t$ $\downarrow$, then $\ty(s \wr t) = \widehat{\ty(s)}$,
	\item 
		if $\langle s, t \rangle$ $\downarrow$, then $\ty(\langle s, t \rangle) = \ty(s) \wedge \ty(t)$,
	\item 
		if $\lambda(x, s)$ $\downarrow$, then $\ty(\lambda(x, s)) = \ty(x) \vdash \ty(s)$,
\Enu
and the validities
\enu
	\item 
		$* \in \sV_\Lambda$,
	\item 
		if $A,B \in \sV_\Lambda$, then $A \wedge B \in \sV_\Lambda$,
	\item 
		(witnesses, propositions-as-types) if $s \in \sV_\Lambda$, then $\ty(s) \in \sV_\Lambda$,
	\item 
		If $A \in \sV_\Lambda$, then $\name{A} \in \sV_\Lambda$.
 	\item 
		$* \in \sV_\Lambda$,
 	\item 
		if $c \in \sV_\Lambda$ and $\pi(c), \pi'(c)$ $\downarrow$, then $\pi(c), \pi'(c) \in \sV_\Lambda$,
	\item 
		if $a,f \in \sV$, then $f \wr a \in \sV$,
	\item 
		$a,b \in \sV_\Lambda$ implies $\langle a, b \rangle \in \sV_\Lambda$,
	\item 
		if $s \in \sV_\Lambda$, then $\lambda(x,s) \in \sV_\Lambda$. 
\Enu
\Dfn

Note that many type theories, e.g. \cite{MoI1}, include {\em function constants} $f: A \to B$ as well as terms and types; in this formalism (guided by the new semantics) function constants are indistinguishable from types, and together with objects they form a category. Types behave as function constants via the derived operation $A \star s := \name{A} \wr s.$ The operation $\name{}$ is used not only here but also in the construction of $\boC \Lambda$ in Definition \ref{d.catcon}. 

Elements of $\sT_\Lambda$ are called {\em types}, and elements of $\sS_\Lambda$ are called {\em terms}. For a term $s$, the element $\ty(s)$ of $\sT_\Lambda$ is called the {\em type of $s$.} %
We may write $s:T$ to denote the relation $\ty(s) = T$.
A term of the form $\alpha^\cdot$ for some $\alpha$ (which may be a type or a term) is called a {\em variable}. %
We may iterate the operation $()^\cdot$, %
and we do not allow $()^\cdot$ to be substitutive in its argument. %
Therefore we may assume that the symbol $x_i$ 
unpacks to %
$((\dots ((A)^\cdot )^\cdot \dots )^\cdot )^\cdot$. %
In this way, we have a countable stock $x_1, x_2, \dots$ of distinct ``standard'' variables for each type $A$. %
For technical reasons (see below, before Definition \ref{d.lam}), we take these standard variables to be the only variables of $\Lambda$, 
and we place the obvious (total) ordering on variables of each type. %
A variable $x_i$ is {\em free} in a term if it appears in the term, unless it appears but only within a well-formed expression of the form $\lambda(x_i, s)$. In this case we say it appears {\em captured} or {\em bound}. 
We define the mapping on terms
$$\FV(s) = \set{x \in \sX \mid x \text{ appears free in $s$, and } x \nin \sV_\Lambda},$$
where the phrase ``appears free'' has its usual meaning, except that we assume that no variable {\em appears free} in any type. So, for example, for all types $A$, $\FV(\name{A})$ is empty. %
If $s$ is a term, $x$ is a variable, and $t$ is a term whose type is the same as the type of $x$ 
we define notation
$$s[x/t]$$
to be the term $s$ with the variable $x$ replaced by $t'$ in each instance where it does not appear bound in $s$, %
where $t'$ is $t$ with any variable $y \in \FV(t)$ that appears captured in $s$, that is, 
$$y \in \FV(t) \cap \CAP(s),$$
where $\CAP(s)$ is the set of variables appearing captured in $s$, %
replaced by a variable of the same type that is not in the set $\VAR(s) \cup \VAR(t)$ of variables appearing in either $s$ or $t$. %
These choices are made in {\em the simplest order-preserving way}, by which is meant that once the set of variables to be changed is found, the entire set is incremented by the smallest positive integer such that the set of variables so generated is not in $\VAR(s) \cup \VAR(t)$. %
These incrementing operations are associative, as is the substitution operation itself. Hence we have
$$s[x/t][y/r] = s[x/t[y/r]]$$
for all terms $s,r,t$ and variables $x,y$. %
We may often ignore the extra step involving $t'$, for it is only necessary because we have not set terms $s$ and $s'$ equal in $\Lambda$ which are the same up to one or more free variables (a form of $\alpha$-conversion) in $\prelambdaCalc$ or in the category $\lambdaCalc$ defined next. %
Note that a morphism in $\prelambdaCalc$ sends closed terms to closed terms. %

\dfn\label{d.lam}
A {\em generalized typed lambda calculus} is %
a pre-generalized typed lambda calculus %
on which there is an equality relation on the terms $\sS_\Lambda$ of $\Lambda$ defined as follows: %
Let $\sP$ be the finite power set $\sP_{fi} (\sX)$ of $\sX$. %
For each finite set $\bar x = \set{x_1, \dots, x_n}$ in $\sP$, %
let 
$$\sR(\Lambda, \bar x) := \set{s \in \sS_\Lambda \mid \FV(s) \lies \bar x}.$$
We define the relation $=_{\bar x}$ on $\sR(\Lambda, \bar x)$ 
to be the smallest equivalence relation that satisfies
\enu
	\item $=_{\bar x}$ is reflexive, symmetric, and transitive,
	\item Substitutivity conditions:
	\enu
		\item 
			if $s =_{\bar x} t$, and $\pi(s)$ $\downarrow$, then $\pi(s) =_{\bar x} \pi(t)$, and $\pi'(s) =_{\bar x} \pi'(t)$,
		\item 
			if $s =_{\bar x} t$, then $s \wr r =_{\bar x} t \wr r$ and $u \wr s =_{\bar x} u \wr t$ whenever these expressions are well-defined, 
		\item 
			if $s =_{\bar x} t$, and $\langle s,r \rangle$ $\downarrow$, then $\langle s,r \rangle =_{\bar x} \langle t,r \rangle$, and similarly in the second argument,
		\item 
			if $s =_{\bar x} t$, then $s \wr r =_{\bar x} t \wr r$ and $u \wr s =_{\bar x} u \wr t$ whenever these expressions are well-defined, 
	\Enu
	\item for all $s:\trut$, $s =_{\bar x} *$,
	\item for all $a : A, b: B,$ 
		$$\pi(\langle a, b \rangle) =_{\bar x} a,$$
		$$\pi'(\langle a, b \rangle) =_{\bar x} b,$$
	\item for all $c: A \wedge B$,
		$$\langle \pi(c), \pi'(c) \rangle =_{\bar x} c,$$
	\item For all terms $s \in \sS_\Lambda$, terms $a \in \sS_\Lambda$, and variable $x$ that may appear in $\bar x$,
	\enu
		\item $(\lambda(x,s)) \wr a =_{\bar x} s[x/a],$
		\item $\lambda(x, s \wr x) =_{\bar x} s$,
		\item if $\FV(s) = \set{x}$, there exists a unique $A \in \sT_\Lambda$ such that $s =_{\set{x}} \name{A} \wr x$. \label{ax.unname}
		\item ($\alpha$-conversion for lambda terms) 
		$$\lambda(y, s) =_{\bar x} \lambda(y', s[y/y'])$$
		if $\ty(y) = \ty(y')$ and $y' \nin \FV(s)$.
	\Enu
\Enu
We observe that $\FV()$ is still well-defined. %
We denote by 
$$s^\natural$$ 
the type $A$ given by Axiom \ref{ax.unname}. 
We impose the condition on the $=_{\bar x}$'s that: 
\enu
	\item if $\bar x \lies \bar y$ then for all $s,t \in \sR(\Lambda, \bar y)$, $s =_{\bar x} t$ implies $s =_{\bar y} t$.
\Enu
Because %
(1) $s =_{\FV s} s$, and %
(2) if $s =_{\bar x} t$ and $s =_{\bar y} t$, then there exists a finite set $\bar z$ such that $s =_{\bar z} t$ and $\bar x, \bar y \lies \bar z$, %
we may define an equivalence relation {\em equality in $\sS_\Lambda$} on the set $\sS_\Lambda$ of terms of $\Lambda$ by
$$s = t \,\,\,\text{ if }\,\, s =_{\bar x} t \text{ for some $\bar x$ in $\sP$}.$$
A {\em morphism} $\Phi: \Lambda \to \Mu$ of generalized typed lambda calculi, also called a {\em translation}, is a mapping 
$$\Phi: \Lambda \to \Mu$$
that satisfies the following, where equalities between terms are interpreted as equality in $\sS_\Lambda$:
\enu
	\item for all $A \in \sT_\Lambda, s \in \sS_\Lambda$, $\Phi(A) \in \sT_\Mu$ and $\Phi(s) \in \sS_\Mu,$
	\item the restriction of $\Phi$ to $\sT_\Lambda$ is a morphism of ideal categories that satisfies
		$$\Phi(\trut_\Lambda) = \trut_\Mu,$$
		$$\Phi(A \wedge B) = \Phi(A) \wedge \Phi(B),$$
	\item if $s =_{\bar x} t$, then $\Phi(s) =_{\Phi(\bar x)} \Phi(t)$. \label{ax.preservesequalities}
	\item 
		$\Phi(\ty(s)) = \ty(\Phi(s)),$
	\item 
		$\Phi(\name{A}) = \name{\Phi(A)},$
	\item 
		$\Phi(*) = *$, 
	\item 
		for all $\alpha \in \Lambda$, $\Phi(\alpha^\cdot) = \Phi(\alpha)^\cdot,$
	\item 
		$\Phi(\pi(c)) = \pi(\Phi(c),$ and $\Phi(\pi'(c)) = \pi'(\Phi(c)),$
	\item 
		$\Phi(s \wr t) = \Phi(s) \wr \Phi(t),$
	\item 
		$\Phi(\langle s, t \rangle) = \langle \Phi(s), \Phi(t) \rangle,$
	\item 
		$\Phi(\lambda(x,s)) = \lambda(\Phi(x), \Phi(s)).$
\Enu
As a consequence of (\ref{ax.preservesequalities}), $\Phi$ preserves equalities in $\Lambda$:
$$s = t \eimplies \Phi(s) = \Phi(t).$$
This gives a category $\lambdaCalc$ of generalized typed lambda calculi.
\Dfn

Given a generalized typed lambda calculus $\Lambda$, we can construct an ideal cartesian closed category using Theorem \ref{t.fCompT}:

\dfn\label{d.catcon}
Let $\Lambda$ be a typed lambda calculus. %
Let $\sB_\Lambda$ be the set of {\em bulletins} in $\Lambda$, that is, the set of terms in $\Lambda$ that have only one free variable. %
Also for $A \in \sT_\Lambda$, let 
$$\settl{\bullet}{A} := \settl{x}{\name{A}\wr x}, \quad x : \bar A,$$
that is, a symbol $\settl{x}{s}$ where $x$ is a variable of type $\bar A$, and $s$ is the term $\name{A} \wr x$. %
By Axiom (\ref{ax.unname}) of Definition \ref{d.lam} we may identify these symbols with types in $\Lambda$. 
We define $\boC\Lambda$ to be the set
$$\boC\Lambda := \set{\settl{x}{s} \mid s \in \sB_\Lambda, \text{ $x$ a variable}},$$ 
of symbols $\settl{x}{s}$ for variable $x$ and bulletin $s$, equipped with the structure
\begin{align*}
		\overline{\settl{x}{s}} &:= \settl{\bullet}{\ty(x)}, \\ 
		\widehat{\settl{x}{s}} &:= \settl{\bullet}{\ty(s)}, \\ 
		\settl{x}{s} \cdot \settl{y}{t} &:= \settl{y'}{s[x/t]}, \\
		\settl{x}{s} \vdash \settl{y}{t} &:= \settl{u}{v}, \quad u:\ty(s), v:\ty(t), 
\end{align*}
where
$$y' = 
	\begin{cases} 	
				y 
						& \text{ if $\FV(t)$ is empty}, \\ 
				inc_n(y)	
						& \text{ if $\FV(t) = \set{u}$ and $\FV(s[x/t]) = \set{inc_n (u)}$,}
	\end{cases}
$$
where $inc_n$ is the modification of the variable described after Definition \ref{d.prelam}. %
Let $\sK_{\boC\Lambda}$ be the set of symbols $\settl{x}{k}$ where $k$ is a constant in $\Lambda$, that is, $\FV(k) = \nll$. %
Let equality of symbols in $\boC\Lambda$ be defined by
\enu
	\item $\settl{x}{s} = \settl{y}{t}$ if $\ty(x) = \ty(y), \ty(s) = \ty(t),$ and there is $z:\ty(x)$ such that $s[x/z] = t$, 
	\item $\settl{x}{s} = \settl{x}{u}, \quad u:\ty{s}, \,\,$ \text{ if $\FV(s) = \set{y}$ and $y \neq x$}, 
	\item if $\ty(s) = \trut$, then $\settl{x}{s} = \settl{x}{*}$, 
	\item for all bulletins $s$ and all variables $x,y$ of the same type, $\settl{x}{s} = \settl{y}{s[x/y]}$. 
\Enu
This gives an ideal category $\boC\Lambda$, where 
the identity of $\settl{x}{s}$ is 
$$1_{\settl{x}{s}} = \settl{y}{y}, \quad y:s^{\natural},$$
terminal arrows are of the form 
$$\settl{y}{*},$$
and types (in the sense of section \ref{s.th}) are of the form
$$\settl{x}{y}, \quad x \neq y.$$
Validities defining $\boC\Lambda$ are the evident ones based on Definition \ref{d.icat}. 
\Dfn

We have an ideal category $\boC\Lambda$, but we have not directly made any assumptions about the category $\sT_\Lambda$. Nevertheless, we have:


\prop\label{p.catcon}
$\boC\Lambda$ is an ideal cartesian closed category. 
\Prop
\prf
Set
\begin{align*}
		\trut &:= \settl{u}{*}, \quad u:\trut_\Lambda,\\
		\settl{x}{s} \wedge \settl{y}{t} &:= \settl{z}{\langle s \wr \pi(z), t \wr \pi'(z) \rangle}, \\
		\langle \settl{x}{s}, \settl{y}{t} \rangle &:= \settl{z}{\langle s[x/z], t[y/z] \rangle}, \\
		\settl{z}{s}^* &:= \settl{x}{\lambda(y,s \wr \langle x,y \rangle) }, \quad \text{ where $z:A \times B, x:A$}, \\
		\pi &:= \settl{z}{\pi(z)}, \\
		\pi' &:= \settl{z}{\pi'(z)},\\
		\epsilon &:= \settl{z}{\pi(z) \wr \pi'(z)}, 
\end{align*}
with validities as needed (Definition \ref{d.iccc}). 
\Prf

%
%
%
%
%
%
%
%
%
%

We can also construct a typed lambda calculus from the data of a cartesian closed category:

\dfn\label{d.internallanguage}
Let $\sC$ be an ideal cartesian closed category. We define the symbol $\boL\sC$ as follows:
\enu
	\item The set of types of $\boL\sC$ is the set of symbols $A_f$ indexed by elements $f \in \sC$:
		$$\sT_{\boL\sC} := \set{A_f \mid f \in \sC},$$
		in fact we set $A_f = f$ and take $\sC$ itself as the set of types (this is needed for the proof of Theorem \ref{t.chl}), however, we use the notation $A_f$ at times when it seems to lessen the potential for confusion.
	\item The set of terms of $\boL\sC$ is the set of polynomials $\phi$ over $\sC$ sourced at $\trut$, that is, 
		$$\sS_{\boL\sC} := \set{\phi \mid \phi \in \sC[\vec x] \text{ for some $\vec x$, } \hat{\phi} \text{ is in $\sC$, and } \bar{\phi} = \trut}, $$
		where we assume 
		that indeterminates have internal structure given by the syntax $()^\cdot$. 
	\item Define
		\begin{align*}
			\ty(\phi) &:= A_{\hat \phi}, \\
		 	s(A_f) &:= A_{sf}, \\
			t(A_f) &:= A_{tf}, \\
			A_f \cdot A_g &:= A_{f \cdot g}, \\
			A_f \wedge A_g &:= A_{f \wedge g}, \\
			A_f \vdash A_g &:= A_{f \vdash g}, \\
			\name{A_f} &:= \name{f}, \quad \text{ the name of $f$,} \\
			\trut_{\boL\sC} &:= \trut_\sC \vdash \trut_\sC, \\
			* &:= \trut_\sC, \\
			\sK_{\boL\sC} &\text{ is the set of constant polynomials.} %
		\end{align*}
	\item if $\phi$ is a bulletin in $x$ over $\sC$, then let 
		$$\phi^\natural := A_{\epsilon_x \phi}.$$
	\item $\sV_{\boL\sC}$ is the set $\set{A_f  \mid f \in \sV_\sC}$ %
	joined with the set of valid constant terms, %
	joined with the set of polynomials valid according to Definition \ref{d.poly}. %
	\item Define
		\begin{align*}	
			\trut_{\boL\sC} 	&:= \trut \vdash \trut, \\
			*_{\boL\sC} 	&:= \trut,
		\end{align*}
\Enu
Then we have a pre-generalized typed lambda calculus. %
We make from this a generalized typed lambda calculus by imposing the equality relation on terms inherited from equality in $\sC(\vec x)$:
the equality relation $=_{\vec x}$ is defined to be equality in $\sC(\vec x)$, along with the usual inclusions of polynomial systems in one another. 
\Dfn


$\boL\sC$ is called the {\em internal language} of the ideal cartesian closed category $\sC$. %
%
%
%
%
%
%
%
%
%
%
Next, we verify that these constructions are functorial: 

\prop\label{p.functor}
We have the following:
\enu
	\item $\boC$ is a functor from $\lambdaCalc$ to ${\bf ICCC}$. %
	\item $\boL$ is a functor from ${\bf ICCC}$ to $\lambdaCalc$. %
\Enu
\Prop
\prf
Given $\Phi: \Lambda \to \Lambda'$, we define $\boC\Phi : \boC\Lambda \to \boC\Lambda'$ by
$$\boC\Phi \settl{x}{s} := \settl{\Phi(x)}{\Phi(s)}$$
for 
$\settl{x}{s} \in \boC\Lambda.$ 
Now we check that $\boC \Phi$ is a morphism in ${\bf ICCC}$, %
and that %
$\boC$ is a functor (Definition \ref{d.functor}). %

Let $F: \sC \to \sD$ in ${\bf ICCC}$. Define a mapping $\boL F: \boL(\sC) \to \boL(\sD)$ by
\begin{align*}
\boL F(A_f) &:= A_{F(f)}, \\
\boL F ( \alpha^\cdot) &:= (\boL F(\alpha))^\cdot, 
\end{align*}
and extend $F$ from $\sC$ to polynomials over $\sC$ in the most straightforward way. %
Now we check that $\boL F$ is a morphism in $\lambdaCalc$, and that $\boL$ is indeed a functor. %
\Prf

%
%
%
%
%
%
%
%
%
%


\dfn\label{d.eta}
Let $\Lambda$ be a generalized typed lambda calculus. Define a mapping $\Lambda$ to $\boL\boC\Lambda$ 
by defining, in the pre-generalized typed lambda calculus $\Lambda_0$ obtained by ignoring equalities in $\sS_\Lambda$, 
\begin{align*}
\eta_\Lambda(A) &:= A_{\settl{\bullet}{A}}	& A \in \sT_\Lambda, 
		 			\\
\eta_\Lambda(k) &:= \settl{x}{k}, 			& k \in \sS_\Lambda, \FV(k) = \nll, \ty(x) = \trut_{\boL\boC\Lambda},	\\
\eta_\Lambda(\alpha^\cdot) &:= (\eta_\Lambda(\alpha))^\cdot,					& \alpha \in \Lambda,		\\
\eta_\Lambda(\pi(\phi)) &:= \pi(\eta_\Lambda \phi ), 												& \\
\eta_\Lambda(\pi'(\phi)) &:= \pi'(\eta_\Lambda \phi ),												& \\
\eta_\Lambda(\langle \phi, \psi \rangle) &:= \langle \eta_\Lambda (\phi) , \eta_\Lambda( \psi) \rangle 			& \\
\eta_\Lambda(\phi \wr \psi) &:= \eta_\Lambda (\phi) \wr \eta_\Lambda (\psi)								& \\
\eta_\Lambda(\lambda(x, \phi)) &:= \lambda(\eta_\Lambda(x), \eta_\Lambda (\phi))			& x \in \sX_\Lambda 
\end{align*}
The map $\eta_\Lambda$ is well-defined upon passage to $\Lambda$, since analogous equalities between polynomials hold in both $\Lambda$ and $\boL\boC\Lambda$. 
\Dfn

An alternative approach (really the same) to Definition \ref{d.eta} is via an isomorphism with a lambda calculus with parameter \cite{LaSc1}: 


\dfn\label{d.indeterminatevariable}
Let $\Lambda$ be a generalized typed lambda calculus, and let $x \in \sX_\Lambda$ be a variable. We define the symbol 
$$\Lambda_x$$
to be the generalized typed lambda calculus is defined exactly as $\Lambda$, except that
$$\sV_{\Lambda_x} := \set{x} \cup \sV,$$
that is, $x$ is taken to be a validating term in $\Lambda_x$. 
\Dfn

Intuitively, $\Lambda_x$ is $\Lambda$ with $x$ treated as a constant instead of as a variable. %

\lem\label{l.indeterminatevariable}
Let $\sC$ be an ideal cartesian closed category, and let $x$ be an indeterminate (with the variable syntax). Then the polynomial category 
$\boC\Lambda(x)$ over $\boC\Lambda$ is isomorphic to $\boC\Lambda_x$ in ${\bf ICCC}$. 
\Lem
\prf
By Proposition \ref{l.fCompT}, we need only check that $\boC\Lambda_{x}$ has the desired universal property of $\boC\Lambda(x)$. %
See \cite{LaSc1}. %
\Prf

Using Lemma \ref{l.indeterminatevariable}, we can identify polynomials $\tilde \phi$ over $\boC\Lambda$ with the corresponding symbol $\settl{u:\trut}{\phi(\vec x)}$ in $\boC\Lambda_{\vec x}$, where $\vec x = \FV(\phi)$ corresponds to the free variables $\ksi_1, \dots, \ksi_n$ of $\tilde \phi$ over $\boC\Lambda$ via the isomorphism. 


%
Finally, we have an extension of Lambek's equivalence between simply typed lambda calculi and cartesian closed categories:

\thm\label{t.chl}
The functors $\boC$ and $\boL$ form an equivalence
$$
\begin{tikzcd}
	\lambdaCalc \arrow[r, "\boC", shift left] \arrow[r, leftarrow, shift right, "\boL" below] & {\bf ICCC} \\
\end{tikzcd}\vspace{-3.5ex} 
$$
between $\lambdaCalc$ and ${\bf ICCC}$. 
\Thm
\prf
For $\sD$ in ${\bf ICCC}$, define $\varepsilon_\sD: \boC\boL \sD \to \sD$ to be the map 
$$\varepsilon_\sD : \settl{x}{\phi} \mapsto 
		\begin{cases} 
			\varepsilon_x \phi, & \text{if $\FV(\phi) = \set{x},$ or $\FV(\phi)$ is empty, or $\ty(\phi) = \trut_\Lambda$,} \\
			\hat x \vdash \hat \phi \,\,\text{ in $\sC$, } & \text{otherwise.} 
		\end{cases}
$$
This map is well-defined since if $\settl{x}{\phi} = \settl{y}{\psi}$, then $\phi[x/z] = \psi[y/z]$, where $z$ does not appear in $\phi$ or $\psi$. Let these be $\phi(z), \psi(z)$. Then $\varepsilon_z \phi(z) = \varepsilon_z \psi(z)$. But $z$ is eliminated by evaluation, so $\varepsilon_x \phi = \varepsilon_z \phi(z) = \varepsilon_z \psi(z) = \varepsilon_x \psi.$
Let $F: \sC \to \sD$ in ${\bf ICCC}$. Then to check that $\varepsilon: \sD \mapsto \varepsilon_\sD$ is a natural transformation, that is,
$$\varepsilon(\sD) \of \boC\boL(F) = F \of \varepsilon(\sC),$$
we check that for every $\settl{x}{\phi}$ in $\boC\boL \sC$, where $\phi$ is a bulletin in $x$ over $\sC$,
$$\varepsilon_\sD (\boC\boL F(\settl{x}{\phi})) = F(\varepsilon_\sC (\settl{x}{\phi})).$$
If $\phi$ is a non-constant bulletin in a variable different than the variable appearing in the symbol, then
\begin{align*}
\varepsilon_\sD (\boC\boL F \settl{x}{\phi}) 
	&= \varepsilon_\sD (\boC\boL F \settl{x}{y} ) \\
	&= \varepsilon_\sD \settl{x' : F(\ty(x))}{y' : F(\ty(y))} \\
	&= F(\ty(x)) \vdash F(\ty(y)) \\
	&= F( \ty(x) \vdash \ty(y) ) \\
	&= F(\varepsilon_\sC \settl{x}{\phi}).
\end{align*}
In the other cases, %
this reduces to checking that %
$$F(\varepsilon_x \phi) = \varepsilon_z \boL F \phi,$$
where $\boL F(x) \threeline z.$ We proceed by cases as in the proof of Theorem \ref{t.dedt}: 
if $\phi$ is a constant $k: \trut \to \hat k$, then 
\begin{align*}
F \varepsilon_x \phi 
	&= F(k) \cdot F(\pi'_{\trut, \trut} \cdot \langle 1_\trut , \trut \vdash \trut \rangle ) \\
	&= F(k) \cdot 1_\trut \\
	&= F(k) \\
	&= \boL F(k) \\
	&= \epsilon_z \boL F(k). 
\end{align*}
If $\phi$ is a variable $x: \trut \to \hat x$ equal to the variable captured by the symbol, then
\begin{align*}
F \varepsilon_x \phi
	&= F \varepsilon_x x \\
	&= F( \pi_{\trut, \hat x} \cdot \langle 1_\trut, \hat x \vdash \trut \rangle ) \\
	&= \pi_{\trut, \widehat{F(x)}} \cdot \langle 1_\trut, \widehat{F(x)} \vdash \trut \rangle ) \\
	&= \varepsilon_{\boL F (x)} \,\boL F(x). 
\end{align*}
The other cases are similar. 

For a generalized typed lambda calculus $\Lambda$ in $\lambdaCalc$, define $\eta(\Lambda) := \eta_\Lambda$ of Definition \ref{d.eta}.
To show that $\eta$ is a natural transformation, let $\Phi: \Lambda \to \Mu$ in $\lambdaCalc$. Then
$$\eta(\Mu) \of \Phi = \boL \boC(\Phi) \of \eta(\Lambda)$$
becomes, for types, 
$$\eta_\Mu (\Phi(A)) = \boL\boC\Phi (\eta_\Lambda (A)),$$
which is easily verified. Indeed,
\begin{align*}
\boL\boC\Phi(\eta_\Lambda (A))
	&= \boL\boC\Phi(A_{\settl{\bullet}{A}}) \\
	&= A_{\boC\Phi{\settl{\bullet}{A}}} \\
	&= A_{\settl{z}{\Phi(\name{A}) \wr z}}, \quad \ty z = \Phi(\bar A) = \overline{\Phi(A)}, \\
	&= A_{\settl{z}{\name{\Phi(A)} \wr z}} \\
	&= \eta_\Mu (\Phi(A)).
\end{align*}
For terms, we proceed by induction on the length of a term $s$ of $\Lambda$. 
If $s = k$ is a constant term (of length zero), 
\begin{align*}
\eta_\Mu (\Phi(k)) 
	&= \settl{u}{ \Phi k } \quad u:\trut \\
	&= \settl{u }{ \Phi k } \quad u: \Phi(\trut) \text{ since $\Phi(\trut) = \trut$} \\
	&= \boC \Phi \settl{x}{k} \\
	&= \boL \boC \Phi \settl{x}{k} \\
	&= \boL\boC\Phi (\eta_\Lambda(k)).
\end{align*}
If $s = x$, a variable of type $A$, then
\begin{align*}
\boL\boC\Phi(\eta_\Lambda(x)). 
	&= \boL \boC \Phi(\ksi), \quad \ksi: \settl{\bullet}{A} \\
	&= \Phi\ksi, \quad \Phi\ksi : \settl{\bullet}{\Phi A} \\
	&= \eta_\Mu \Phi (x). 
\end{align*}
We can similarly check the other cases $\pi(t), \pi'(t), t \wr r, \langle t, r \rangle, \lambda(y, t)$. 

Both $\eta_\Lambda$ and $\varepsilon_\Lambda$ are invertible as maps. %
Indeed, by Theorem \ref{t.fCompT}, $\varepsilon$ is injective, and also surjective (since $g \cdot y$ is itself a polynomial). %
To show that $\eta_\Lambda$ is invertible, we use Lemma \ref{l.indeterminatevariable}: if $\phi$ is a polynomial over $\boC\Lambda$ in variables $x_1, \dots, x_n$, we pass via the isomorphism of Lemma \ref{l.indeterminatevariable} from $\phi$ to an element $\phi'$ in $\boC\Lambda_{x_1, \dots, x_n}$ of the form $\settl{y}{t}$. Now note that $\eta_\Lambda(t) = \phi$, so $\eta_\Lambda$ is surjective. On the other hand if $\eta_\Lambda s = \eta_\Lambda t$, for two terms $s,t \in \sS_\Lambda$, then 
$\settl{u:\trut}{s} = \settl{u:\trut}{t} \quad \text{ in $\boC\Lambda_{x_1, \dots, x_n}$.} $ %
so $s = t$ as terms over $\Lambda$, by definition of equality in $\sS_{\boC\Lambda_{x_1, \dots, x_n}}$. %

Next we check (cf. Definition \ref{d.adjunction}) that the triangle laws hold. 
Let $\sC$ be in ${\bf ICCC}$. For a type $A_f$ in $\boL\sC$, 
\begin{align*}
\boL\varepsilon_\sC (\eta_{\boL\sC} (A_f)) 
	&= \boL\varepsilon_\sC (A_{\settl{\bullet}{A_f}}) \\
	&= A_{\varepsilon_\sC (\settl{\bullet}{A_f}) } \\
	&= A_{\varepsilon_z \name{f} \wr z}, \quad z:A_{\bar f}, \\
	&= A_{\varepsilon_z f\cdot z} \\
	&= A_f.
\end{align*}
Next, let $\phi$ be a term of $\boL\sC$, that is, a polynomial over $\sC$ in variables $x_1, \dots, x_n$, say. Then 
\begin{align*}
\boL\varepsilon_\sC (\eta_{\boL\sC} (\phi) ) 
	&= \varepsilon_{\sC(x_1, \dots, x_n)} (\settl{u:\trut}{\phi})  & \text{by Lemma \ref{l.indeterminatevariable}} \\
	&= \varepsilon_{u} \phi 		\\ 
	&= \phi.
\end{align*}
Next, if $\settl{x}{s}$ is an element of $\boC\Lambda$, then we must verify: %
$$\varepsilon_{\boC\Lambda} ( \boC\eta_{\Lambda} \settl{x}{s} ) = \settl{x}{s}.$$
The first case we check is that where $s$ is a bulletin in a variable not equal to that appearing in the symbol. %
Then 
$\settl{x}{s} = \settl{x}{y} = \ty(x) \vdash \ty(y)$ 
for some variable $y$, with $\ty(y) = \ty(x)$. We have
\begin{align*}
\varepsilon_{\boC\Lambda} (\boC \eta_\Lambda \settl{x}{s} ) 
	&= \varepsilon_{\boC\Lambda} ( \settl{\ksi}{\ksi'} ) , \quad \text{where $\ksi:\settl{\bullet}{\ty(s)}$, $\ksi':\settl{\bullet}{\ty(x)}$} \\
	&= \hat{\ksi} \vdash \hat{\ksi'} \\
	&= \settl{x}{y} \\
	&= \settl{x}{s}.
\end{align*}
Next, we check when %
$s = k$ is a constant term of type $B$ in $\Lambda$, and $x$ is a variable of type $A$ in $\Lambda$. Then 
\begin{align*}
\varepsilon_{\boC\Lambda} (\boC\eta_\Lambda \settl{x}{k} )
	&= \varepsilon_{\boC\Lambda} \settl{\eta_\Lambda x}{\eta_\Lambda k} \\
	&= \varepsilon_{\boC\Lambda} \settl{\ksi}{\settl{u}{k}}, 
		\text{ where $\ksi$ has type $\eta_\Lambda \ty(x)  = \settl{\bullet}{A} = \settl{v:\trut}{A \wr v}$, and $u:\trut_\Lambda$} \\
	&= \varepsilon_\ksi \settl{u}{k} \\
	&= \settl{u}{k} \cdot \arter_{\hat \ksi} \\
	&= \settl{u}{k} \cdot \settl{w:A}{*} \\
	&= \settl{w:A}{k[x/*]} \\
	&= \settl{w}{k} \\
	&= \settl{x}{k}.
\end{align*}
Next, if $s = x$ is a variable of type $A$ and is the same variable as that appearing in the symbol, then
\begin{align*}
\varepsilon_{\boC\Lambda} (\boC\eta_\Lambda \settl{x}{x} )
	&= \varepsilon_{\boC\Lambda} \settl{\eta_\Lambda x}{\eta_\Lambda x} \\
	&= \varepsilon_{\boC\Lambda} \settl{\ksi}{\ksi}, \quad \text{ where $\ksi$ has type $\settl{\bullet}{A} = A$},\\
	&= \varepsilon_\ksi \ksi \\
	&= \settl{\ksi}{\ksi} \\
	&= \settl{x}{x}.
\end{align*}
The other cases are proved similarly. Hence the triangle laws hold, and the theorem is proved. 
\Prf




\section{Conclusion and Future Work}\label{s.th-conc}


In this chapter we have shown that cartesian closed structure can be modified to include mappings on the set of objects that recover the base and the power of an exponential. 
We have indicated that the mathematics of cartesian closed categories is not affected by this addition, and moreover, by making this modification, we widen further the class of admissible functors (for some purposes relevant to categorical logic and type theory) 
to include arbitrary cartesian functors. We have also shown that this calculus extends beyond categories, to generalized categories. 
We have also presented a lambda calculus which permits the extension of the Curry-Howard-Lambek correspondence to the general case. 
Our work suggests that polynomials over categories and terms over types are in fact essentially the same thing. This can also be seen also in the ordinary categorical case, but in the generalized setting, the observation is made unavoidable. 
The fundamental insight of the Curry-Howard correspondence is thus that the cartesian closed structure on a cartesian closed category can be expressed almost entirely in terms of properties of objects in the space of polynomials. This seems (to the author) to be the mathematical content of the theorem. 

Because of the rich variety of subject matter in categorical logic and related subjects, there are a number of directions in which this work can be continued. 
%
For example, the work of Moggi on computational effects \cite{MoI1} has had an influence on much subsequent work (see for example Wadler\cite{WaR2}, Mulry, \cite{MuY1,MuY2}, Kobayashi \cite{Kobayashi1}, Semmelroth and Sabry \cite{SemmelrothSabry1}). In \cite{MoI1}, an extension of the lambda calculus is introduced and it is shown that it is possible to provide categorical semantics for computational effects by making use of monads. 
In fact, two constructions are presented. The first relates a cartesian closed category equipped with a monad to a monadic equational theory (one in which contexts consist of a unique typed variable) extended by a computational effect (he calls this the {\em simple metalanguage}), and the second relates a {\em strong} monad to a general equational theory (what Moggi calls an {\em algebraic} equational theory, this one called the {\em metalanguage}), i.e., one in which contexts may be arbitrary finite lists of typed variables. %

Let $\sC$ be a category with a monad $T = (T, \eta, \mu)$. Then $T$ is a {\em strong monad} if it is equipped with a natural transformation $t$ from the functor $(-) \times T(-):  \sC \times \sC \to \sC \times \sC$ to the functor $T(- \times -) :  \sC \times \sC \to \sC \times \sC$ (where $\times$ denotes both the product in Cat and the product in $\sC$). This $t$, called a {\em strength}, must additionally satisfy the identities:
$$(T \of \pi_{1,A}) \vertof t_{1, A} = \pi'_{1, TA},$$
$$(T \of \alpha_{A,B,C} ) \vertof t_{A \times B, C} = t_{A, B \times C} \vertof (1_A \times t_{B, C}) \vertof \alpha_{A, B, TC}$$
$$t_{A,B} \vertof (1_A \times \eta_B) = \eta_{A \times B}$$
$$t_{A,B} \vertof (1_A \times \mu_B) = \mu_{A \times B} \vertof (T \of t_{A,B}) \vertof t_{A, TB}$$
where notation is the same as in the preceding sections, except $\times$ denotes the product in $\sC$. 
%
%
In Chapter \ref{c.genm}, the fundamental parts of the theory of monads are extended to the setting of generalized categories in two ways, one via a generalized triple, and the other via a generalized Kleisli construction. 
Questions remain about how the present work is connected to Moggi's, since the Kleisli category in the generalized setting is a more subtle construction than in the one-categorical setting. 
It is also possible to extend our work in this paper to the setting of topos theory. We may define:


\dfn
An {\em ideal elementary topos} 
is an ideal cartesian closed category with 
\enu
	\item all finite limits and colimits,
	\item a subobject classifier.
\Enu
\Dfn

\noindent An investigation into topos theory in the generalized setting (including several sheaf theoretical constructions) is made in Chapter \ref{c.it}. 

\pagebreak
\singlespacing
\chapter{Ideal Toposes}\label{c.it}
\doublespacing
\vspace{10ex}

\section{Introduction}\label{s.it-i}

In Chapter \ref{c.th} it was shown that the framework of generalized categories may be used introduce an ideal of types to an ordinary category. 
Thus generalized categories (of which ordinary categories are a special case) can be given a structure possessing all the desirable properties of a cartesian closed structure, in spite of there being no way to put a true cartesian closed structure on a generalized category, since the mapping $- \times X$ 
is not a functor in general. The structure that arises in its stead is an ideal cartesian closed category: %
a generalized category equipped with cartesian structure and an evaluation operation, along with types and a set of constants. %
In this Chapter, we prepare the way for the study of topos theory in the generalized setting. %
We study both the elementary case and the case in which the topos structure arises from the existence of a site of definition. %

In section \ref{s.it} we define ideal elementary toposes and develop basic facts and constructions in that setting, following \cite{sga4,MaMo1,JoN1,FrD1,LaSc1}. An ideal topos is an ideal cartesian closed category equipped with finite limits, finite colimits, and a subobject classifier. After establishing preliminary results we show that the slice theorem can be generalized, and that for this reason, the generalized topos theory has, perhaps surprisingly, much of the same character as the usual one-categorical topos theory, arising from the Heyting algebra structure present and the rich system of adjoints generated by an arrow. %
A major feature of ordinary topos theory is the relationship with sheaf theory through which there arises a rich world of examples and applications. %
In section \ref{s.sht}, we develop a generalization of sheaf theory, and show that the connection between the categorical theory of sheaves developed by Grothendieck and his school \cite{sga4} and the theory of toposes developed by Lawvere, Tierney and others \cite{LaTi1} are in relation with one another in the generalized setting much as in the one-categorical case. 
An approach to sheaves that differs from the usual functorial approach is necessary in order to (in particular) carry %
the Yoneda embedding through in the generalized setting. 
The approach we adopt uses pointed profunctors. 


Aside from yielding information about the extent of applicability of topos theory, a result shown by our work is that it is often easier to work in the generalized setting than in the one-categorical setting. In part this may be due to the same phenomenon that occurs with other kinds of generalizations: in the general setting there are sometimes fewer possible directions for an argument to go in, thus reasoning in the general setting is psychologically more straightforward than in settings that arise from applications. 

The style in which this chapter is written is a departure from the style of most works in category theory. In part, this is a response to the needs of the subject matter. The author, much of whose training in mathematics was in analysis, %
has experienced first hand the reality that %
category theory is far better appreciated and far better understood in some circles than in others. %
He believes that an understanding and appreciation of the power and concreteness of categorical methods is possible among all mathematicians. %
There is, in the author's view, more than enough room for a wide variety of stylistic approaches to category theory.


\section{Basic Ideal Topos Theory}\label{s.it}

In a cartesian closed category, corresponding to a fragment of intuitionistic logic, there is one truth value, $\trut$, possibly accompanied by a symbol for absurdity $\bottom$. %
In a topos, this single truth value $\trut$ is expanded via a monic $\trut \armonic \Obsoc$. %
The intuition for the role of this monic, coming from sheaf theory, is {\em truth-in-place}, or localized truth. %
This means that what is true in a topos $\sE$ (or in the internal language of $\sE$) depends in general both on {\em location} in the topos and on {\em how far one can see} in the topos. Such an intuition might be a useful guide as we build the formalism.

\dfn\label{d.soc}
A {\em subobject classifier} for a generalized category $\sE$ is a pair $(\arsoc, \Obsoc)$ where $\Obsoc$ is an object in $\sE$ and where $\arsoc: \trut \armonic \Omega$ is in $\sE$, such that for every monic $m:a \armonic b$ in $\sE$, there is a unique arrow $p: b \to \Obsoc$ such that the following square is a pullback:
\[
\begin{tikzcd}
a \arrow[tail]{d}[swap]{m} \arrow{r}{\arter} 		& \trut \arrow[tail]{d}{\arsoc} \\
b \arrow{r}[swap]{p} 					& \Obsoc 
\end{tikzcd}
\]
\Dfn

\dfn\label{d.it}
An {\em ideal elementary topos} is an ideal cartesian closed category $\sE$ with
\enu
	\item all finite limits and colimits,\footnote{As in the one-category case, it suffices to have a terminal/equalizers/products, or a terminal/pullbacks (dually for colimits) \cite{ScMi}.} 
	\item a subobject classifier. 
\Enu
\Dfn



An ideal elementary topos $\sE$ has equalizers: for every pair of elements $f,g : a \to b$ in $\sE$, there is an arrow $e:e_0 \to a$ which has the property that $f \cdot e = g \cdot e$, and which is universal with respect to this property \cite{MacCW}. 
An equalizer is denoted $\Eq(f,g)$. This notation denotes an arrow well-defined in $\sE$ only up to isomorphism. 
It is also convenient to define, for element $a \in \sE$,  
$$\arsoc_a \threeline \arsoc \cdot \arter_a.$$

Because an ideal elementary topos is an ideal cartesian closed category, it carries the bijection:
\begin{equation}\label{eq.cccbij}
\hom(a,b) \bij \hom(1, a \vdash b).
\end{equation}
The element on the right, in terms of the element on the left, is denoted \cite{ScM3,LaSc1} 
$$\name{f}$$
This is called the {\em name} of $f$. 
%
In the case $b = \Omega$ the bijection (\ref{eq.cccbij}) becomes:
\begin{equation}\label{eq.namebijpower}
\hom(a,\Obsoc) \bij \hom(1, a \vdash \Omega).
\end{equation}
The element $a \vdash \Omega$ is viewed as a power of $a$ (cf. Theorem \ref{t.sliceT}). %
For example, in the category of sets, it reduces to the set of maps from a set $A$ to the two-element set $\set{0,1}$. %
The bijection (\ref{eq.namebijpower}) is augmented in a topos. Let $\Subobj (a)$ denote the subobjects of an element $a \in \sE$, that is, the monics of target $a$ taken only up to isomorphisms between sources \cite{MacCW}. 

\prop\label{p.predbij}
For every $a \in \sE$, there is a bijection
\begin{equation}\label{eq.predbij}
\Subobj (a) \bij \hom(a, \Obsoc) 
\end{equation}
\Prop
\prf
Consider the map which sends a monic $m$ to the element in the lower row of the diagram of Definition \ref{d.soc}. %
It lifts to a well-defined map on subobjects, by a simple diagram chase. %
It is equally clear that it is injective, and if $\phi:a \to \Obsoc$ is given, then since $\arsoc$ is monic (Definition \ref{d.it}), the pullback of $\phi$ and $\arsoc$ is also monic, so it is also surjective. %
\Prf

The bijections (\ref{eq.namebijpower},\ref{eq.predbij}) are natural in $a$. %
A {\em predicate} in $\sE$ is an element $\phi$ in an ideal elementary topos $\sE$ with target $\Obsoc$. As in the proof of Proposition \ref{p.predbij}, we let 
$$\pred(m), \quad \pred([m]), \quad \pred(s)$$
denote the predicate associated to monics $m: \bar m \rightarrowtail \hat m$, subobjects $[m]$, and elements $s: (a \vdash \Omega) \to \trut$ via Proposition \ref{p.predbij} and Equation \ref{eq.namebijpower} (it does not lead to any confusion to overload the symbol $\pred$ in this way). The predicate $\pred(m)$ is also called the {\em characteristic function}, the {\em characteristic map}, or the {\em classifying map} of $m$, or of the subobject defined by $m$, or of $s: 1 \to (a \vdash \Omega)$. %
Given a predicate, we can also apply the bijection (\ref{eq.predbij}) in the other direction to obtain a subobject from a predicate. 
In general, however, we cannot obtain a canonical monic given a predicate. We can do so in some cases, for example, in the category of sets we may use subsets. Analogous canonical monics are given more generally in generalized presheaf categories by subpointed profunctors (see section \ref{s.sht}). %
For many purposes, however, such monics, though appealing, are not really necessary. 

Because we will make limited use of elementary ideal toposes
, we introduce only Lambek's {\em strict} logical morphisms \cite{LaSc1}: 

\dfn\label{d.catidealtop}
A {\em morphism} $F:\sE \to \sE'$ of ideal elementary toposes is a functor \cite{ScM3} of ideal cartesian closed categories $F:\sE \to \sE'$ such that
such that 
$$F(\Obsoc) = \Obsoc,$$
$$F(\pred(m)) = \pred(F(m))$$
for all monics $m$ in $\sE$. %
This gives a category $\idealTop$ of elementary ideal toposes. %
\Dfn



\prop\label{p.epimoniciso}
Let $\sE$ be an ideal elementary topos. Then every monic in $\sE$ is an equalizer, and every epi monic is an isomorphism. 
\Prop
\prf
Because $\sE$ has a subobject classifier, every monic in $\sE$ is an equalizer. Indeed, 
$$\pred(m) \cdot m = \arsoc \cdot \arter_a = \arsoc \cdot \arter_b \cdot m,$$
and from this observation it is not difficult to see that a monic $m:a \to b$ is the equalizer of $\pred(m)$ and $\arsoc_b$. %
Therefore a monic is an equalizer of, say, $f$ and $g$. If it is also epi, then it follows that $f = g$. Therefore the identity $1_{\hat m}$ also has the equalizing property, giving a universal arrow $k: \hat m \to \bar m$. So $k\cdot m = 1_{\hat m}$. Moreover $m \cdot k = 1_{\bar m}$, since 
$$m \cdot k \cdot m = m = m \cdot 1_{\bar m}$$
and $m$ is monic. %
\Prf

Let $a \in \sE$. Let $\Delta_a$ be the diagonal map
$$\Delta_a := \langle 1_a, 1_a \rangle.$$
This map is monic, hence it has an associated predicate
$$\delta_a := \pred (\Delta_a).$$
A functor of ideal cartesian closed categories $F: \sE \to \sE'$ between ideal elementary toposes preserves $\pred()$ if and only if it preserves the $\delta$:
$$F(\delta_a) = \delta_{F(a)},$$
for all $a \in \sE$. %
A pullback diagram
$$
\begin{tikzcd}
\cdot \arrow[r, "\pullb'_{f,g}"] \arrow[d, swap, "\pullb_{f,g}"] & 
b \arrow[d, "g"] 
\\
a \arrow[r, "f"] &
c
\end{tikzcd}
$$
gives rise to another pullback diagram
$$
\begin{tikzcd}
\cdot \arrow[d, swap, "\langle \pullb_{f,g} {,} \pullb'_{f,g} \rangle"] \arrow[r, "f \cdot \pullb_{f,g} 
	"]  & 
c \arrow[d, "\Delta_c"]
\\
a \wedge b \arrow[r, "f \wedge g"] &
c \wedge c,
\end{tikzcd}
$$
and conversely. Here (as in 
previous chapters) $a \wedge b$ denotes the product of $a$ and $b$. 

We can also define 
$$\set{\cdot}_a := (\delta_a)^*,$$
where $()^*$ is the exponential transpose in $\sE$. %
We also take the predicate
$$\sigma_a := \pred(\set{\cdot}_a)$$
once we show that: 

\prop\label{p.singletonmonic}
$\set{\cdot}_a$ is monic.
\Prop
\prf
Suppose that $\set{\cdot}_a \cdot b = \set{\cdot}_a \cdot b'$. Then, in the notation of \cite{ScM3},
\[
(\delta_a \cdot \langle b \cdot \pi_{\bar b, a} , \pi'_{\bar b, a} \rangle )^* = (\delta_a \cdot \langle b' \pi_{\bar b, a} , \pi'_{\bar b, a} \rangle )^* 
\]
so
\[
\delta_a \cdot \langle b \cdot \pi_{\bar b, a} , \pi'_{\bar b, a} \rangle = \delta_a \cdot \langle b' \pi_{\bar b, a} , \pi'_{\bar b, a} \rangle.
\]
So 
\begin{equation}\label{eq.deltawedge1}
\delta_{a} \cdot (b \wedge 1_a) = \delta_a \cdot (b' \wedge 1_a).
\end{equation}
In the commutative diagram
$$
\begin{tikzcd}
\bar b \arrow[r, "b"] \arrow[d, swap, "\langle 1_{\bar b} {,} b \rangle"] & 
a \arrow[r, "ter_{a}"] \arrow[d, "\Delta_{a}"] & 
\trut \arrow[d, "\arsoc"] 
\\
\bar b \wedge a \arrow[r, "b \wedge 1_{a}"] &
a \wedge a \arrow[r, "\delta_{a}"] &
\Obsoc
\end{tikzcd}
$$
both squares are pullbacks. %
Hence the square formed by combining the two squares is also a pullback square. %
By equation (\ref{eq.deltawedge1}), 
the pair $(\langle 1_{\bar b}, b \rangle, ter_a \cdot b)$ is also a pullback for the diagram in which the bottom row is replaced with 
$\bar b \wedge a \overset{b' \wedge 1_a}{\to} a \wedge a \overset{\delta_a}{\to} \Obsoc$. %
So there is an invertible $\theta$ such that $\langle 1_{\bar b}, b \rangle \cdot \theta = \langle 1_{\bar b}, b' \rangle$. %
So $\langle \theta, b \cdot \theta \rangle = \langle 1_{\bar b}, b' \rangle$. %
So $\theta = 1_{\bar b}$ and $b \cdot \theta = b'$, and so $b = b'$, as desired. 
\Prf

Hence $\sigma_a$ is well-defined. The predicates $\sigma_a$ and $\delta_a$ may be interpreted, respectively, 
as tests: 
$$\text{$f: c \to (a \vdash \Omega)$ is a singleton} \quad \becomes \quad \sigma_a \cdot f = \trut_{c} $$
$$\text{$b,c :d \to a$ are equal} \quad \becomes \quad \delta_a \langle b,c \rangle = \trut_{d} $$
These interpretations can be extended to polynomials over $\sE$ \cite{ScM3}. 
In the latter case, the interpretation is justified in the formalism of $\sE$:

\prop
Let $\sE$ be an ideal elementary topos, and let $b,c: d \to a$ in $\sE$. Then 
$$b = c \eifaoif \delta_a \langle b,c \rangle = \trut_d.$$
\Prop
\prf
We follow the proof found in \cite{LaSc1}, which applies without change in the generalized setting. 
If $b = c$, then 
\begin{align*}
\delta_a \cdot \langle b, c \rangle 
&= \delta_a \cdot \langle b, b \rangle \\
&= \delta_a \langle 1_a, 1_a \rangle b \\
&= \arsoc_a \cdot b = \arsoc_{\bar b}.
\end{align*}
Conversely, if $\delta_a \cdot \langle b,c \rangle = \arsoc_a$, 
then the pair $(\langle b, c \rangle, ter_a)$ form a cone on the diagram 
$a \wedge a \overset{\delta_a}{\to} \Obsoc \overset{\arsoc}{\leftarrow} \trut$. 
But this diagram has pullback $a$. So there is a unique $k: d \to a$ such that 
$$
\langle b,c \rangle = \Delta_a \cdot k = \langle k, k \rangle.
$$
By taking projections, we have $b = c = k$. 
\Prf

Another important result, supplying monic ``images'' in toposes, continues to hold in the generalized setting:

\prop\label{p.imagee}
Let $\sE$ be an ideal elementary topos. For every $f \in \sE$, there is a factorization of $f$
$$f = m \cdot e$$
where $e$ is epi and $m$ is a monic that is universal in the sense that for every monic $\tilde m$ such that there exists $e$ such that $f = \tilde m \cdot e$, there exists $e'$ such that $m = \tilde m \cdot e'$. 
\Prop
\prf
Using the existence of finite colimits, take the pushout of the diagram $\hat f \overset{f}{\leftarrow} \bar f \overset{f}{\to} \hat f$ and call the injections $x$ and $y$. 
These form a diagram 
$
\begin{tikzcd}
\hat f \arrow[r, shift right, swap, "y"] \arrow[r, shift left, "x"] &
\cdot
\end{tikzcd}
$; 
let $m$ denote their equalizer. 
Since $f$ has the equalizing property (i.e., $x \cdot f = y \cdot f$), $f$ factors through $m$. 
Call this factorization $f = m \cdot e$. 
See the figures:
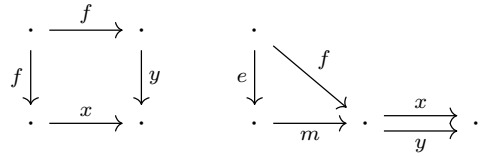
\[
\begin{tikzcd}
\cdot \arrow[r, "f"] \arrow[d, swap, "f"] &
\cdot \arrow[d, "y"] \\
\cdot \arrow[r, "x"] &
\cdot 
\end{tikzcd}
\quad\quad
\begin{tikzcd}
\cdot \arrow[d, swap, "e"] \arrow[dr, "f"] \\
\cdot \arrow[r, swap, "m"] &
\cdot \arrow[r, shift right, swap, "y"] \arrow[r, shift left, "x"] &
\cdot 
\end{tikzcd}
\]
Say that $m$ has the {\em image property} if
$$\text{for all } \tilde m, \tilde e, f = \tilde m \cdot \tilde e, \tilde m \text{ monic} \eimplies \text{ there exists } u \text{ such that } m = \tilde m \cdot u.$$
We claim that $m$ has the image property. Indeed, if 
$f = \tilde m \cdot \tilde e$, $\tilde m$ monic, then by Proposition \ref{p.epimoniciso} $\tilde m$ is an equalizer, say, of $s$ and $t$ in $\sE$. But if 
$s \cdot \tilde m = t \cdot \tilde m$, then 
$s \cdot \tilde m \cdot \tilde e = t \cdot \tilde m \cdot \tilde e$, so $s \cdot f = t \cdot f$. 
This gives a unique $u$ such that $s = u \cdot x, t = u \cdot y$, see figures:
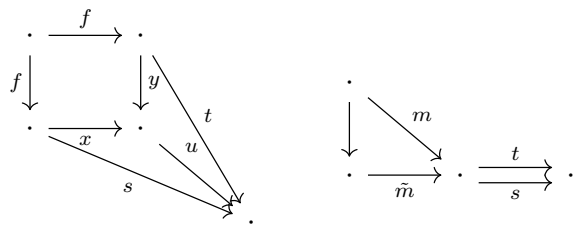
\[
\begin{tikzcd}
\cdot \arrow[r, "f"] \arrow[d, swap, "f"] &
\cdot \arrow[d, "y"] \arrow[ddr, "t"] \\
\cdot \arrow[r, swap, "x"] \arrow[drr, swap, "s"] &
\cdot \arrow[dr, "u" near start] \\
&
&
\cdot 
\end{tikzcd}
\quad\quad
\begin{tikzcd}
\cdot \arrow[d, swap] \arrow[dr, "m"] \\
\cdot \arrow[r, swap, "\tilde m"] &
\cdot \arrow[r, shift right, swap, "s"] \arrow[r, shift left, "t"] &
\cdot 
\end{tikzcd}
\]
Therefore $u \cdot x \cdot m = u \cdot y \cdot m$, so $s \cdot m = t \cdot m$,
so $m$ is an imposter cone for $\tilde m$. So $m$ factors through $\tilde m$, as desired. 

Now $m$ is monic since it is an equalizer; it remains to show that $e$ is epi. 
Suppose that we repeat the construction 
(equalizer of pushout) on the arrow $e$ instead of $f$; %
this yields a factorization $e = m' \cdot e'$ of $e$. The composition 
$m \cdot m'$ is monic, and $f$ factors through it, hence $m$ factors through it 
as well, say, $m = m \cdot m' \cdot v$. Therefore $1_{\hat{m'}} = m' \cdot v$. 
Therefore the monic $m'$ is epi, hence $m'$ is an isomorphism, 
again by using Proposition \ref{p.epimoniciso}. 
Now if the equalizer of the pushout of $\cdot \overset{e'}{\leftarrow} \cdot \overset{e'}{\to} \cdot$ is an isomorphism, then the projections of the pushout are equal. %
It can be checked that $e'$ is therefore epi. %
It follows immediately that $e$ is epi as well, as was to be proved. 
\Prf

The $m$ of Proposition \ref{p.imagee} (the equalizer of the pushout of $\cdot \overset{f}{\leftarrow} \cdot \overset{f}{\to} \cdot$) is called the {\em image} of $f$ and we denote it by $\im(f)$, following \cite{FrD1}. 

\cor\label{c.imagee}
The subobjects below $a \in \sE$ form a lattice. 
\Cor
\prf
For monics $m,n$ below $a$, define:
$$[m] \cap [n] := [n \cdot \pullb_{m,n}],$$
$$[m] \cup [n] := [\im [m,n]],$$
where $[m,n]$ is the universal arrow with respect to the coproduct $\bar m + \bar n$, following \cite{LaSc1}. %
These operations are well-defined, and provide a greatest lower bound (least upper bound, respectively). %
\Prf

The lattice $\Sub(a)$ has maximal element given by $1_a: a \to a$. In fact, in toposes $\Sub(a)$ has still more order-theoretic structure, but this takes time to prove (cf. Corollary \ref{c.sliceT}). 
Moreover, there is a morphism of lattices (i.e., an order-preserving map) $\Sub(a) \to \Sub(b)$ given an arrow $k: b \to a$, given by the pullback
$$k^{-1}([m]) := [\pullb_{k,m}].$$
Next we make the following observation:

\prop\label{p.subterob}
Let $a$ be an element of an ideal elementary topos $\sE$. 
The following are equivalent: 
 \enu
 	\item The terminal arrow $ter_a: a \to \trut$ is monic. 
	\item For all $b \in \sE$, there is at most arrow $b \to a$. 
\Enu
\Prop

If $a \in \sE$ has the property of Proposition \ref{p.subterob}, it is said to be {\em open}.

A result 
that more than any other single result gives a special character to toposes among categories says that the slice category of a topos is again a topos. 
The ``Fundamental Theorem of Topos Theory'' \cite{FrD1}, also called the Slice Theorem, is next to be proved in the generalized setting. 
%
We first introduce the slice category, and reintroduce subobjects, this time as a subcategory of the slice category. 


\dfn\label{d.subandslice}
Let $b$ be an element in an ideal elementary topos $\sE$. We define the {\em slice category} above or below $b$ to be the set of all triples 
$$(f,s,t)$$
that satisfy: $s: \bar s \to b$, $t: \bar t \to b$, $f: \bar s \to \bar t$, and $t \cdot f = s$. 
Composition and source/target are defined by:
$$(g,t,r) \cdot (f,s,t) := (g \cdot f,s,r),$$
$$\overline{(f,s,t)} := (1_{\bar s}, s, s),$$
$$\widehat{(f,s,t)} := (1_{\bar t}, t, t).$$
These operations define the structure of a one-category on $\sE/b$. 
\Dfn

Let $b$ be an element of a generalized category $\sE$. 
Consider the set of all monics with target $b$, along with a choice of map $f: m_1 \to m_2$ between each pair $m_1, m_2$ of monics with target $b$, whenever such a choice can be made (that is, whenever $\hom(m_1, m_2)$ is nonempty). In this way, the monics are taken as objects of a preorder-category and a subcategory of the slice category, denoted $\Sub_\sE (b)$ and called (herein) the {\em subobject category}. %
We must be somewhat careful to keep in mind that the (semi-)lattice structure of Corollary \ref{c.imagee} is only defined up to isomorphism in the subobject {\em category} defined in this section. As usual, this does not lead to any difficulty. 

The slice category is the category-theoretic analog of descending sets in preorders. In that case, there is nothing that distinguishes the subobject category from the slice category, though in general they are distinct (for example even in the case of the category of sets). 

\prop\label{p.subandslice}
Let $\sE$ be an ideal category or generalized category, and let $b \in \sE$.
\enu
	\item $\sE/b$ has terminal object $1_{\sE/b} := (1_b, 1_b, 1_b)$.
	\item If $\sE$ has finite limits, then $\sE/b$ has finite limits. 
	\item If $\sE$ has products and a subobject classifier, then $\sE/b$ has a subobject classifier. 
	\item $(m,s,t)$ is monic in $\sE/b$ if and only if $m$ is monic in $\sE$. \label{i.monicslice}
	\item The inclusion 
		$$i: \Sub(b) \hookrightarrow \sE/b$$
		of the subobject category in the slice category has left adjoint given by $\im(-)$. 
	\item There is an isomorphism of lattices
		$$\Sub_{\sE} (b) \bij \Sub_{\sE/b} (1_{\sE/b}).$$ \label{i.latticeiso}
\Enu
\Prop
\prf
(\ref{i.latticeiso}) follows from (\ref{i.monicslice}). %
The product arrows in $\sE/b$ are
$$\pi_{(1,s,s), (1,t,t)} := (\pullb_{s,t}, t \cdot \pullb_{s,t}, s)$$
$$\pi'_{(1,s,s), (1,t,t)} := (\pullb'_{s,t}, t \cdot \pullb'_{s,t}, s)$$
Also set
$$\Omega_{\sE/b} := (1_{\Omega \wedge b}, \pi'_{\Omega, b}, \pi'_{\Omega, b}),$$
\[
\arsoc_{\sE/b} := (\langle \delta_b \cdot \Delta_b, 1_b \rangle, 1_b, \pi'_{\Omega,b}). \qedhere
\]
\Prf

Note that the Grothendieck construction, or ``category of elements'' of a presheaf, is also a slice category. 
More precisely, using the language of section \ref{s.sht}, let $(P, \S_P)$ be a pointed profunctor over a generalized category $\sC$. Then the Grothendieck construction $\el(P)$ is $P/\S_P$, that is, it is the one-category whose underlying set is the set of triples $(f,s,t)$ with $f \in \sC$, $s,t$ in $P$ with target $\S_P$ (i.e., sections), and satisfying $f: \bar s \to \bar t$, and $t \cdot f = s$. Composition in $\el(P)$ is given by 
$$(g,t,r) \cdot (f,s,t) = (g \cdot f, s, r).$$
and source and target $\overline{(f,s,t)} = (1_{\bar s}, s, s)$, $\widehat{(f,s,t)} = (1_{\bar t}, t, t)$. 
As with any slice, the map $(f, s, t) \mapsto f$ defines a functor $\pi = \pi_P$ from $\el(P)$ to $\sC$:
$$\pi(f,s,t) := f.$$
Indeed, 
$$\pi( \overline{(f,s,t)}) = \pi((1_{\bar s}, s, s)) = \overline{1_{\bar s}} = \bar s = \bar f = \overline{\pi((f,s,t))},$$
and similarly for target/composition. 


The ``sub and slice'' structure, i.e., the structure of the subobject lattice/category and the slice category, take on their usual role familiar from the one-categorical case. However, it is even more fitting in the generalized situation than in the one-categorical situation to call the slice category a ``slice'': it is a one-dimensional structure within the generalized category (by ``one-dimensional'' we mean it is a one-category). This one-category has in it 
the subobject lattice (a true lattice, for a topos) represented as a system of monics. 
The fundamental result is: 

\thm\label{t.sliceT}
Let $\sE$ be an ideal elementary topos, $b \in \sE$. Then the slice category 
$\sE/b$ 
is a cartesian closed category. 
\Thm
\prf
First observe that, unlike in the one-categorical case, the power object operator $P(-)$ is not a functor. Rather, it is given by two separate operators, which we denote by $P(-)$ and $\sP(-)$, defined for $a \in \sE$ by
$$P(a) := a \vdash \Omega,$$
$$\sP(a) := (\varepsilon_{\Omega, \hat a} \langle \pi_{\hat a \vdash \Omega, \bar a}, a \cdot \pi'_{\hat a \vdash \Omega, \bar a} \rangle )^*.$$
Then $\sP(a): a \to \Omega$, and $\sP(a) : P(\hat a) \to P(\bar a).$ Moreover $\sP(a)$ is the unique element in $\sE$ that satisfies (omitting the subscripts of projections)
$$\varepsilon_{\Omega, \bar a} \cdot \langle \sP(a) \cdot \pi, \pi' \rangle = \varepsilon_{\Omega, \hat a} \cdot \langle \pi, a \cdot \pi' \rangle.$$
or what is the same again
$$ \varepsilon_{\Omega, \bar a} \cdot (\sP(a) \wedge 1) = \varepsilon_{\Omega, \hat a} \cdot (a \wedge 1). $$
The rest of the proof is to verify that this change in the nature of the $P(-)$ operator does not invalidate the proof in the one-categorical case. 

Now we wish to prove that there is a natural bijection between two hom sets, as in:
\begin{equation}\label{e.powerslice}
\hom_{\sE/b}(\fs \wedge \ft, \Omega_{\sE/b} ) \bij \hom_{\sE/b} (\ft, P \fs) 
\end{equation}
where (provisionally) we use the letters $\fs, \ft,$ etc. for objects $(1_{\bar a},a,a)$ of the slice, and where $P(-)$, the power object operator in the category $\sE/b$, is still to be defined. %
In order to do this, we convert hom sets to subobject lattices, as we may always do in toposes. 
For the rest of the proof of (\ref{e.powerslice}), let $\fs = (1_{\bar s}, s, s)$ and $\ft = (1_{\bar t}, t, t)$, for $s,t \in \sE$ with target $b$. %

First, $\fs \wedge \ft$ is the pullback $(1, t \pullb_{s,t}, t \pullb_{s,t})$ with source $\bar s \times_b \bar t$. %
By Proposition \ref{p.subandslice}, $\Omega_{\sE/b}$ is $(1_{\Omega \wedge b}, \pi', \pi')$. 
Therefore, by inspection, we have a bijection
$$\hom_{\sE/b} (\fs \wedge \ft, \Omega_{\sE/b}) \bij \hom_\sE ( \bar s \times_b \bar t , \Omega)$$
that is natural in $\fs, \ft$. 
Now in the topos $\sE$ we use the natural bijection 
$$\hom(\bar s \times_b \bar t, \Omega) \bij \Sub(\bar s \times_b \bar t).$$
The fiber product $\bar s \times_b \bar t$ is a subobject of $\bar s \wedge \bar t$, whence we obtain
$$\Sub(\bar s \times_b \bar t) \bij \set{[m] \mid [m] \text{ is a subobject of } \bar s \wedge \bar t \text{ and } [m] \lies [ s \cdot \pullb_{s,t}]}$$
Next we translate to predicates and obtain 
$$\hom_{\sE/b}(\fs \wedge \ft, \Omega_{\sE/b} ) \bij \set{h: \bar s \wedge \bar t \to \Omega \mid h \wedge_b \pred( \langle \pullb_{s,t} , \pullb_{s,t}' \rangle ) = h}$$
where here we write an infix operator $\wedge_b$ for the arrow $\hom(b, \Omega \wedge \Omega) \to \hom(b, \Omega)$ 
arising from the meet $\cap$ in $\Sub(b)$ via $\hom(b, \Omega) \iso \Sub(b)$, 
and where we use the standard expression of the fiber product as a subobject of the product:
$$\bar s \times_b \bar t \overset{\langle \pullb_{s,t} , \pullb_{s,t}' \rangle}{\armonic} \bar s \wedge \bar t.$$
Now we take star (transpose) of the $h$'s to obtain 
$$\hom_{\sE/b}(\fs \wedge \ft, \Omega_{\sE/b}) \bij \set{k: \bar t \to P(\bar s) \mid k \wedge_{int} \sP(s) \cdot \set{\cdot}_b \cdot t = k}$$
where $\wedge_{int}$ (``internal meet'') denotes the infix operator defined by the arrow $P(b) \wedge P(b) \to P(b)$ again arising from the meet $\cap$ in $\Sub(b)$, 
and where we have used the fact (result of inspection) that 
$$\pred(\langle \pullb_{s,t} , \pullb'_{s,t} \rangle) = \varepsilon_{\Omega, c} (1 \wedge (\sP(s) \cdot \set{\cdot}_b \cdot t)).$$
This now becomes (writing $\wedge_{int}$ for the internal meet arrow, not the infix operator) 
$$\hom_{\sE/b}(\fs \wedge \ft, \Omega_{\sE/b}) \bij \set{k: \bar t \to P(\bar s) \mid \wedge_{int} \cdot (1 \wedge (\sP(s) \cdot \set{\cdot}_b \cdot t)) \cdot \langle k,t \rangle = \pi_{k,t} \cdot \langle k,t \rangle}$$
Now let 
$$t_1 := \wedge_{int} \cdot (1 \wedge (\sP(s) \cdot \set{\cdot}_b \cdot t).$$
and define, for objects $(1,s,s)$ in $\sE/b$,
$$P_{\sE/b} ( (1,s,s)) := (1, \pi'_{P(s), b} \cdot \Eq(t_1, \pi_{k,t}), \pi'_{P(s), b} \cdot \Eq(t_1, \pi_{k,t}) ),$$
we have the desired bijection (\ref{e.powerslice}) using the universal property of the equalizer $\Eq(t_1, \pi_{k,t})$ of $t_1$ and the first projection $\pi_{k,t}$. 
%
%
\Prf

We now make note of some of the most significant consequences of the slice theorem (for more see, for example, \cite{JoN1,MaMo1}). The proofs of these results now follow closely the one-categorical case. 

\cor\label{c.sliceT}
Let $\sE$ be an ideal elementary topos. 
\enu
	\item The slice $\sE/b$ over an element $b \in \sE$ is an elementary topos. \label{i.slicetopos}
	\item The subobject category $\Sub_a (\sE)$ below every element $a \in \sE$ is (up to isomorphism) a Heyting algebra. \label{i.subheyting}
	\item If $k: b \to a$, then the arrow $k^*: \sE/a \to \sE/b$ given by 
		$$k^*(f,s,t) := (\langle t \cdot \pullb_{k,t} , \pullb'_{k,t} \rangle_{s,k} , \pullb'_{k,s} , \pullb'_{k,t} )$$
		has a left adjoint $\Sigma_k$ and a right adjoint $\Pi_k$. 
	\label{i.pisigma}
\Enu
\Cor





There are many other consequences of Theorem \ref{t.sliceT}, see the references for others. %
In particular, the slice theorem is used in the proof of the theorem of Giraud \cite{sga4}, relating Grothendieck toposes to a list of conditions on an ordinary one-category. %
In order to extend this result it must be determined whether a generalized sheaf (to be defined below) is the colimit of representable ones. 


%

\section{Sheaves over Generalized Categories}\label{s.sht}

If we shunt away the task of providing motivation, it is remarkably easy, due to the work of Lawvere and Tierney, to introduce a notion of sheaf:

\dfn\label{d.latitop}
Let $\sE$ be an ideal elementary topos. A {\em Lawvere-Tierney topology,} or {\em topology} on $\sE$ is a choice of element $j \in \sE$ such that
\enu
	\item[(0)] $j : \Obsoc \to \Obsoc$,
	\item[(1)] $j \cdot \arsoc = \arsoc$,
	\item[(2)] $j \cdot j = j$,
	\item[(3)] $j \cdot \delta_{\arsoc} = \delta_{\arsoc} \cdot (j \wedge j).$
\Enu
A monic $m$ is {\em dense} with respect to $j$ if %
$$[m] = \subob (j\cdot \pred(m)),$$
where $[m]$ is the subobject defined by $m$ and $\subob(\phi)$ is the subobject defined by a predicate $\phi$. Dense subobjects are defined similarly. 
\Dfn

\dfn\label{d.sheaflatitop}
Let $(\sE,j)$ be an elementary ideal topos $\sE$ equipped with a Lawvere-Tierney topology $j$. %
A {\em sheaf} for $j$ on $\sE$ is an element $P$ of $\sE$ such that %
for every monic $m$ in $\sE$, and for every $f:\bar m \to P$ in $\sE$, if $m$ is dense, there exists a unique $g \in \sE$ such that %
$$g \cdot m = f.$$
\Dfn

To motivate these definitions, however, we must relate them to a generalization of geometric sheaf theory. %
This is the task for the rest of this section. 

\subsection{Pointed Profunctors}\label{ss.pprofunctor}
The profunctor abstraction \cite{BeU2} provides a convenient language for beginning a study of generalized sheaves. 

\dfn\label{d.pointedprofunctor}
A {\em pointed profunctor} is a structure $(P, \S)$ consisting of a generalized category $P$ equipped with a distinguished object $\S$ of $P$ satisfying: for all elements $f$ in $P$,
$$\dom(f) = \S \eimplies f = \id_\S.$$
\Dfn

The full generalized category $\sC$ on the objects $U \in P$, $U \neq \S$ is the {\em base category} of the pointed profunctor $P$. %
In this case $P$ is said to be defined {\em over $\sC$}. %
The element $\S$ is the {\em section object}. %
If $U$ is an element of $\sC$, the base category of $P$, the set of elements $\hom(U, \S)$ are called {\em sections over $U$} and the set $\hom(U, \S)$ is also called the {\em fiber over $U$}, denoted $\Fib_P (U)$ or $P(U)$. %

\dfn
A {\em morphism} of pointed profunctors $P$ and $Q$ is a functor $\phi: P \to Q$ satisfying: %
\enu
	\item $\phi(\S_P) = \S_Q,$
	\item for all elements $U \in \sC$, $\phi(U) = \S \eimplies U = \S.$
\Enu
If the induced map between base categories is the identity, then $F$ is said to be {\em base-preserving}. 
\Dfn








A pointed profunctor over a category is not a presheaf, because it may have fibers over morphisms in the base. 

\prop\label{p.presheafiso}
Let $\sC$ be a generalized category. 
The category of presheaves over $\sC$ is isomorphic to the category of pointed profunctors and base-preserving morphisms over $\sC$ with the property that
\begin{equation}\label{eq.presheafcondition}
\text{proper arrows of $\sC$ remain proper arrows of $P$.}
\end{equation}
\Prop
\prf
Fix a generalized category $\sC$ and let $(A,\S)$ be a pointed profunctor over $\sC$ with the property \ref{eq.presheafcondition}. %
Define %
$$F:\sC^{op} \to \text{Set}$$
by sending $a \in \sC$ to the map $s \mapsto sa$ from the set $\set{s \in A \mid \bar s = \hat a}$ to the set $\set{s \in A \mid \bar s = \bar a}$. Then we can check that $F(ab) = F(b) \of F(a),$ and $F(a)$ is an identity in the category of sets if and only if $a$ is an identity in $\sC$. %
So $F$ is a presheaf on $\sC$. 
Now let $F:\sC^{op} \to \text{Set}$ be a limit functor. For $U \in \Ob(\sC),$ let $\Sect(U)$ %
be the set which is the domain (and codomain) of $F(1_a)$, and let $\Sect(\sC) := \bigcup_{U \in \Ob(\sC)} \Sect(a)$.
Let 
$$P := \sC \sqcup \Sect(\sC) \sqcup \set{\S},$$
where $\sqcup$ is disjoint union, and extend the operations of $\sC$ to $P$: for $a \in \sC$, $s \in \Sect(\sC),$ 
\begin{align*}
	\S &\in \Ob(A),			\\
	\cod(s) &:= \S,			\\
	\dom(s) &:= a, \quad \call{s \in \Sect(a)}, \\
	s\cdot a &:= F(a)(s), \text{ if} \cod(s) = \dom(a), \text{ not def. otherwise.}
\end{align*}
Then $P$ is a category, and moreover a pointed profunctor over $\sC$. %



We thus obtain a presheaf, and this construction is clearly inverse to the preceding one. %
Now we have to show that this bijection $\Phi$ extends to maps. Let $A,B$ be two presheaves, and let $\phi:A \to B$ be a morphism of presheaves. Then because
$$\phi(sa) = \phi(s)\phi(a),$$
if we let $\Phi(P) := F$, $\Phi(Q) := G$, we have $\phi(F(a)(s)) = G(\phi(a)) (\phi(s)).$ %
We are given that $\phi$ leaves $\sC$ fixed, so we have, for all $s \in \dom(F(a))$,
$$(\phi \of F(a))(s) = (G(a) \of \phi) (s).$$
For $U \in \sC$, let $\phi_{nat}(U) = \phi|_{\Sect(U)}$. Then for all $a \in \sC$, 
$$\phi_{nat}(\hat a) \of F(a) = G(a) \of \phi_{nat} (\bar a) \,\downarrow.$$
So we have a natural transformation $F \Rightarrow G$.

Conversely, suppose $\phi_{nat}$ is a natural map for a natural transformation $\phi:F \to G$, where $F,G$ are two presheaves over $\sC$. Then for all $a \in \sC$, %
$$F(a): \Sect_F (\bar a) \to \Sect_F (\hat a),$$
$$G(a): \Sect_G (\bar a) \to \Sect_G (\hat a).$$
So $\phi_{nat} (a)$ is a function 
$$ \phi_{nat}(a):\Sect_F (a) \to \Sect_G (a).$$
We can therefore define a function
$$\phi: \Sect_F(\sC) \to \Sect_G (\sC)$$
by mapping $s$ (over $a$) to $\phi_{nat}(a)(s)$. Extend this map to $A$ by setting $\phi(U) = U$ for $U \in \sC$, and $\phi(\S) = \S$. %
We obtain a base-preserving morphism of pointed profunctors. %
This is $\Phi(\phi)$, and $\Phi$ is still invertible after extending to maps. %
Functoriality of $\Phi$ is easily verified. %
So we obtain the desired isomorphism. %
\Prf

This result can be interpretted as indication that the category of sets is not a natural player in the generalized setting. %
This result can be extended to the two-categorical level as well, where the same phenomenon recurs. 
We omit the details, but make mention of the following, whose proof is a more involved rendition of the steps in the preceding one. The bicategorical notion coinciding with pointed profunctor is a {\em coweighted bicategory}: 

\thm\label{t.bicat}
Let $\sC$ be a category. Then the category of (normalized) fibered categories and optransformations is isomorphic to the category of coweighted bicategories over the categorical base $\sC$ and base-preserving opmorphisms. 
\Thm

\subsection{Generalized Sheaves}\label{ss.gensheaf}
We have obtained a category of pointed profunctors in the preceding section. We could apply the construction of \cite{ScM3} to obtain an ideal cartesian closed category. %
We could then pursue the investigation further and ask: 
is it an ideal elementary topos? %
However, this line of inquiry is likely only to yield a dead letter, because there is no hope of establishing a Yoneda embedding in the generalized setting via this construction. %
The Yoneda embedding is a full-fledged functorial embedding of the base generalized category provided that there is a way to assign domains and codomains to presheaves (pointed profunctors). To develop such machinery, we must modify our approach. We turn to the language of bipartite categories. 


\dfn\label{d.bipartitegencat}
A {\em bipartite generalized category} is a generalized category $\sE$ equipped with a subcategory $\sC$, the {\em base} 
of the generalized category, such that every element $\rho$ of $\sE$ such that $\cod(\rho) \in \sC$ satisfies $\rho \in \sC$. 
\Dfn

We say that $\sE$ is {\em over the base $\sC$} in this case. 
The elements of $\sE$ not in $\sC$ whose target is in $\sC$ is called a {\em section} of $\sC$. 
The set of elements of $\sE$ that are neither base elements nor sections is called, collectively, the {\em upper region} of $\sE$. 


Let $\sE$ be a bipartite generalized category. It is easy to see that fixing a choice of element $P \in \sE$ in the upper region produces, by a direct construction, a pointed profunctor $\tilde P$ over the base of $\sE$. Moreover, a choice of element $P \in \sE$ in the upper region produces a morphism of pointed profunctors $\bar P : \tilde P_1 \to \tilde P_2$ from its source to its target. 

\dfn\label{d.maximalbipartitegencat}
A bipartite generalized category $\sE$ is {\em maximal} over its base if every construction of a pointed profunctor yields an element of $\sE$, and if every construction of a pointed profunctor homomorphism also produces an element of $\sE$. 
\Dfn

\prop\label{p.maxunique}
There exists a unique maximal bipartite generalized category over a given category $\sC$. 
\Prop
\prf
Suppose that there were two distinct maximal bipartite categories $\sE, \sE'$ over a base $\sC$. Then 
$$\sE'' := \sC \sqcup \Sect(\sE) \sqcup \Sect(\sE') \sqcup \UR(\sE) \sqcup \UR(\sE')$$
is a bipartite category over $\sC$ extending both, contradicting the maximality of $\sE$ and $\sE'$. 
\Prf

This proof is (classically) valid, assuming we work in a universe $\sU_{univ}$ closed under the usual set theoretical operations. %


\thm\label{l.yoneda}
Let $\sE$ be the maximal bipartite generalized category over a base generalized category $\sC$. Then there exists a functor $\yoneda:\sC \to \sE$ with the property that for every element $U \in \sC$ there is a bijection between the fiber of any element $P$ and the homomorphisms from $\yoneda (U)$ to $P$.
\Thm
\prf
The proof is by construction: we produce a copy of $\sC$ (an ``exact carbon copy'') along with all sections of $\sE$, and a third copy of each element $f: a \to b$ in $\sC$, changing the target of the third copy to the second copy of the target. By maximality, this system of pointed profunctors is already in $\sE$. 
\Prf

The same construction, carried out by copying the base and moving the copy to the upper region, gives an extension of any bipartite generalized category. %
In some sense, this is the meaning of Yoneda's lemma. %
We can go further on the basis we have now laid out, but instead, we shall continue on the basis of a notion of presheaf in the generalized setting that is in agreement with the notion of ideal elementary topos used in section \ref{s.it}. %
The reason we must shift is if $P \vdash Q$ is an ideal element \cite{ScM3} of the upper region of a bipartite category, %
the composition of $P \vdash Q$ with a section element yields a section element that must itself be ideal. %
Therefore we define: %

\dfn\label{d.imbgencat}
An {\em ideal maximal bipartite generalized category} $\sE$ is a bipartite category that satisfies:
\enu
	\item $\sE$ has the maximality property of Definition \ref{d.maximalbipartitegencat},
	\item For every $P$ in the upper region of $\sE$ and for every $U$ in the base of $\sE$, there is a section $s_U :U \to P$ that has the property that for all $f, g \in \sE$ that are not identities, if $s \cdot f$ $\downarrow$ then $s_U \cdot f = s_{\bar f}$, and if $g \cdot s$ $\downarrow$ then $g \cdot s = s_U$. 
\Enu
\Dfn

By the same argument as proved Proposition \ref{p.maxunique}, there is a unique ideal maximal bipartite generalized category over a base $\sC$. The $s_U$ of Definition \ref{d.imbgencat} is the ``ideal element'' in the fiber (cf. \cite{ScM3}). 

\dfn\label{d.gencatpresheaf}
Let $\sC$ be a generalized category. The ideal maximal bipartite generalized category $\sE$ over $\sC$ is called the {\em generalized category of generalized presheaves over $\sC$}. %
An element $P$ of the upper region of $\sE$ 
is called a {\em generalized presheaf} over $\sC$. 
\Dfn

Thus a generalized presheaf is never considered ``in a vacuum'', but always in the context of a system of presheaves.

Next, we have a result that gives a construction of an ideal cartesian closed category \cite{ScM3} from the data of a generalized category:

\thm\label{t.pshiccc}
The generalized presheaf category over a base $\sC$ 
is an ideal cartesian closed category. 
\Thm
\prf
We check the definition \cite{ScM3}. %
Let $\sE$ be the generalized presheaf category over $\sC$. %
It is clearly a generalized category, and we can for simplicity (not being concerned with any logical interpretation at present) assume that all elements of $\sE$ are valid. %
We have an element $\trut$ in $\sE$ given by the presheaf ${\bf pt}$ with itself as source and target and with a single element $pt$ in each fiber. %
We take the constants of $\sE$ to be the terminal arrows $\arter_P$ that send every section of a generalized presheaf $P$ to $pt$. %
(This element $pt$ must be the ideal section in each fiber.) %
The binary product operation $\times$, as noted above, is given by fiber-wise union and component-wise composition, recursively defined on sources and targets. %
The operation $\langle P,Q \rangle$ is defined when $\bar P = \bar Q$ by sending $s$ in $\bar P = \bar Q$ to $( P\cdot s, Q \cdot s )$. 

Now we can no longer put off the ideal element $P \vdash Q$ for given $P,Q$ in $\sE$. 
We must construct given $P$ and $Q$ an element 
$$P \vdash Q$$
sourced at $P$ and targeted at $Q$. Let $\S_{P \vdash Q}$ be a fresh symbol, and consider the set of arrows $s$,
$$s:P \times \yoneda(U) \to Q$$
For each such $s$ we define an element 
$$\tilde s : U \to \S_{P \vdash Q}.$$
This defines a pointed profunctor with section object $\S_{P \vdash Q}$. %
We make a generalized presheaf out of this pointed profunctor by setting %
$$s(P \vdash Q) = P,$$
$$t(P \vdash Q) = Q.$$
If $\phi: R \to P$, $\psi: Q \to S$ in $\sE$, we define 
$$(P \vdash Q) \cdot \phi := R \vdash Q,$$
$$\psi \cdot (P \vdash Q) := P \vdash S,$$
unless $\phi$ is a constant in $\sK$. In that case, $\phi$ is a terminal, hence $P$ is ${\bf pt}$. In this case, we may define $(P \vdash Q) \cdot \phi$ on sections to be be the constant map $s \mapsto s_U$, where $s_U$ is the ideal section, on each fiber. 
We can also check that if $(P,Q)$ is a pair of elements in $\sE$, then we have the projections $\pi_{P,Q}, \pi'{P,Q}: P \times Q \to P, P \times Q \to Q$, respectively, and $(\pi_{P,Q}, \pi'_{P,Q})$ satisfy the axioms for a good pair for $(P,Q)$. 


The operation $()^*$ must now be defined. Let $\phi: R \times P \to Q$ in $\sE$. We define a pointed profunctor homomorphism $\phi^*$ from $R$ to $P \vdash Q$ (which we identify with $\hom(\yoneda(-) \times P , Q)$) now containing an ideal element $\yoneda(-))$ by defining $\phi^*(s)$, for $s:U \to R$, to be the the map sending a section $(f,x)$ of $\yoneda(U) \times P$, where $f$ is identified with an element of $\hom(-,U)$, to be
$$\phi \cdot \langle s \cdot f, t \rangle.$$
This defines a pointed profunctor homomorphism. %
This establishes that $\sE$ defines a positive intuitionistic generalized deductive system. 
Before checking good evaluations, we pause to observe that there is a unique $f: a \to \trut$, hence $f$ must be $a \vdash \trut$. 
In particular, for $a = \trut$, this shows that $\trut \vdash \trut = \trut$. 
This shows that $\sK$ is closed under 
$\sK$ closed under $\wedge, \langle,\rangle,$ and $()^*$, as desired. %

It remains only to show that there are good evaluations. 
Let $P,Q$ be a pair of elements in $\sE$. Let $(\pi, \pi')$ be a good pair. Since good pairs are unique if they exist, these are the usual projections on $P \times Q$. 
For $U$ in $\sC$, for $\theta$ a section of $P \vdash Q$ over $U$, and $s$ a section of $P$ over $U$, define 
$$\epsilon_U (\theta, s) = \theta_U ( 1_U, s),$$
where we again identify a section $\theta$ of $P \vdash Q$ over $U$ with its associated element $\yoneda(U) \times P \to Q$. %
We pass to this morphism the identity $1_U$ on $U$ along with the given section $s$ of $P$. %
Finally, we check that Lambek's axioms
$$\epsilon \cdot \langle f^* \cdot \pi_{c,a}, \pi'_{c,a} \rangle = f,$$
$$(\epsilon \cdot \langle f \cdot \pi_{c,a}, \pi'_{c,a} \rangle)^* = f.$$
follow as a consequence of the definitions of $\epsilon$ and $()^*$. 
\Prf


We will prove that the generalized presheaves form an ideal elementary topos after introducing sieves. 



\dfn\label{d.sieve}
A {\em sieve} on a generalized category $\sC$ is a subset $\sS$ of $\sC$ such that there exists $U \in \sC$ such that for every $s \in \sS$, $\hat s = U$. We say that $\sS$ is a {\em sieve at $U$}. We write $\bar \sS = U$. 
\Dfn

The notation $\bar \sS$ is intended to evoke the definition $\sS: U \to \Omega$, where the target $\Omega$ will be defined below. Note the possible source of confusion: $\overline{\sS}$ is the common {\em target} of elements of $s$. 

\prop\label{p.sieve}
We have the following:
\enu
	\item The set $\set{f \in \sC \mid \hat f = U}$ is a sieve, called the {\em maximal sieve} at $U$. 
	\item If $\sS$ is a sieve at $U$ and $\rho$ is an element of $\sC$ with $\hat \rho = U$, the set
$$\stem{\sS}{\rho} := \set{\sigma \in \sC \mid \rho \cdot \sigma \downarrow \eand \rho \cdot \sigma \in \sS}$$ 
		is a sieve at $\bar \rho$. 
	\item If $\sS$ is sieve, $f \in \sC$, and $\hat f = \bar \sS$, then $f \in \sS$ if and only if $\stem{\sS}{f}$ is maximal. 
\Enu
\Prop

The order-theoretic notion corresponding to a sieve is a lower set, sometimes called a downset or descending set. %
%
The symbol $\stem{\sS}{\rho}$, denoted $\rho^*(\sS)$ in \cite{MaMo1}, denotes a generalization of the notion of restriction of a cover: $\set{U_\alpha}|_U := \set{U_\alpha \cap U}$, where $\set{U_\alpha}$ denotes an open cover in a topological space. 
Sieves are closely related to the notion of a subpointed profunctor: %

\dfn\label{d.subpointedprofunctor}
A {\em subpointed profunctor} $P'$ of a pointed profunctor $P$ on a generalized category $\sC$ is a generalized category of $P$ with the same section object $\S$ and the same base generalized category $\sC$. 
\Dfn

The subpointed profunctor is closed under composition, hence $s \cdot \rho$ is in $P'$ for every $\rho$ in $\sC$ and $s$ in $\Sect (P')$. 
It can be explicitly observed that the subpointed profunctors of a given pointed profunctor $P$ have the structure of a Heyting algebra, in fact, a frame. Indeed, the meet and join of subpointed profunctors $P_1$ and $P_2$ are given by 
$$P_1 \vee P_2 = P_1 \cup P_2,$$
$$P_1 \wedge P_2 = P_1 \cap P_2,$$
where the symbols $\cup$ and $\cap$ are used to denote the corresponding sets, along with the obvious pointed profunctor structure. %
The maximum element of the Heyting algebra is $P$ itself, and the minimum is the subpointed profunctor $(\sC, \S)$ with no sections. The difference or residual of $P_1$ and $P_2$ is
$$P_1 \Leftarrow P_2 = \set{s \in P_1 \mid \text{for all } f \in P, s \cdot f \,\downarrow, s \cdot f \in P_2 \eimplies s \cdot f \in P_1}.$$
This structure entails a {\em pseudocomplement} to any subpointed profunctor:
$$\lnot P := \set{s \in P \mid s \nin P, \eand s\cdot f \nin P \text{ for all } f \text{ such that } s \cdot f \,\downarrow}.$$
The data of a sieve at $U \in \sC$ is also given by a subpointed profunctor of $\yoneda(U)$. 

\dfn\label{d.gtop}
A {\em topology on $\sC$} is a collection $J$ of subsets in $\sC$, satisfying:
\enu
	\item Every $\sS$ in $J$ is a sieve in $\sC$.
	\item For every $U \in \sC$, the maximal sieve at $U$ is in $J$. 
	\item If $\sS$ is in $J$ and $\rho \in \sC$ with $\hat \rho = \bar{\sS}$, then $\stem{\sS}{\rho}$ is in $J$. 
	\item If $\sS$ is in $J$ and $\sR$ is any sieve satisfying $\bar \sR = \bar \sS$, 
	and if for every $f \in \sS$, $\stem{f}{\sR}$ is in $J$, then $\sR$ is in $J$. 
\Enu
Elements of the subcollection of $J$ consisting of sieves at $U \in \sC$ are called the sieves that {\em cover $U$} or {\em covering sieves} for $U$.  
A generalized category $\sC$ equipped with a topology $J$, denoted $\sC$ or $(\sC, J)$, is a {\em site}. 
\Dfn

\prop\label{p.gtop}
Let $(\sC, J)$ be a site. 
\enu
	\item If $\sS, \sS'$ are covering sieves at a common target $U \in \sC$, then the intersection $\sS \cap \sS'$ of $\sS$ and $\sS'$ (regarded as subsets of $\sE$) is a covering sieve. 
	\item If $\sS$ is a covering sieve and if $\sR$ is a sieve such that $\sR \rise \sS$ (that is, $\sR$ is larger than $\sS$), then $\sR$ is a covering sieve. 
	\item If $\sS$ is a covering sieve and if for every $f \in \sS$ there is given a covering sieve $\sS_f$ at $\bar f$,
	then the ``sum''
	$$\underset{f \in \sS}{\fS} \sS_f := \set{f \cdot g \mid f \in \sS, g \in \sS_f}$$
	is a covering sieve. 
\Enu
\Prop







We can repeat the same constructions as usual, leading to well-known examples of topologies \cite{sga4,MaMo1}. 

\dfn\label{d.closedsieve}
Let $(\sC, J)$ be a site. A {\em closed sieve} $\sS$ on $\sC$ is a sieve $\sS$ with the property that for all elements $f \in \sC$,
$$\stem{\sS}{f} \in J \eimplies \stem{\sS}{f} \text{ is the maximal sieve}. $$
In particular, $\sS$ is maximal, thus to be closed is somehow to be super-maximal. 
The {\em closure} of a sieve $\sS$ is the sieve
$$\closure{\sS} := \set{f \in \sC \mid \hat f = \overline{\sS} \eand \stem{\sS}{f} \text{ is in } J}.$$
\Dfn

\prop\label{p.closedsieve} $\phantom{v}$
\enu
	\item The closure of a sieve $\sS$ is a closed sieve. The closure of a sieve is the smallest closed sieve containing it. 
	\item If $\sS$ is closed, then for all composable $f$, $\stem{\sS}{f}$ is closed. 
	\item $\stem{\closure(\sS)}{f} = \closure(\stem{\sS}{f})$.
	\item A closed covering sieve is maximal, and conversely. 
\Enu
\Prop

Closed sieves need not be covers, and covering sieves need not be closed sieves. %
Authors usually write $\overline{\sS}$ for the closure, following the usual convention in topology for the same terminology (for a different notion). We prefer to write $\overline{\sS}$ for the element of $\sC$ that $\sS$ covers, since this emphasizes its role in the subobject classifier: $\sS : \overline{\sS} \to \Omega$, even though this unfortunately prevents us from using the commonplace notation for the closure (in the sense of sieves). In the ordinary localic case in which $\sC$ is the category of elements $U$ of a locale, and covers are epimorphic families, the closed sieves are precisely the principal lower sets $\downarrow U$. 
These are precisely those lower sets which are closed under joins of their elements. This is also a degenerate case in the sense that these closed sieves are (covering, hence) maximal. 

\dfn\label{d.matchingfamily}
Let $\sC$ be a generalized category and let $J$ be a topology on $\sC$. Let $P$ be a generalized presheaf on $\sC$, and let $\sS$ be a covering sieve in $\sC$. %
A {\em matching family in $P$ with respect to $\sS$} is %
a mapping $x:\sS \to \Sect(P),$ where $\Sect(P)$ is the set of all sections of $P$, that satisfies:
\enu
	\item $x$ assigns a section of $P$ to each element $\rho$ in $\sS$,
	\item (compatibility) For all $\sigma$ in $\sC$, $x(\rho) \cdot \sigma \downarrow$ implies
		$$x(\rho) \cdot \sigma = x(\rho \cdot \sigma).$$ %
\Enu
If $x = \set{x_\rho}$ is a matching family, an {\em amalgamation} of $x$ is an element $s \in \overline{\sS}$ such that
such that 
$$ \text{For every $\rho \in \sS$, } x(\rho) = s \cdot \rho.$$ %
that is, $x$ is completely determined by $s$. %
An amalgamation need not exist, or be uniquely defined if it exists. %
A {\em generalized sheaf} on $\sC$ or on $(\sC, J)$ is a generalized presheaf $P$ such that for every element $U$ of $\sC$, and for every cover $\sS$ of $U$, every matching family for $\sS$ has a unique amalgamation. %

The {\em generalized category of generalized sheaves on a site $\sC$} 
is the subcategory $\Sheaf(\sC,J)$ of the generalized category of generalized presheaves consisting of generalized presheaves $P$ that satisfy: 
$$P \text{ is a sheaf with respect to $J$}, \eand \source(P), \target(P) \text{ is in $\Sheaf(\sC,J)$.}$$
along with composition induced from the category of generalized presheaves. 
\Dfn

\thm\label{t.gensh}
The generalized category of generalized sheaves over a site $(\sC, J)$ is an ideal cartesian closed category, and has a subobject classifier. 
\Thm
\prf
Let $P, Q$ be sheaves; we claim that $P \vdash Q$ is also a sheaf, which is all we need in order to show that $\Sh(\sC,J)$ is an ideal cartesian closed category in light of Theorem \ref{t.pshiccc}. 
But indeed, we can show that $P \vdash Q$ is separated, and that 
$P \vdash Q$ amalgamates, much as in the one-categorical setting, see for example \cite{MaMo1}. %
Let $\S_\Omega$ be the target of a pointed profunctor $\Omega$ over $\sC$, and define sections by setting
$$\Fib_\Omega (U) := \set{\sS \mid \sS \text{ is a closed sieve at $U$}}$$
and complete the categorical structure by setting
$$\sS \cdot \rho := \stem{\sS}{\rho}$$
for $\rho: U \to V$ and $\sS$ a closed sieve at $U$. %
We incorporate $\Omega$ into the system of presheaves by making $\Omega$ an object: 
$$\hat\Omega = \bar\Omega = \Omega.$$
The terminal element $\trut$ is the subpointed profunctor whose fibers consist only of maximal sieves. %
An arrow $\arsoc: \trut \to \Omega$ in $\Sh(\sC, J)$ is given by inclusion. %
As in the one-categorical case, we can now show that $\Obsoc$ is a sheaf, and that $(\Obsoc, \arsoc)$ is the desired subobject classifier of $\Sh(\sC, J)$. 
\Prf


\cor\label{t.genpsh}
Every generalized category of generalized presheaves is an ideal elementary topos. 
\Cor
\prf
We have finite limits and colimits in this case, exactly as in the one-categorical case. 
\Prf




The following result extending \cite{LaTi1} says that in the setting of generalized categories, Lawvere-Tierney topologies on ideal elementary toposes extend the notion of Grothendieck topology to settings in which there is no site of definition. 

\thm\label{t.latiT}
Let $\sC$ be a generalized category, and let $\sE = \PSh(\sC)$ be the generalized category of generalized presheaves over $\sC$. 
For every Grothendieck topology $J$ on $\sC$, there is a Lawvere-Tierney topology $j$ on $\sE$ such that $P \in \sE$ is a generalized sheaf with respect to $J$ if and only if $P$ is a generalized sheaf with respect to $j$. Conversely, given a Lawvere-Tierney topology $j$ on a generalized category of generalized sheaves, there is a Grothendieck topology $J$ on $\sE$ such that the same property holds. 
\Thm
\prf
Essentially no changes are needed to the proof in the generalized setting. %
First, let $(\sC, J)$ be a site and let $\sE$ be as in the statement of \ref{t.latiT}. %
Define %
$$j = \set{g \mid \stem{\sS}{g} \in J}.$$
We can check from the definitions that $j$ is a Lawvere-Tierney topology. 
Let $j$ be a Lawvere-Tierney topology on an ideal elementary topos $\sE$, and assume that $\sE$ is $\PSh(\sC)$ for some base generalized category $\sC$. %
In this case, $\Omega$ is a particular pointed profunctor, and we can consider subpointed profunctors of it. 
Consider the subpointed profunctor of $\Omega$ given by 
$$J = \set{\text{sections } \sS \text{ of } \Omega \mid j(\sS) \text{ is the maximal sieve at $\overline{\sS}$} }.$$
We can check from the definitions that $J$ is a Grothendieck topology. 
The procedures yielding $j$ from $J$ and vice versa are inverse to one another. Moreover, a presheaf $P$ is a sheaf with respect to $J$ if and only if it is a sheaf with respect to $j$. 
\Prf

\pagebreak
\singlespacing
\addtocontents{toc}{\vspace{12pt}}
\addcontentsline{toc}{chapter}{\hspace{-1.6em} {\bf \large References }}
\bibliographystyle{abbrv}
\bibliography{mathmain}

\begin{thebibliography}{10}

\bibitem{AcAdMiVe1}
P.~Aczel, J.~Adamek, S.~Milius, and J.~Velebil.
\newblock Infinite trees and completely iterative theories: a coalgebraic view.
\newblock {\em Theoretical Computer Science}, 300:1--45, 2003.

\bibitem{sga4}
M.~Artin, A.~Grothendieck, and J.-L. Verdier.
\newblock {\em S{\'e}minaire de G{\'e}ometrie Alg{\'e}brique du Bois-Marie
  1963-1964 (SGA4)}, volume 270 of {\em Lecture Notes in Mathematics}.
\newblock Springer-Verlag, 1972.

\bibitem{Barendregt1}
H.~Barendregt.
\newblock Introduction to generalized type systems.
\newblock {\em Journal of Functional Programming}, 1(2):125--154, 1991.

\bibitem{BaWe1}
M.~Barr and C.~Wells.
\newblock {\em Category Theory for Computing Science}.
\newblock Prentice Hall, 1st edition, 1990.

\bibitem{BeU2}
J.~B{\'e}nabou.
\newblock Distributors at work.
\newblock based on notes by Thomas Streicher, 2000.

\bibitem{CockettConstellations}
R.~A.~G. Cockett.
\newblock Constellations.
\newblock private communication, 2016.

\bibitem{deBruijn1}
N.~G. de~{B}ruijn.
\newblock Telescope mappings in typed lambda calculus.
\newblock {\em Information and Computation}, 91(2):189--204, 1991.

\bibitem{Ehresmann1}
C.~Ehresmann.
\newblock {\em Cat{\'e}gories et structures}.
\newblock Paris: Dunod, 1965.

\bibitem{EiMa1}
S.~Eilenberg and S.~MacLane.
\newblock General theory of natural equivalences.
\newblock {\em Transactions of the American Mathematical Society}, 58:231--294,
  1945.

\bibitem{ElBlTi1}
C.~Elgot, S.~Bloom, and R.~Tindell.
\newblock On the algebraic structure of rooted trees.
\newblock {\em Journal of Computer and System Sciences}, 16:362--399, 1978.

\bibitem{FrD1}
P.~Freyd.
\newblock Aspects of topoi.
\newblock {\em Bulletin of the Australian Mathematical Society}, 7:1--76, 1972.

\bibitem{GrZ+}
G.~Gierz, K.~Hofmann, J.~D. Lawson, M.~Mislove, and D.~S. Scott.
\newblock {\em Continuous Lattices and Domains}.
\newblock Number~93 in Encyclopedia of Mathematics and its Applications.
  Cambridge University Press, 2003.

\bibitem{Howard1969}
W.~A. Howard.
\newblock {\em To {H}. {B}. {C}urry: Essays on Combinatory Logic, Lambda
  Calculus, and Formalism}, chapter The formulas-as-types notion of
  construction, pages 479--490.
\newblock Academic Press, 1980 (Reprint of 1969 article).

\bibitem{JoN1}
P.~T. Johnstone.
\newblock {\em Sketches of an Elephant (2 Volumes)}.
\newblock Clarendon Press, 2002.

\bibitem{KaSc1}
M.~Kashiwara and P.~Schapira.
\newblock {\em Categories and Sheaves}.
\newblock Springer, 2006.

\bibitem{Kobayashi1}
S.~Kobayashi.
\newblock Monad as modality.
\newblock {\em Theoretical Computer Science}, 175:29--74, 1997.

\bibitem{Kolmogorov1}
A.~Kolmogorov.
\newblock Zur {D}eutung der intuitionistischen {L}ogik.
\newblock {\em Mathematische Zeitschrift}, 35:58--65, 1932.

\bibitem{LaK1c}
J.~Lambek.
\newblock Deductive systems and categories {III}.
\newblock In {\em Lecture Notes in Mathematics}, volume 274, pages 57--82.
  Springer, 1972.

\bibitem{LaK2}
J.~Lambek.
\newblock Functional completeness of cartesian categories.
\newblock {\em Annals of Mathematical Logic}, 6:259--292, 1974.

\bibitem{LaSc1}
J.~Lambek and P.~J. Scott.
\newblock {\em Introduction to Higher Order Categorical Logic}.
\newblock Cambridge Studies in Advanced Mathematics. Cambridge University
  Press, 1986.

\bibitem{LaTi1}
F.~W. Lawvere.
\newblock Quantifiers as sheaves.
\newblock {\em Actes du {C}ongr{\'e}s {I}nternational des {M}ath{\'e}maticiens
  ({N}ice, 1970)}, 1:329--334, 1971.

\bibitem{MacCW}
S.~MacLane.
\newblock {\em Categories for the Working Mathematician}.
\newblock Springer, 2nd edition, 1998.

\bibitem{MaMo1}
S.~MacLane and I.~Moerdijk.
\newblock {\em Sheaves in Geometry and Logic}.
\newblock Springer, 1992.

\bibitem{MoI1}
E.~Moggi.
\newblock Notions of computation and monads.
\newblock {\em Information and Computation}, 93:55--92, 1991.

\bibitem{MuY2}
P.~Mulry.
\newblock {\em Applications of Categories in Computer Science}, volume 177 of
  {\em London Mathematical Society Lecture Note Series}, chapter Strong Monads,
  Algebras, and Fixed Points.
\newblock Cambridge University Press, 1992.

\bibitem{MuY1}
P.~Mulry.
\newblock Lifting theorems for {K}leisli categories.
\newblock In {\em Mathematical Foundations of Programming Semantics}, volume
  802 of {\em Lecture Notes in Computer Science}, pages 304--319. Springer,
  2005.

\bibitem{Pierce1}
B.~Pierce.
\newblock {\em Types and Programming Languages}.
\newblock MIT Press, 2002.

\bibitem{Reynolds1980}
J.~C. Reynolds.
\newblock {\em Proceedings of the Aarhus Workshop on Semantics-Directed
  Compiler Generation}, volume~94 of {\em Lecture Notes in Computer Science},
  chapter Using category theory to design implicit conversions and generic
  operators.
\newblock Springer-Verlag, 1980.

\bibitem{ScM3}
L.~T. Schoenbaum.
\newblock A generalization of the {C}urry-{H}oward correspondence.
\newblock arXiv:1612.02816, 2016.

\bibitem{ScMi}
L.~T. Schoenbaum.
\newblock Introduction to generalized categories.
\newblock arXiv:1612.02885, 2016.

\bibitem{SemmelrothSabry1}
M.~Semmelroth and A.~Sabry.
\newblock Monadic encapsulation in {ML}.
\newblock {\em {ACM} {SIGPLAN} International Conference on Functional
  Programming, Paris, France}, pages 8--17, 1999.

\bibitem{SmPl1}
M.~B. Smyth and G.~D. Plotkin.
\newblock The category-theoretic solution of recursive domain equations.
\newblock {\em SIAM Journal of Computing}, 11(4):761--783, 1982.

\bibitem{SpK3}
D.~I. Spivak.
\newblock Simplicial databases.
\newblock {\em arXiv:0904.2012}, 2009.

\bibitem{StT1}
R.~Street.
\newblock The geometry of oriented simplexes.
\newblock {\em Journal of Pure and Applied Algebra}, 49(3):283--335, 1987.

\bibitem{TrSc1}
A.~S. Troelstra and H.~Schwichtenberg.
\newblock {\em Basic Proof Theory}, volume~43 of {\em Cambridge Tracts in
  Theoretical Computer Science}.
\newblock Cambridge University Press, 2nd edition, 2000 (1st ed. 1996).

\bibitem{WaR1}
P.~Wadler.
\newblock Comprehending monads.
\newblock {\em Mathematical Structures in Computer Science}, 2:461--493, 1992.

\bibitem{WaR2}
P.~Wadler.
\newblock The marriage of effects and monads.
\newblock {\em {ACM} Transactions on Computational Logic}, 4(1):1--32, 2003.

\bibitem{WaR5}
P.~Wadler.
\newblock Propositions as types.
\newblock {\em Communications of the {ACM}}, 58(12):75--84, 2015.

\bibitem{WaD1}
M.~Wand.
\newblock Fixed-point constructions in order-enriched categories.
\newblock {\em Theoretical Computer Science}, 8:13--30, 1979.

\end{thebibliography}


\chapter*{Vita}
\doublespacing
\setlength{\parindent}{1.75em}
\vspace{0.2em}
\vspace{2ex}
\addtocontents{toc}{\vspace{12pt}}
\addcontentsline{toc}{chapter}{\hspace{-1.5em} {\bf \large Vita}}
Lucius Schoenbaum was born in New Orleans, Louisiana. He finished his undergraduate degree at the University of Georgia in 2005 with the major of Mathematics, and finished his Master's degree in Philosophy, also at the University of Georgia, under the supervision of O.B. Bassler in 2008, with a thesis in philosophical logic entitled ``Philosophy, Mathematics, and Proof Theory'' containing an axiomatization of typed set theory and a critique and attempted reconciliation of the contentious debate over foundations between Hilbert and Brouwer. Wishing to improve his understanding of the philosophical issues he wrestled with, he came to Louisiana State University to  continue his studies in mathematics in August 2009. He earned a Master's degree in Mathematics studying harmonic analysis under the supervision of Gestur Olafsson in 2012. He is currently a candidate for the degree of Doctor of Philosophy in Mathematics, which will be awarded in Spring 2017. His dissertation research has been carried out under the supervision of Daniel S. Sage and Jimmie D. Lawson. 

\end{document}